\documentclass[11pt,twoside,leqno,a4paper]{article}

\usepackage{geometry, amsfonts, amsthm, amsmath, indentfirst, graphicx, enumerate,xcolor, subcaption, lineno, cite, float}
\usepackage{paralist}
\usepackage{amssymb}

\usepackage{extarrows}% for making wide = symbols. See the command \xlongequal

\usepackage[nottoc,notlot,notlof]{tocbibind}%for appearing references in the list of contents

\RequirePackage[breaklinks,colorlinks]{hyperref}
\hypersetup{linkcolor=blue,citecolor=blue,urlcolor=blue}

\usepackage{titlesec}
\titleformat{\section}{\filcenter\Large\bfseries}{\normalfont\thesection.}{0.5em}{}
\titleformat{\subsection}{\filcenter\bfseries}{\normalfont\thesubsection.}{.5em}{}

\newtheorem{mythm}{\,\,\,\,\,\, \normalfont\scshape Theorem}[section]
\newtheorem{mycor}[mythm]{\,\,\,\,\,\, \normalfont\scshape Corollary}
\newtheorem{mylem}[mythm]{\,\,\,\,\,\, \normalfont\scshape Lemma}
\newtheorem{myprop}[mythm]{\,\,\,\,\,\, \normalfont\scshape Proposition}
\newtheorem{myconjecture}[mythm]{\,\,\,\,\,\, \normalfont\scshape Conjecture}
\newtheorem{mydefn}[mythm]{\,\,\,\,\,\, \normalfont\it Definition}
\newtheorem{myrem}[mythm]{\,\,\,\,\,\, \normalfont\it Remark}
\newtheorem{mynotation}[mythm]{\,\,\,\,\,\, \normalfont\it Notation}

\newtheorem{myassumption}{\,\,\,\,\,\, \normalfont\it Assumption}

\newtheorem{theoremA}{\,\,\,\,\,\, \normalfont\scshape Theorem}

\newtheorem{theoremB}{\,\,\,\,\,\, \normalfont\scshape Theorem}

\newtheorem{theoremC}{\,\,\,\,\,\, \normalfont\scshape Theorem}

\definecolor{mygreen}{rgb}{0.50,0,0.50}

\numberwithin{equation}{section}

\providecommand{\keywords}[1]{\small	Keywords: #1}

\newcommand\blfootnote[1]{%
\begingroup
\renewcommand\thefootnote{}\footnote{#1}%
\addtocounter{footnote}{-1}%
\endgroup
}

\usepackage{etoolbox}
\makeatletter
\patchcmd{\endmyassumption}{\@endpefalse}{}{}{}
\makeatother

\usepackage{authblk}

\AtEndDocument{%
\par
\medskip
\begin{tabular}{@{}l@{}}%
\textsc{\small Department of Mathematics, Imperial College London, London SW7 2AZ, United Kingdom}\\
\textsc{\small Faculty of Engineering and Natural Sciences, Kadir Has University, 34083 Istanbul, Turkey}\\
\textit{\small E-mail address}: \href{mailto:s.bakrani-balani@imperial.ac.uk}{\texttt{\small s.bakrani-balani@imperial.ac.uk}}
\end{tabular}}

\AtEndDocument{%
\par
\medskip
\begin{tabular}{@{}l@{}}%
\textsc{\small Department of Mathematics, Imperial College London, London SW7 2AZ, United Kingdom}\\
\textit{\small E-mail address}: \href{mailto:jeroen.lamb@imperial.ac.uk}{\texttt{\small jeroen.lamb@imperial.ac.uk}}
\end{tabular}}

\AtEndDocument{%
\par
\medskip
\begin{tabular}{@{}l@{}}%
\textsc{\small Department of Mathematics, Imperial College London, London SW7 2AZ, United Kingdom}\\
\textit{\small E-mail address}: \href{mailto:d.turaev@imperial.ac.uk}{\texttt{\small d.turaev@imperial.ac.uk}}
\end{tabular}}

\begin{document}

\title{\bf Invariant manifolds of homoclinic orbits and the dynamical consequences of a super-homoclinic:
A case study in \(\mathbb{R}^4\) with \(\mathbb{Z}_2\)-symmetry and integral of motion}
\author{\textsc{Sajjad Bakrani}, \textsc{Jeroen S. W. Lamb} and \textsc{Dmitry Turaev}}
\date{August 8, 2022}
\maketitle

\begin{abstract}
We consider a \(\mathbb{Z}_2\)-equivariant flow in \(\mathbb{R}^{4}\) with an integral of motion and a hyperbolic equilibrium with a transverse homoclinic orbit \(\Gamma\). We provide criteria for the existence of stable and unstable invariant manifolds of \(\Gamma\). We prove that if these manifolds intersect transversely, creating a so-called super-homoclinic, then in any neighborhood of this super-homoclinic there exist infinitely many multi-pulse homoclinic loops. An application to a system of coupled nonlinear Schr\"odinger equations is considered.
\end{abstract}

\blfootnote{\keywords{homoclinic, super-homoclinic, invariant manifold, coupled Schr\"odinger equations}}
\blfootnote{SB was supported by EU Marie Sklodowska-Curie ITN Critical Transitions in Complex Systems (H2020-MSCA-ITN-2014 643073 CRITICS), European Union ERC Advanced Grant of Sebastian Van Strien (339523 RGDD), and TUBITAK grant (No. 118C236). JSWL thanks the London Mathematical Laboratory, for support through its fellowship programme. DT was supported by Leverhulme Trust (RPG-2021-072), by the grants 19-11-00280 and 19-71-10048 of the Russian Science Foundation (RSF), by the Mathematical Center at the Lobachevsky University of Nizhny Novgorod, and by the grant 075-15-2019-1931 of Russian
Ministry of Science and Higher Education.}

\tableofcontents

%\listoffigures

%\addtocontents{toc}{\protect\setcounter{tocdepth}{2}}
% This command control what appears on the table of content.

\section{Introduction}

\subsection{Background}

Consider a Hamiltonian system (or more generally, a system with a smooth first integral) defined for \(x\in \mathbb{R}^{2n}\), \(n\geq 2\) with an integral \(H\), and a hyperbolic equilibrium \(O\) at the origin. An orbit \(\Gamma = \{x\left(t\right): t\in \mathbb{R}\}\) of this system is said to be '{\it homoclinic to \(O\)}' or '{\it a homoclinic loop}' if it belongs to both stable and unstable invariant manifolds of \(O\), or equivalently, \(x\left(t\right) \rightarrow O\) as \(t\rightarrow \pm \infty\). Existence of homoclinic orbits for systems with a smooth first integral is known to be a robust phenomenon. This is due to the fact that the \(n\)-dimensional stable and unstable invariant manifolds of \(O\) lie in the same (\(2n-1\))-dimensional level \(H = \text{constant}\), and they may intersect transversely in that level along the homoclinic orbits. A natural question which arises here is the possible dynamics near homoclinic orbits in the level \(H\left(x\right) = H\left(O\right)\).

When \(x(t)\) oscillates as it converges to \(O\) (this happens when the leading, i.e. the nearest to the imaginary axis, eigenvalues of the linear part of the system at \(O\) are complex), the dynamics in \(H^{-1}\left(O\right)\) is highly non-trivial \cite{Devaney1976, BelyakovShilnikov1990, Lerman1991, Lerman2000, Lerman1997, Buffoni1996, Barrientos2016robust}. On the other hand, when the leading eigenvalues are real and \(\|x\left(t\right)\|\) decays to zero monotonically as \(t\rightarrow\pm\infty\), the generic dynamics in \(H^{-1}\left(O\right)\) are very simple. Thus, the only orbits staying in a small neighborhood of a finite bunch of generic homoclinic loops to a saddle with real leading eigenvalues are only the homoclinic loops themselves and the point \(O\) \cite{DimaShilnikov1989,Turaev2014}.

In this case, to have interesting behavior, we need degeneracies or symmetries. It was shown in \cite{ShilnikovDimaSuperhomoclinic1997} (based on an earlier work \cite{EleonskiiKulagin1989}) that symmetries can lead to the emergence of the so-called super-homoclinic orbits which, in turn, serve as limits of infinite series of multi-pulse
homoclinic loops. Namely, as shown in \cite{ShilnikovDimaSuperhomoclinic1997,DimaMultipulse2001}, in certain situations, the homoclinic loops or bunches of homoclinic loops can have stable or unstable invariant manifolds; the super-homoclinics correspond to the intersection of these manifolds.

In the non-conservative setting, super-homoclinic orbits and the non-trivial dynamics associated with them were discovered and studied by Homburg \cite{Homburg1996memoirs}. Eleonsky et al. \cite{EleonskiiKulagin1989} spotted super-homoclinic orbits in the numerical investigation of an electromagnetic field in a nonlinear medium. Barrientos et al. \cite{Rodrigues2019} found super-homoclinics near a homoclinic loop to a saddle-focus in the context of reversible systems, similar to the structure described for the Hamiltonian case by Belyakov and Shilnikov in \cite{BelyakovShilnikov1990}. Chawanya and Ashwin \cite{Ashwin2010} built an example of a heteroclinic network that possesses a super-homoclinic in the sense of an orbit which connects sub-networks. In general, super-homoclinic orbits may potentially appear in heteroclinic networks, especially if the network undergoes a chaotic behavior, see e.g. \cite{NajafiDuartePeixe2020}.

In this paper, we consider the simplest case of \(\mathbb{Z}_{2}\)-symmetry which results in the emergence of stable and unstable invariant manifolds for homoclinic loops, enabling the creation of super-homoclinic orbits, and describe the multi-pulse homoclinics associated to them.

\subsection{Problem setting and results}\label{Setting and results}

Consider a \(\mathcal{C}^{\infty}\)-smooth 4-dimensional system of differential equations
\begin{equation}\label{eq100}
\dot{x}=X(x), \quad x\in \mathbb{R}^{4},
\end{equation}
with a \(\mathcal{C}^{\infty}\)-smooth first integral \(H: \mathbb{R}^{4}\rightarrow \mathbb{R}\), i.e.
\begin{equation}\label{eq200}
\nabla H(x) \cdot X(x) \equiv 0.
\end{equation}
\begin{myassumption}\label{assumption10}
\(X\) has a hyperbolic equilibrium state \(O\) at the origin.
\end{myassumption}
By (\ref{eq200}), we have \(H^{\prime}(0) X^{\prime}(0) \equiv 0\).
Since \(X^{\prime}(0)\) is nonsingular by Assumption \ref{assumption10}, the linear part of \(H\) at \(O\) vanishes.
\begin{myassumption}\label{assumption20}
The quadratic part of \(H\) at \(O\) is a nondegenerate quadratic form.
\end{myassumption}
It is easy to see (see e.g. \cite{BakraniPhDthesis}) that when Assumptions \ref{assumption10} and \ref{assumption20} are satisfied, system (\ref{eq100}) near \(O\) can be brought to the following form by a linear transformation:
\begin{equation}\label{eq300}
\dot{u}=-Au + o\left(\lvert u\rvert, \lvert v\rvert\right), \quad \dot{v}=A^{T}v + o\left(\lvert u\rvert, \lvert v\rvert\right),
\end{equation}
where \(u=\left(u_{1},u_{2}\right)\in \mathbb{R}^{2}\), \(v=\left(v_{1},v_{2}\right) \in \mathbb{R}^{2}\) and \(A\) is a matrix whose eigenvalues have positive real parts. Moreover, the first integral takes the form:
\begin{equation}\label{eq400}
H=\langle v, Au\rangle + o\left(u^{2} + v^{2}\right),
\end{equation}
where \(\langle\cdot, \cdot\rangle\) is the standard inner product on \(\mathbb{R}^{2}\).
\begin{myassumption}\label{assumption30}
System (\ref{eq300}) is invariant with respect to the symmetry
\begin{equation}\label{eq425}
(u_{1}, v_{1})\leftrightarrow(-u_{1}, -v_{1}).
\end{equation}
\end{myassumption}
Assumption \ref{assumption30} implies that the plane \(\lbrace u_{1}=v_{1}=0 \rbrace\) is invariant with respect to the flow of system (\ref{eq300}). 

Note that we can always assume that \(H\) is invariant with respect to symmetry (\ref{eq425}), i.e.
\begin{equation}\label{eq426}
H\left(-u_{1}, u_{2}, -v_{1}, v_{2}\right) = H\left(u_{1}, u_{2}, v_{1}, v_{2}\right).
\end{equation}
Otherwise, \(\widetilde{H}\left(u_{1}, u_{2}, v_{1}, v_{2}\right):= \frac{1}{2} \left[H\left(u_{1}, u_{2}, v_{1}, v_{2}\right) + H\left(-u_{1}, u_{2}, -v_{1}, v_{2}\right)\right]\) can be taken as the first integral.

\begin{figure}
\centering
\includegraphics[scale=.15]{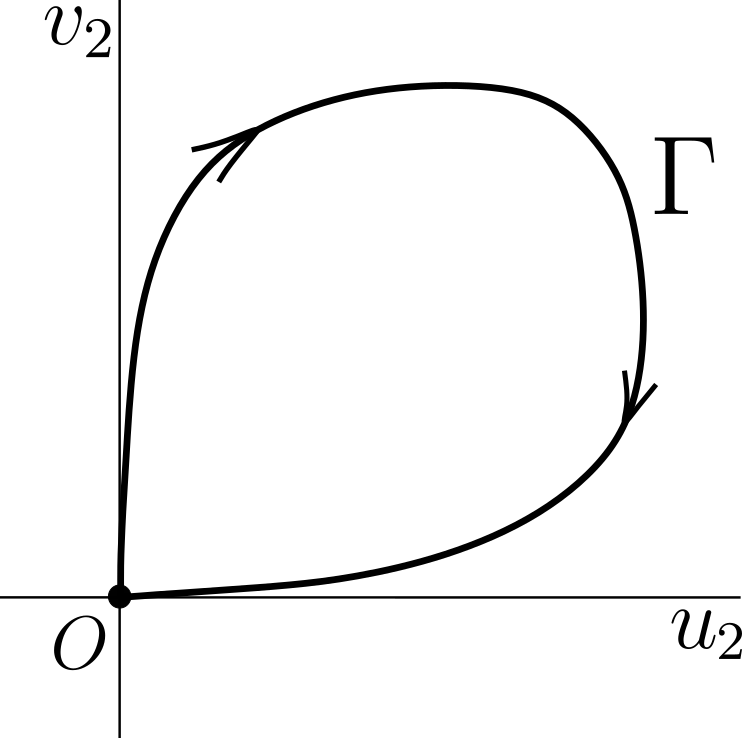}
\caption{The transverse homoclinic loop \(\Gamma\) in the invariant plane \(\lbrace u_{1}=v_{1}=0 \rbrace\).}
\label{Figure8o98b87rv8i6}
\end{figure}

The equilibrium state \(O\) is a saddle with 2-dimensional stable and unstable invariant manifolds \(W^{s}\left(O\right)\) and \(W^{u}\left(O\right)\) which are tangent at \(O\) to the \(u\)-plane and \(v\)-plane respectively. Both the invariant manifolds lie in the 3-dimensional level \(\lbrace H = 0 \rbrace\) and may intersect transversely in that level, producing a number of homoclinic loops. We consider the following specific case:
\begin{myassumption}\label{assumption50}
In the invariant plane \(\lbrace u_{1}=v_{1}=0 \rbrace\), there exists a homoclinic loop \(\Gamma\) of the transverse intersection of \(W^{s}\left(O\right)\) and \(W^{u}\left(O\right)\) (see Figure \ref{Figure8o98b87rv8i6}).
\end{myassumption}

Since the action of this symmetry commutes with the linear part of system (\ref{eq300}), the matrix \(A\) is diagonal and takes the form
\begin{equation*}
A = \left(
\begin{array}{ccc}
\lambda_{1} & 0 \\
0 & \lambda_{2}
\end{array} \right),
\end{equation*}
for some positive real numbers \(\lambda_{1}\) and \(\lambda_{2}\). Without loss of generality, let \(\lambda_{1} \leq \lambda_{2}\).
\begin{myassumption}\label{assumption40}
\(\lambda_{2} \neq 2\lambda_{1}\).
\end{myassumption}
This is not a technical assumption. Indeed, we will see that the cases \(\lambda_{2}<2\lambda_{1}\) and \(2\lambda_{1} < \lambda_{2}\) are dynamically different.

Let \(\mathcal{U}\) be a sufficiently small neighborhood of \(\Gamma \cup \{O\}\) in the zero level-set \(\lbrace H=0\rbrace\). The main issue which is addressed in this paper is giving a complete description of dynamics in \(\mathcal{U}\). 

\begin{mydefn}\label{Defno9uo87ictcuyghxfhgfgdr}
Let \(\mathcal{A} = \Gamma_{1} \cup \Gamma_{2}\cup \cdots\cup \Gamma_{m}\), where \(\Gamma_{i}\) are homoclinic to an equilibrium \(O\), and \(m\geq 1\). Consider a sufficiently small open neighborhood \(\mathcal{U}\) of \(\mathcal{A} \cup \{O\}\) in the energy level of \(O\). The local stable (resp. unstable) set of \(\mathcal{A}\), denoted by \(W^{s}_{\text{loc}}(\mathcal{A}, \mathcal{U})\) (resp. \(W^{u}_{\text{loc}}(\mathcal{A}, \mathcal{U})\)), is the union of \(\mathcal{A}\) itself and the set of the points in \(\mathcal{U}\) whose forward (resp. backward) orbits lie in \(\mathcal{U}\) and their \(\omega\)-limit sets (resp. \(\alpha\)-limit sets) coincide with \(\mathcal{A} \cup \lbrace O\rbrace\). We may use the notations \(W^{s}_{\text{loc}}(\mathcal{A})\) and \(W^{u}_{\text{loc}}(\mathcal{A})\) for the stable and unstable sets of \(\mathcal{A}\) when the neighborhood \(\mathcal{U}\) is clear from the context.
\end{mydefn}

By this definition, the local stable and unstable sets of \(\Gamma\) always contain \(\Gamma\). Note that these sets lie in the zero-level set \(\lbrace H=0 \rbrace\). Denote by \(W^{s}_{\mathcal{U}}\left(O\right)\) (resp. \(W^{u}_{\mathcal{U}}\left(O\right)\)) the set of the points in \(W^{s}_{\text{glo}}(O)\) (resp. \(W^{u}_{\text{glo}}(O)\)) whose forward (resp. backward) orbits lie entirely in \(\mathcal{U}\). Obviously,
\begin{equation*}
W^{s}_{\mathcal{U}}\left(O\right) \cap W^{s}_{\text{loc}}\left(\Gamma\right) = W^{u}_{\mathcal{U}}\left(O\right) \cap W^{u}_{\text{loc}}\left(\Gamma\right) = \Gamma.
\end{equation*}

\subsubsection{Dynamics near a single homoclinic orbit}

Our first result is the following:
\begin{theoremA}\label{Th782ju2iu22}
Under Assumptions \ref{assumption10}-\ref{assumption40}, the forward (backward) orbit of a point in \(\mathcal{U}\) lies entirely in \(\mathcal{U}\) if and only if it belongs to \(W^{s}_{\mathcal{U}}\left(O\right) \cup W^{s}_{\text{loc}}\left(\Gamma\right)\) (resp. \(W^{u}_{\mathcal{U}}\left(O\right) \cup W^{u}_{\text{loc}}\left(\Gamma\right)\)).
\end{theoremA}
By this theorem, to understand the dynamics near the homoclinic orbit \(\Gamma\), we need to study the local stable and unstable sets of this orbit. This is done in the following two theorems.
\begin{theoremA}\label{thmkuyvuyuy090988}
If \(\lambda_{2} < 2\lambda_{1}\) and Assumptions \ref{assumption10}-\ref{assumption50} hold, then \(W^{s}_{\text{loc}}\left(\Gamma\right) = W^{u}_{\text{loc}}\left(\Gamma\right) = \Gamma\). Thus, by Theorem \ref{Th782ju2iu22}, the forward and backward orbits of any point in \(\mathcal{U}\) either leave \(\mathcal{U}\) or converge to \(O\).
\end{theoremA}

\begin{figure}
\centering
\includegraphics[scale=0.20]{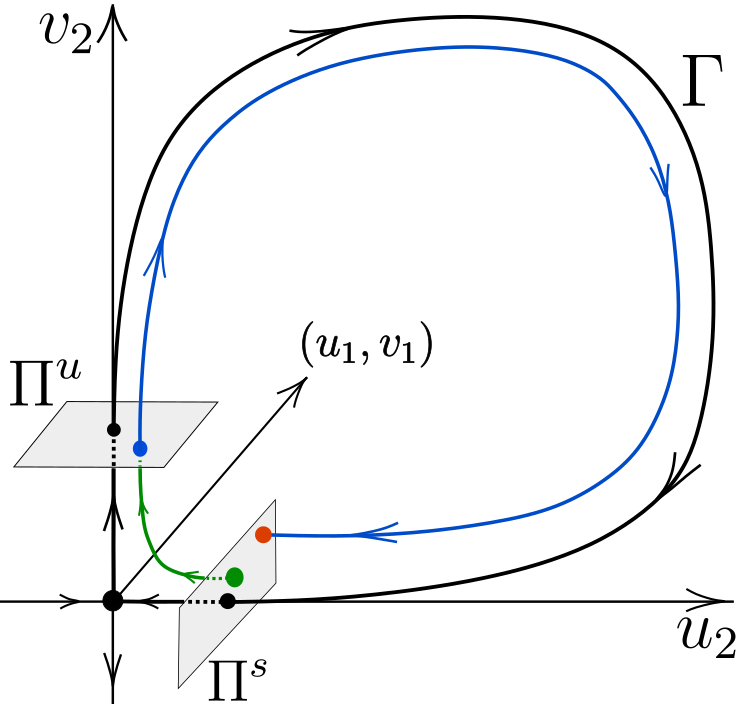}
\caption{\small This figure shows the positions of the cross-sections \(\Pi^{s}\) and \(\Pi^{u}\) to the homoclinic loop \(\Gamma\). The green and blue curves correspond to the maps \(T^{\text{loc}}\) and \(T^{\text{glo}}\), respectively. Namely, \(T^{\text{loc}}\) maps the green point on \(\Pi^{s}\) to the blue point on \(\Pi^{u}\) and then \(T^{\text{glo}}\) maps the blue point to the red point on \(\Pi^{s}\). The red point is the image of the green point by the Poincar\'e map \(T = T^{\text{glo}}\circ T^{\text{loc}}\).}
\label{Figure00m9n8nb4y2bc}
\end{figure}

The next theorem describes the local stable and local unstable sets of \(\Gamma\) when \(2\lambda_{1} < \lambda_{2}\). The formulation uses a specific choice of coordinates near the equilibrium \(O\). We introduce this coordinate system in Section \ref{localdynamicsnearO} (see normal form (\ref{eq23000})). For this choice of coordinates, system (\ref{eq300}) keeps its form and its invariance with respect to symmetry (\ref{eq425}). Moreover, the first integral takes the form
\begin{equation}\label{eq430}
H\left(u_{1}, u_{2}, v_{1}, v_{2}\right)=\lambda_{1} u_{1} v_{1} - \lambda_{2} u_{2} v_{2} + o\left(u^{2} + v^{2}\right),
\end{equation}
and satisfies (\ref{eq426}). The local stable and unstable, as well as the local strong stable and strong unstable, invariant manifolds of \(O\) are straightened (i.e. \(W^{s}_{\text{loc}}\left(O\right) = \lbrace v_{1} = v_{2} = 0\rbrace\), \(W^{u}_{\text{loc}}\left(O\right) = \lbrace u_{1} = u_{2} = 0\rbrace\), \(W^{ss}_{\text{loc}}\left(O\right) = \lbrace u_{1} = v_{1} = v_{2} = 0\rbrace\), \(W^{uu}_{\text{loc}}\left(O\right) = \lbrace u_{1} = u_{2} = v_{1} = 0\rbrace\)), and the loop \(\Gamma\) leaves \(O\) along \(v_{2}\)-axis toward positive \(v_{2}\) and enters \(O\) along \(u_{2}\)-axis toward positive \(u_{2}\) (see Figure \ref{Figure00m9n8nb4y2bc}).

Take a small \(\delta > 0\) and consider two small 2-dimensional cross-sections to the loop \(\Gamma\) inside the level \(\lbrace H=0 \rbrace\): \(\Pi^{s}=\lbrace u_{2}=\delta\rbrace \cap \lbrace H=0\rbrace\) and \(\Pi^{u} =\lbrace v_{2}=\delta\rbrace \cap \lbrace H=0\rbrace\) (see Figure \ref{Figure00m9n8nb4y2bc}). On each of the cross-sections \(\Pi^{s}\) and \(\Pi^{u}\), the variables \(u_{2}\) and \(v_{2}\) are uniquely determined by \(\left(u_{1}, v_{1}\right)\) (see Lemma \ref{Cor518b8b89or48v10010}). This allows us to parametrize \(\Pi^{s}\) and \(\Pi^{u}\) by the \(\left(u_{1}, v_{1}\right)\)-coordinates.

Orbits which lie in \(\mathcal{U}\) define a Poincar\'e map \(T\) that takes a subset of \(\Pi^{s}\) to \(\Pi^{s}\). This map can be written as a composition of a local map \(T^{\text{loc}}\) from a subset of \(\Pi^{s}\) to \(\Pi^{u}\) which is defined by the flow inside the \(\delta\)-neighborhood of \(O\), and a global map \(T^{\text{glo}}\) from \(\Pi^{u}\) to \(\Pi^{s}\) which is defined by the flow near the global piece of \(\Gamma\) outside the \(\delta\)-neighborhood of \(O\) (see Figure \ref{Figure00m9n8nb4y2bc}). Since the flight time from \(\Pi^{u}\) to \(\Pi^{s}\) is bounded, the global map \(T^{\text{glo}}\) is a diffeomorphism. Consider the points \(M^{s} = \Gamma \cap \Pi^{s}\) and \(M^{u} = \Gamma \cap \Pi^{u}\) (note that both points correspond to \(\left(0,0\right)\) in \(\Pi^{s,u}\)). The Taylor expansion of \(T^{\text{glo}}\) at \(M^{u}\) has the form
\begin{equation}\label{eq63000}
T^{\text{glo}}\left(u_{1}, v_{1}\right) = \left(a u_{1} + b v_{1} + o\left(u_{1}, v_{1}\right),\, c u_{1} + d v_{1} + o\left(u_{1}, v_{1}\right)\right),
\end{equation}
for some \(a, b, c, d \in \mathbb{R}\). Since the local unstable manifold of \(O\) corresponds to \(\{\mathrm{u}_{1}=0\}\), and the local stable manifold corresponds to \(\{\mathrm{v}_{1}=0\}\), the transversality assumption (see Assumption \ref{assumption50}) is equivalent to \(d\neq 0\).

\begin{theoremA}\label{Invariantmanifoldthm}
Let \(2\lambda_{1} < \lambda_{2}\) and Assumptions \ref{assumption10}-\ref{assumption50} hold. Suppose that system (\ref{eq300}) near the equilibrium \(O\) is brought to the form (\ref{eq23000}) and let \(b\), \(c\) and \(d\) in (\ref{eq63000}) be non-zero.
\begin{enumerate}[(i)]
\item If \(cd > 0\), then \(W^{s}_{\text{loc}}(\Gamma)= \Gamma\). If \(cd<0\), then \(W^{s}_{\text{loc}}\left(\Gamma\right)\) is a \(\mathcal{C}^1\)-smooth 2-dimensional invariant manifold which is tangent to \(W_{\text{glo}}^{s}\left(O\right)\) at every point of \(\Gamma\).
\item If \(bd<0\), then \(W^{u}_{\text{loc}}(\Gamma) = \Gamma\). If \(bd > 0\), then \(W^{u}_{\text{loc}}\left(\Gamma\right)\) is a \(\mathcal{C}^1\)-smooth 2-dimensional invariant manifold which is tangent to \(W_{\text{glo}}^{u}\left(O\right)\) at every point of \(\Gamma\).
\end{enumerate}
\end{theoremA}

To stress the smoothness of \(W^{s}_{\text{loc}}(\Gamma)\) and \(W^{u}_{\text{loc}}(\Gamma)\), we further call them local stable and unstable invariant manifolds.

We remark a parallel to the case of a general homoclinic loop, considered in detail by Homburg \cite{Homburg1996memoirs}, see also \cite{Turaev1984case,Sandstede2000center,Shashkov1999existence}, for which the existence of stable/unstable manifolds depends on the sign of a certain coefficient in the Poincar\'e map (the so-called separatrix value \cite{Dimabook}).

\subsubsection{Dynamics near a homoclinic figure-eight}

\begin{figure}
\centering
\includegraphics[scale=.18]{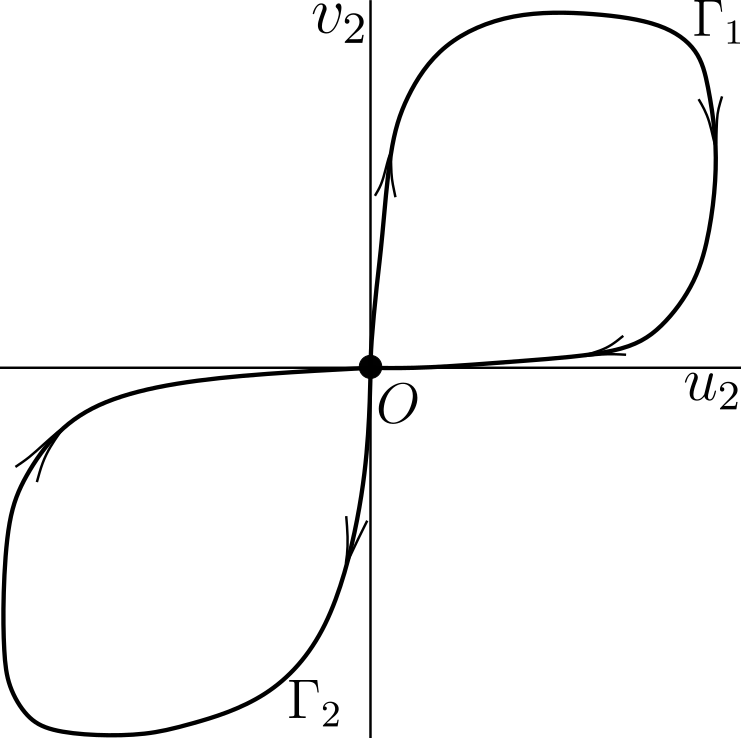}
\caption{\small A pair of transverse homoclinic loops (\(\Gamma_{1}\cup\Gamma_{2}\)) in the invariant plane \(\lbrace u_{1}=v_{1}=0 \rbrace\).}
\label{Figureyuybbilzqbe564}
\end{figure}

Next, we consider the existence of a pair of homoclinic loops in the invariant plane \(\lbrace u_{1}= v_{1} = 0\rbrace\):
\begin{myassumption}\label{assumption60}
There exist two homoclinic loops \(\Gamma_{1}\) and \(\Gamma_{2}\) of transverse intersection of \(W^{s}\left(O\right)\) and \(W^{u}\left(O\right)\) in the invariant plane \(\lbrace u_{1}=v_{1}=0 \rbrace\) such that they leave and enter \(O\) along opposite directions (see Figure \ref{Figureyuybbilzqbe564}).
\end{myassumption}

Such scenario happens generically, when the level-set \(\lbrace H = 0\rbrace\) is compact. Let \(\mathcal{V}\) be a small neighborhood of \(\Gamma_{1} \cup \{O\} \cup \Gamma_{2}\) in the level-set \(\lbrace H = 0\rbrace\) and denote by \(W^{s}_{\mathcal{V}}\left(O\right)\) (resp. \(W^{u}_{\mathcal{V}}\left(O\right)\)) the set of the points in \(W^{s}_{\text{glo}}(O)\) (resp. \(W^{u}_{\text{glo}}(O)\)) whose forward (resp. backward) orbits lie entirely in \(\mathcal{V}\). Then
\begin{theoremB}\label{Thmyipoinoib53}
Under Assumptions \ref{assumption10}-\ref{assumption30}, \ref{assumption40} and \ref{assumption60}, the forward (resp. backward) orbit of a point in \(\mathcal{V}\) lies entirely in \(\mathcal{V}\) if and only if it belongs to \(W^{s}_{\mathcal{V}}\left(O\right) \cup W^{s}_{\text{loc}}\left(\Gamma_{1}\right)\cup W^{s}_{\text{loc}}\left(\Gamma_{1} \cup \Gamma_{2}\right)\cup W^{s}_{\text{loc}}\left(\Gamma_{2}\right)\) (resp. \(W^{u}_{\mathcal{V}}\left(O\right) \cup W^{u}_{\text{loc}}\left(\Gamma_{1}\right)\cup W^{u}_{\text{loc}}\left(\Gamma_{1} \cup \Gamma_{2}\right)\cup W^{u}_{\text{loc}}\left(\Gamma_{2}\right)\)).
\end{theoremB}

The next two theorems give analogues of Theorems \ref{thmkuyvuyuy090988} and \ref{Invariantmanifoldthm} for the case of homoclinic figure-eight.
\begin{theoremB}\label{thm89bqyvrtyvtv}
If \(\lambda_{2} < 2\lambda_{1}\), and Assumptions \ref{assumption10}-\ref{assumption30} and \ref{assumption60} hold, then \(W^{s}_{\text{loc}}\left(\Gamma_{1} \cup \Gamma_{2}\right) = W^{u}_{\text{loc}}\left(\Gamma_{1} \cup\Gamma_{2}\right) = \Gamma_{1}\cup \Gamma_{2}\). Thus, by Theorem \ref{Thmyipoinoib53}, the forward and backward orbits of any point in \(\mathcal{V}\) either leave \(\mathcal{V}\) or converge to \(O\).
\end{theoremB}

Consider the cross-sections \(\Pi^{s}_{1} = \lbrace u_{2} = \delta\rbrace \cap \lbrace H=0\rbrace\) and \(\Pi^{u}_{1} = \lbrace v_{2} = \delta\rbrace \cap \lbrace H=0\rbrace\) on \(\Gamma_{1}\), and \(\Pi^{s}_{2} = \lbrace u_{2} = -\delta\rbrace \cap \lbrace H=0\rbrace\) and \(\Pi^{u}_{2} =\lbrace v_{2} = -\delta\rbrace \cap \lbrace H=0\rbrace\) on \(\Gamma_{2}\) (see Figure \ref{Figure673bob8i7vrq7cvraa}). We can choose \((u_{1}, v_{1})\)-coordinates on each of these cross-sections (see Lemma \ref{Cor518b8b89or48v10010}). Let \(T_{i}\), \(T^{\text{loc}}_{i}\) and \(T^{\text{glo}}_{i}\) be the associated maps along \(\Gamma_{i}\). For \(i=1,2\), consider the points \(M^{s,u}_{i} = \Gamma_{i} \cap \Pi^{s,u}_{i}\), and let \(a_{i}\), \(b_{i}\), \(c_{i}\) and \(d_{i}\) be the corresponding Taylor coefficients of $T^{\text{glo}}_{i}$, as in (\ref{eq63000}).%\textcolor{red}{???????}
\begin{theoremB}\label{Thm6892jjdibbea}
Assume \(2\lambda_{1} < \lambda_{2}\) and Assumptions \ref{assumption10}-\ref{assumption30} and \ref{assumption60}. Suppose that system (\ref{eq300}) near the equilibrium \(O\) is brought to the form (\ref{eq23000}) and let \(b_{i}\), \(c_{i}\) and \(d_{i}\) (\(i=1,2\)) be non-zero.
\begin{enumerate}[(i)]
\item If \(c_{1}d_{1} > 0\) and \(c_{2}d_{2} > 0\), then \(W^{s}_{\text{loc}}\left(\Gamma_{1} \cup \Gamma_{2}\right)\) is a \(\mathcal{C}^1\)-smooth 2-dimensional invariant manifold which is tangent to \(W_{\text{glo}}^{s}\left(O\right)\) at every point of \(\Gamma_{1} \cup \Gamma_{2}\).
\item If \(b_{1} d_{1} <0\) and \(b_{2}d_{2} <0\), then \(W^{u}_{\text{loc}}\left(\Gamma_{1} \cup \Gamma_{2}\right)\) is a \(\mathcal{C}^1\)-smooth 2-dimensional invariant manifold which is tangent to \(W_{\text{glo}}^{u}\left(O\right)\) at every point of \(\Gamma_{1} \cup \Gamma_{2}\).
\item Otherwise, we have \(W^{s}_{\text{loc}}\left(\Gamma_{1} \cup\Gamma_{2}\right) = W^{u}_{\text{loc}}\left(\Gamma_{1} \cup\Gamma_{2}\right) = \Gamma_{1} \cup \Gamma_{2}\).
\end{enumerate}
\end{theoremB}

Let \(0 < \gamma = \frac{\lambda_{1}}{\lambda_{2}}\leq 1\), and consider a one-parameter family \(\{X_{\gamma}\}\) of the vector fields of the form (\ref{eq300}) that satisfy the assumptions stated above. In particular, suppose that \(X_{\gamma}\) possesses a homoclinic orbit \(\Gamma_{\gamma}\) which persists as \(\gamma\) varies. Then, according to our results, when \(\gamma > 0.5\), there is no dynamics near the homoclinic orbit \(\Gamma_{\gamma}\) in its energy level, while when \(\gamma < 0.5\), depending on how the global map behaves (i.e. what the coefficients \(a\), \(b\), \(c\) and \(d\) are), there may exist stable and unstable invariant manifolds to the homoclinic loop \(\Gamma_{\gamma}\). This lets us to conjecture that saddle periodic orbits can be born in the level \(H=0\) as \(\gamma\) crosses \(\frac{1}{2}\). This question requires a further investigation.

\subsubsection{Dynamics near a super-homoclinic orbit}

Coming back to the case of the single homoclinic loop \(\Gamma\), we consider the case in which both \(W^{s}_{\text{loc}}\left(\Gamma\right)\) and \(W^{u}_{\text{loc}}\left(\Gamma\right)\) are non-trivial. Notice that, according to Theorem \ref{Invariantmanifoldthm}, in order for these two manifolds to coexist, we require \(cd < 0\) and \(bd > 0\). Continuing these two local manifolds by the flow of the system gives the global stable and unstable invariant manifolds of \(\Gamma\), denoted by \(W^{s}_{\text{glo}}\left(\Gamma\right)\) and \(W^{u}_{\text{glo}}\left(\Gamma\right)\), respectively. These manifolds lie in the 3-dimensional level \(\{H=0\}\) which means that it would be reasonable if we assume that they intersect transversely in that level. Any orbit at this intersection is bi-asymptotic, or in other words, homoclinic to the union of \(\Gamma\) and the equilibrium \(O\), i.e. converges to \(\Gamma \cup \{O\}\) as \(t\rightarrow \pm\infty\). We refer to such an orbit as '{\it homoclinic to homoclinic}' or '{\it super-homoclinic}' orbit.

\begin{myassumption}\label{assumption80}
There exists a super-homoclinic orbit \(\mathcal{S}\) of the transverse intersection of \(W^{s}_{\text{glo}}\left(\Gamma\right)\) and \(W^{u}_{\text{glo}}\left(\Gamma\right)\).
\end{myassumption}

\begin{figure}
\centering
\includegraphics[scale=.19]{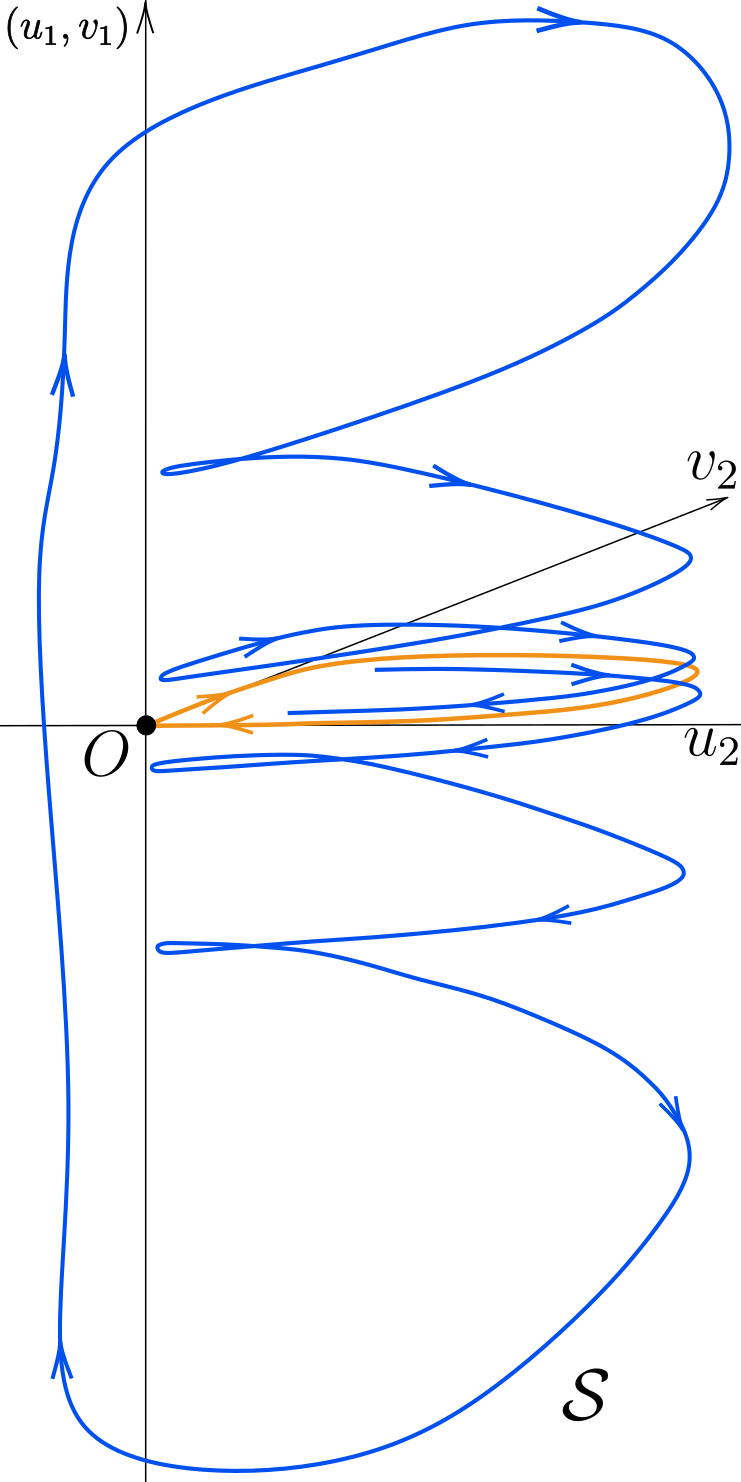}
\caption{\small The orbit \(\Gamma\) (brown) is homoclinic to the saddle equilibrium \(O\). The super-homoclinic orbit \(\mathcal{S}\) (blue) is homoclinic to \(\Gamma \cup \{O\}\).}
\label{Figurejuyv32utwve287f73}
\end{figure}

\begin{theoremC}\label{superhomoclinicthm}
Under Assumptions \ref{assumption10}-\ref{assumption50} and \ref{assumption80}, there exist infinitely many multi-pulse homoclinic loops in a small neighborhood of \(\mathcal{S} \cup \Gamma \cup \{O\}\) (i.e. the closure of \(\mathcal{S}\)).
\end{theoremC}

A similar result holds for a {\it homoclinic to homoclinic figure-eight}:

\begin{myassumption}\label{assumption90}
There exists a super-homoclinic orbit \(\mathcal{S}\) of the transverse intersection of \(W^{s}_{\text{glo}}\left(\Gamma_{1}\cup\Gamma_{2}\right)\) and \(W^{u}_{\text{glo}}\left(\Gamma_{1}\cup\Gamma_{2}\right)\).
\end{myassumption}

\begin{theoremC}\label{superhomoclinicthmforfigure8}
Under Assumptions \ref{assumption10}-\ref{assumption30}, \ref{assumption40}, \ref{assumption60} and \ref{assumption90}, there exist infinitely many multi-pulse homoclinic loops in a small neighborhood of \(\mathcal{S} \cup \{O\} \cup \Gamma_{1} \cup \Gamma_{2}\) (i.e. the closure of  \(\mathcal{S}\)).
\end{theoremC}

The multi-pulse homoclinic orbits in Theorem \ref{superhomoclinicthm} (and Theorem \ref{superhomoclinicthmforfigure8}) refer to homoclinic loops that intersect the cross-section \(\Pi^{s}\) (resp. \(\Pi_1^s \cup \Pi_2^s\)) in \(n\) points.
% \textcolor{red}{\(\Pi_1^s \cup \Pi_2^s\)???, 
We call such orbits \(n\)-pulse homoclinic. We prove that the existence of super-homoclinic orbits implies the existence of \(n\)-pulse homoclinic orbits for arbitrarily large \(n\).

\subsection{Coupled nonlinear Schr\"odinger equations}

The coupled nonlinear Schr\"odinger equation (CNLSE) is one of the basic models for light propagation. This equation also has various applications in different branches of physics since it appears as a universal model of behavior near a threshold of instability (see e.g. \cite{Kirrmann1992validity}). In this section, we discuss this equation as an application of our theory.

The CNLSE is written as
\begin{equation}\label{CNLSE0}
\begin{gathered}
i \Psi_{t} + \Psi_{xx} + 2\left(\alpha \left\vert \Psi\right\vert^{2} + \left\vert \Phi\right\vert^{2}\right)\Psi = 0,\\
i \Phi_{t} + \Phi_{xx} + 2\left(\left\vert \Psi\right\vert^{2} + \beta\left\vert \Phi\right\vert^{2}\right)\Phi = 0,
\end{gathered}
\end{equation}
where \(\Psi\) and \(\Phi\) are complex-valued functions of \(\left(t, x\right)\). We consider the case where \(\alpha\) and \(\beta\) are positive real constants. We consider the steady-state solutions of (\ref{CNLSE0}) which are of the form
\begin{equation*}
\Psi\left(t, x\right) = e^{i\omega_{1}^{2}t} \psi\left(x\right), \qquad \Phi\left(t, x\right) = e^{i\omega_{2}^{2}t} \phi\left(x\right),
\end{equation*}
for some real valued functions \(\psi\) and \(\phi\). By a rescaling, we can assume \(\omega_{1} = 1\) and \(\omega_{2} = \omega\) (\(\omega > 0\)). Thus, the stationary solutions of CNLSE satisfy
\begin{equation*}
\psi_{xx} = \psi - 2\left(\alpha \psi^{2} + \phi^{2}\right)\psi,\qquad
\phi_{xx} = \omega^{2} \phi - 2\left(\psi^{2} + \beta \phi^{2}\right)\phi.
\end{equation*}
Define \(\psi_{1}\left(x\right) = \psi\), \(\psi_{2}\left(x\right) = \psi_{x}\), \(\phi_{1}\left(x\right) = \phi\) and \(\phi_{2}\left(x\right) = \phi_{x}\). Then,
\begin{equation}\label{CNLSE2}
\begin{aligned}
\dot{\psi}_{1} &= \frac{\partial H}{\partial \psi_{2}} = \psi_{2},\qquad && \dot{\psi}_{2} = -\frac{\partial H}{\partial \psi_{1}} = \psi_{1} - 2\left(\alpha \psi_{1}^{2} + \phi_{1}^{2}\right)\psi_{1},\\
\dot{\phi}_{1} &= \frac{\partial H}{\partial \phi_{2}} = \phi_{2},\qquad && \dot{\phi}_{2} = -\frac{\partial H}{\partial \phi_{1}} = \omega^{2}\phi_{1} - 2\left(\psi_{1}^{2} + \beta \phi_{1}^{2}\right)\phi_{1},
\end{aligned}
\end{equation}
where \(H = \frac{1}{2}\left[\psi_{2}^{2} + \phi_{2}^{2} - \psi_{1}^{2} - \omega^{2} \phi_{1}^{2} + \alpha \psi_{1}^{4} + 2\psi_{1}^{2}\phi_{1}^{2} + \beta \phi_{1}^{4}\right]\). This system is Hamiltonian with two degrees of freedom. Diagonalizing the linear part reduces this system to
\begin{equation}\label{eq66eerio9ob31u48xvi}
\begin{aligned}
\dot{u}_{1} &= -u_{1} + E_{1}\left(u, v\right), \qquad &&\dot{v}_{1} = +v_{1} + \frac{1}{2} E_{1}\left(u, v\right),\\
\dot{u}_{2} &= -\omega u_{2} + E_{2}\left(u, v\right),\qquad &&\dot{v}_{2} = +\omega v_{2} - \frac{\omega}{2} E_{2}\left(u, v\right),
\end{aligned}
\end{equation}
where \(E_{1}\) and \(E_{2}\) are cubic functions of \(\left(u, v\right)\), and transforms the Hamiltonian \(H\) to the form \(H = u_{1}v_{1} - \omega u_{2}v_{2} + O\left(\|\left(u, v\right)\|^{4}\right)\). This system is invariant with respect to symmetry (\ref{eq425}) and the symmetry \(\left(u_{2}, v_{2}\right)\leftrightarrow \left(-u_{2}, -v_{2}\right)\). Therefore, assuming \(1\leq\omega\neq 2\), system (\ref{eq66eerio9ob31u48xvi}) meets all Assumptions \ref{assumption10}-\ref{assumption30} and \ref{assumption40}. In addition, it possesses a pair of homoclinic solutions (homoclinic figure-eight):
\begin{equation}\label{eq3435yyvyevydevf}
u_{1}\left(x\right) = 0,\quad
u_{2}\left(x\right) = \frac{\kappa\omega e^{\omega x}}{\sqrt{\beta} \cosh^{2}\left(\omega x\right)},\quad
v_{1}\left(x\right) = 0,\quad
v_{2}\left(x\right) = \frac{\kappa\omega^{2} e^{-\omega x}}{2\sqrt{\beta} \cosh^{2}\left(\omega x\right)}, \qquad \left(\kappa = \pm 1\right).
\end{equation}
These solutions correspond to the following solutions of (\ref{CNLSE0}):
\begin{equation}\label{eq9u8no8byik6tvjutygjyfwtf}
\Psi\left(t, x\right) = 0, \qquad \Phi\left(t, x\right) =  \pm\frac{\omega e^{i\omega^{2} t}}{\sqrt{\beta}\cosh\left(\omega x\right)}.
\end{equation}

We consider the case where the homoclinic figure-eight (\ref{eq3435yyvyevydevf}) is transverse, i.e. Assumption \ref{assumption60} is met. Therefore, the dynamics near this homoclinic figure-eight in the level \(\{H = 0\}\) can be analyzed by Theorems \ref{thm89bqyvrtyvtv} and \ref{Thm6892jjdibbea}. For \(\omega < 2\), Theorem \ref{thm89bqyvrtyvtv} implies that both forward and backward orbits of any point close to the homoclinic figure-eight leave a small neighborhood of it (in the level \(\{H = 0\}\)) unless it lies on the stable or unstable invariant manifolds of \(O\). For the case of \(\omega > 2\), in order to apply Theorem \ref{Thm6892jjdibbea}, one needs to reduce system (\ref{eq66eerio9ob31u48xvi}) to normal form (\ref{eq23000}) and compute the coefficients \(a_{i}\), \(b_{i}\), \(c_{i}\) and \(d_{i}\) (\(i=1, 2\)).

System (\ref{eq66eerio9ob31u48xvi}) is reversible with respect to the  linear involution \(u_{1} \leftrightarrow v_{1}\), \(u_{2} \leftrightarrow v_{2}\). In general, a system \(\dot{x} = f\left(x\right)\) on \(\mathbb{R}^{n}\) is said to be reversible with respect to an involution \(R\), i.e. a diffeomorphism on \(\mathbb{R}^{n}\) with the property \(R^{2} = \text{id}\), if \(dR \circ f = - f\circ R\). It is easily seen that when \(x\left(t\right)\) is a solution, so does \(R\circ x\left(-t\right)\). The reversibility of system (\ref{eq66eerio9ob31u48xvi}) implies
\begin{myprop}\label{Propiyvk7tcnyrxteh}
For \(A = \Gamma_{1}, \Gamma_{2}, \Gamma_{1} \cup \Gamma_{2}\), the manifold \(W^{s}_{loc}\left(A\right)\) is non-trivial if and only if \(W^{u}_{loc}\left(A\right)\) is non-trivial.
\end{myprop}

Reducing system (\ref{eq66eerio9ob31u48xvi}) to normal form (\ref{eq23000}) preserves the invariance of the system with respect to the symmetry \((u_{2}, v_{2})\leftrightarrow (-u_{2}, -v_{2})\). This implies that the loops \(\Gamma_{1}\) and \(\Gamma_{2}\) are symmetric, and \(a_{1} = a_{2} = a\), \(b_{1} = b_{2} = b\), \(c_{1} = c_{2} = c\) and \(d_{1} = d_{2} = d\). Because of the symmetry, Proposition \ref{Propiyvk7tcnyrxteh} implies %\textcolor{red}{(implies????)}
\begin{myprop}\label{Proployvutncrzkyiyolufu7675d}
Simultaneously, all the manifolds \(W^{u}_{loc}\left(\Gamma_{1}\right)\), \(W^{s}_{loc}\left(\Gamma_{1}\right)\), \(W^{u}_{loc}\left(\Gamma_{2}\right)\) and \(W^{s}_{loc}\left(\Gamma_{2}\right)\) are either trivial or non-trivial.
\end{myprop}

Concerning the coefficients \(a\), \(b\), \(c\) and \(d\), notice that the Hamiltonian structure of the equations implies that the map \(T^{\text{glo}}\) is area- and orientation-preserving; hence \(ad-bc = 1\). Therefore, \(\left(\begin{smallmatrix} a & b\\ c & d \end{smallmatrix}\right)^{-1} = \left(\begin{smallmatrix} d & -b\\ -c & a \end{smallmatrix}\right)\), so the reversibility implies \(b = -c\). Recall also that the transversality condition implies \(d\neq 0\).

It follows from Theorem \ref{Thm6892jjdibbea} that if \(bd >0\), then the local unstable invariant manifold of each of the loops \(\Gamma_{1}\) and \(\Gamma_{2}\) is non-trivial, while the local unstable invariant manifold of the homoclinic figure-eight \(\Gamma_{1} \cup \Gamma_{2}\) is trivial (i.e. coincides with \(\Gamma_{1} \cup \{O\} \cup \Gamma_{2}\)). In contrast, when \(bd <0\), the local unstable invariant manifold of the homoclinic figure-eight is non-trivial, while the local unstable invariant manifold of each of the loops \(\Gamma_{1}\) and \(\Gamma_{2}\) is trivial. The same conclusion holds for the corresponding stable manifolds. This analysis together with Propositions \ref{Propiyvk7tcnyrxteh} and \ref{Proployvutncrzkyiyolufu7675d} yields
\begin{myprop}\label{Prop9809byvt6eruwriayeuy}
Let \(\omega > 2\) and suppose that the coefficients \(b\), \(c\) and \(d\) are non-zero. Then, one (and only one) of the following two scenarios holds:
\begin{enumerate}[(i)]
\item The manifolds \(W^{u}_{loc}\left(\Gamma_{1}\cup\Gamma_{2}\right)\) and \(W^{s}_{loc}\left(\Gamma_{1}\cup\Gamma_{2}\right)\) are non-trivial, i.e. \(bd = -cd < 0\).
\item The manifolds \(W^{u}_{loc}\left(\Gamma_{1}\right)\), \(W^{s}_{loc}\left(\Gamma_{1}\right)\), \(W^{u}_{loc}\left(\Gamma_{2}\right)\) and \(W^{s}_{loc}\left(\Gamma_{2}\right)\) are non-trivial, i.e. \(bd = -cd > 0\).
\end{enumerate}
\end{myprop}

To figure out which of the scenarios above happens for CNLSE, one needs to find the corresponding coefficients \(a\), \(b\), \(c\) and \(d\). For any particular values of \(\alpha\), \(\beta\) and \(\omega\), this can be done numerically by solving the linearization of system (\ref{CNLSE2}) along the solution (\ref{eq9u8no8byik6tvjutygjyfwtf}). It is also easy to show that as \(\omega \rightarrow \infty\), the coefficient \(d\) changes sign infinitely many times. When \(d=0\), we have \(b^{2} = 1\), so crossing \(d=0\) leads to a change of sign of \(bd\). Thus, both cases of the above proposition are realized in CNLSE. 
In either case, 
%Proposition \ref{Prop9809byvt6eruwriayeuy} states that 
there are non-trivial local stable and unstable invariant manifolds of the homoclinic orbits in the CNLSE. Globalizing these manifolds, we conjecture that they intersect transversely along some super-homoclinic orbits:
\begin{myconjecture}
The coupled nonlinear Schr\"odinger equations given by (\ref{CNLSE0}) possess transverse super-homoclinic orbits.
\end{myconjecture}

%We expect this conjecture to be true since 
Our conjecture is based on the fact that the stable and unstable manifolds of the homoclinic loops are 2-dimensional manifolds  lying in the same compact 3-dimensional energy level, and - by Poincar\'e recurrence - come infinitely close to each other. Hence, they are likely to intersect transversely along super-homoclinic orbits. Moreover, numerical evidence \cite{EleonskyKorolev1993,EleonskyCodinghomoclinic1996,ClassificationSolitaryWaves1997,Yang1998} points to the existence of infinitely many multi-pulse homoclinic orbits in the CNLSE. This supports our conjecture since, by Theorems \ref{superhomoclinicthm} and \ref{superhomoclinicthmforfigure8}, the existence of these multi-pulse homoclinics might be a bi-product of the existence of super-homoclinic orbits.

\subsection{The methods and organization of the paper}

The standard approach for investigating the dynamics near homoclinic orbits is to study the Poincar\'e maps along these orbits. The main difficulty in dealing with these maps is that the Poincar\'e map along a homoclinic orbit is a singular map defined on a non-trivial domain. Let \(\Gamma\) be a homoclinic orbit, \(\Sigma\) be a small cross-section to it, \(T\) be the Poincar\'e map defined on some domain \(\mathcal{D}\subset \Sigma\), and \(M\) be the intersection point of \(\Gamma\) and \(\Sigma\). The point \(M\) does not belong to the domain \(\mathcal{D}\), however, it is in the closure of this set. The domain \(\mathcal{D}\) in our case consists of several (at least two) connected components each with non-empty interior, while the set \(\mathcal{D}\cup \{M\}\) is connected, but \(\mathcal{D} \cup \{M\}\) does not contain an open neighborhood of \(M\) in \(\Sigma\). The point \(M\) is a singularity for the Poincar\'e map \(T\): as \(z\rightarrow M\) for \(z\in \mathcal{D}\), we have \(\|dT\left(z\right)\|\rightarrow\infty\), so, extending the map \(T\) to \(\mathcal{D} \cup \{M\}\) by defining \(T\left(M\right) = M\) and turning \(M\) to the fixed point of \(T\), does not remove this singularity.

Due to these properties, the smooth theory of Hadamard-Perron cannot be applied directly to study the invariant manifolds of the Poincar\'e maps along homoclinics. Our approach for investigating the invariant manifolds is applying the method of Shilnikov cross-maps \cite{Dimabook,GonchenkoCrossmap2010}. % \textcolor{red}{references?}. 
To describe this method, assume that the Poincar\'e map \(T\) is written by
\begin{equation*}
\begin{aligned}
\overline{x} &= f\left(x, y\right),\\
\overline{y} &= g\left(x, y\right)
\end{aligned}
\end{equation*}
where \(f\) and \(g\) are some functions defined on \(\mathcal{D}\) such that
\begin{equation*}
\lim_{\left(x, y\right)\rightarrow\left(0, 0\right)} f\left(x, y\right) = \lim_{\left(x, y\right)\rightarrow\left(0, 0\right)} g\left(x, y\right) = \left(0, 0\right),
\end{equation*}
and the point \(M\) corresponds to \(\left(x, y\right) = \left(0, 0\right)\). Suppose that \(y\) in the second equation can be solved in terms of \(\left(x,\overline{y}\right)\), i.e. \(y = G\left(x,\overline{y}\right)\) for some function \(G\). This leads us to introduce the map \(T^{\times}\) defined by
\begin{equation*}
\begin{aligned}
\overline{x} &= F\left(x, \overline{y}\right),\\
y &= G\left(x,\overline{y}\right),
\end{aligned}
\end{equation*}
where \(F\left(x, \overline{y}\right) = f\left(x, G\left(x,\overline{y}\right)\right)\). In other words, the Poincar\'e map \(T\) takes \(\left(x, y\right)\) to \(\left(\overline{x}, \overline{y}\right)\) if and only if the cross-map \(T^{\times}\) maps \(\left(x, \overline{y}\right)\) to \(\left(\overline{x}, y\right)\). Denote the domain of \(T^{\times}\) by \(\mathcal{D}^{\times}\). Note that the point \(M = \left(x, \overline{y}\right) = \left(0, 0\right)\) is not in the domain \(\mathcal{D}^{\times}\), however, it lies in the closure of \(\mathcal{D}^{\times}\).

The advantage of dealing with the cross-map \(T^{\times}\) over the Poincar\'e map \(T\) is that despite \(\|dT\left(x, y\right)\| \rightarrow\infty\) as \(\left(x, y\right)\rightarrow M\), the limit \(\lim_{\left(x, \overline{y}\right)\rightarrow M} T^{\times}\left(x, \overline{y}\right)\) may exist. This property, if the \(\left(x, y\right)\) coordinates system on \(\Sigma\) is chosen appropriately, enables us to extend the cross-map \(T^{\times}\) to an open neighborhood of \(M\) smoothly. Then, according to Theorem \ref{thm3000} (see Appendix \ref{Invariantmanifoldscrossmaps}), if \textit{extended} \(T^{\times}\) satisfies certain properties, the Poincar\'e map \(T\) possesses an invariant manifold that contains the \(\omega\)-limit points of every forward orbit of the domain. This is exactly the procedure that we follow in this paper to prove the existence of invariant manifolds of the Poincar\'e maps along the homoclinic orbits.

In order to obtain necassary estimates for the Poincar\'e map, we first bring our system near the equilibrium state \(O\) to a normal form. Notice that our system is not necessarily linearizable. Indeed, since the spectrum of the linear part of the system is \(\{-\lambda_{2}, -\lambda_{1}, \lambda_{1}, \lambda_{2}\}\), for some \(0 < \lambda_{1} \leq \lambda_{2}\), resonances prevent the linearization. In this paper, we follow the approach of Ovsyaninkov and Shilnikov \cite{OvsyannikovandShilnikov1986} (see also \cite{DimaMultipulse2001} and \cite{Dimabook}) for the normal form reduction. In contrast to the standard normal form approach in which the system near a hyperbolic equilibrium reduces to a polynomial vector field that consists of only resonant terms up to some order (see e.g. \cite{NonlocalbifurcationsIlyashenko} and \cite{BronsteinKopanskii}), in the approach of \cite{OvsyannikovandShilnikov1986}, some non-resonant terms remain in the normal form, while some resonant terms of low-orders are eliminated.

Once the system is brought to a normal form, we need to investigate the behavior of the orbits near the equilibrium state \(O\). This enables us to compute the Poincar\'e maps. To do this, we apply the method of successive approximations (more specifically, Shilnikov's method of solving boundary value problems, see \cite{Dimabook}) to estimate the flow near the equilibrium \(O\).

This paper is organized as follows: Section \ref{localdynamicsnearO} is dedicated to the study of the local map \(T^{\text{loc}}\). In Section \ref{Setupsection76q5y6d45ljbiywbux}, we define this map and its domain precisely. Then, in Section \ref{choiceofcoordinatesnearO}, we bring our system near the equilibrium state \(O\) to a normal form. In Section \ref{trajectoriesnearO}, we investigate the behavior of the orbits near the equilibrium state \(O\). Finally, in Section \ref{localmapsection}, we analyze the domain and the behavior of the local map.

In Section \ref{Dynamics nearhomoclinic}, we use the results of Section \ref{localdynamicsnearO} to study the dynamics near the homoclinic orbits. In Section \ref{Somenotationsandsetting}, we introduce some notations. In Section \ref{Dynamicsnearasinglehomoclinictrivialcase}, we study the dynamics near the homoclinic orbit \(\Gamma\) when \(\lambda_{2} < 2\lambda_{1}\). Theorems \ref{thmkuyvuyuy090988} is proved in this section. The dynamics near \(\Gamma\) when \(2\lambda_{1} < \lambda_{2}\) is studied in Section \ref{kjhliub812ftf878g3x7c}. We prove Theorem \ref{Invariantmanifoldthm} in this section. Theorem \ref{Th782ju2iu22} is also proved in these two sections. The case of homoclinic figure-eight is studied in Section \ref{Dynamicsnearfigure8homoclinic}. The proofs of Theorems \ref{Thmyipoinoib53}, \ref{thm89bqyvrtyvtv} and \ref{Thm6892jjdibbea} are provided in this section. Finally, we discuss the case of superhomoclinics in Section \ref{Superhomoclinicsubsection}. We prove Theorems \ref{superhomoclinicthm} and \ref{superhomoclinicthmforfigure8} in this section.

Most of the technical lemmas and calculations are postponed to appendices. We also give a brief introduction to the method of cross-maps in Appendix \ref{Invariantmanifoldscrossmaps}.

\section{Analysis near the equilibrium state \texorpdfstring{\(O\)}{Lg}}\label{localdynamicsnearO}

\subsection{Set-up and notations}\label{Setupsection76q5y6d45ljbiywbux}

Our approach for studying the dynamics near the homoclinic loop \(\Gamma\) (and homoclinic figure-eight \(\Gamma_{1} \cup \Gamma_{2}\)) is based on the study of the behavior of the corresponding Poincar\'e map(s). As was mentioned earlier, the Poincar\'e map \(T\) along the homoclinic loop \(\Gamma\) can be written as the composition of a global and a local map. This section is dedicated to the study of the behavior of the local map \(T^{loc}\). To this end, we first need to choose appropriate coordinates near the equilibrium state \(O\) of system (\ref{eq300}). This is done in Section \ref{choiceofcoordinatesnearO} below. We consider three different cases of \(\lambda_{1} = \lambda_{2}\), \(\lambda_{1} < \lambda_{2} < 2\lambda_{1}\) and \(2\lambda_{1} < \lambda_{2}\), and introduce a specific normal form for each case. In Section \ref{trajectoriesnearO}, we employ the Shilnikov technique for solving boundary value problems to compute the flow near the equilibrium \(O\). This allows us to find an approximation for the local map. Finally, in Section \ref{localmapsection}, we study the behavior of this map and investigate some of its properties.

In comparison to the global map, the local map has more complicated behavior. Indeed, \(T^{glo}\) is a diffeomorphism and can be approximated by its Taylor polynomial while the local map \(T^{loc}\) is a singular map with a non-trivial domain.

Let us now give a more precise meaning to the above terminologies. Recall the cross-sections \(\Pi^{s}\) and \(\Pi^{u}\). In all of the normal forms considered in Section \ref{choiceofcoordinatesnearO}, the local stable and local unstable as well as the local strong stable and local strong unstable invariant manifolds of \(O\) are straightened. Therefore, the homoclinic loop \(\Gamma\) intersects \(\Pi^{s}\) and \(\Pi^{u}\) at \(M^{s} = \left(0,\delta,0,0\right)\) and \(M^{u} = \left(0,0,0,\delta\right)\), respectively. As it is proved later (see Section \ref{Coordinates-on-cross-sections}), we can choose a \(\left(u_{1}, v_{1}\right)\) coordinate-system on each of these cross-sections. Both \(M^{s}\) and \(M^{u}\) correspond to \(\left(0, 0\right)\) in this coordinate-system.

\begin{figure}
\centering
\begin{subfigure}{0.4\textwidth}
\centering
\includegraphics[scale=0.14]{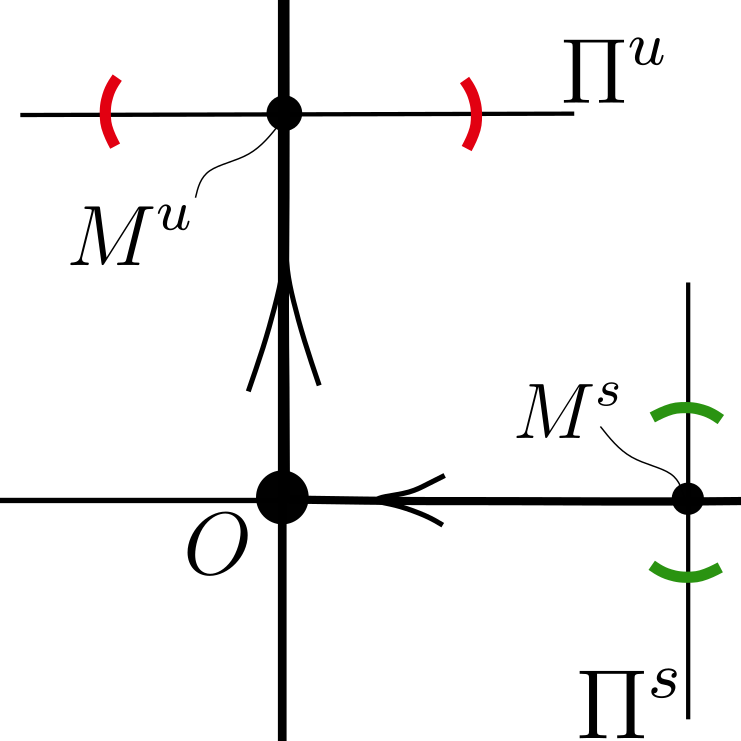}
\end{subfigure}
\begin{subfigure}{0.4\textwidth}
\centering
\includegraphics[scale=0.14]{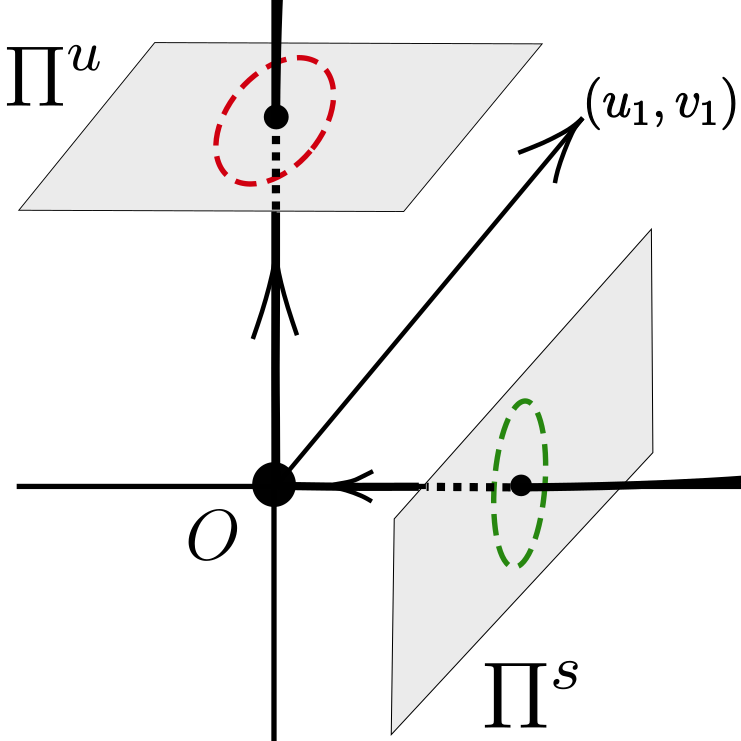}
\end{subfigure}
\caption{\small The \(\epsilon\)-ball around \(M^{s}\) in \(\Pi^{s}\) and the \(\epsilon_{u}\)-ball around \(M^{u}\) in \(\Pi^{u}\) are shown by green and red colors, respectively. The domain \(\mathcal{D}\) of the Poincar\'e map \(T\) is the set of the points \(\left(u_{10}, v_{10}\right)\) in the green ball whose forward orbits intersect \(\Pi^{u}\) at \(\left(u_{1\tau}, v_{1\tau}\right)\) in the red ball (see relation (\ref{conditionfordomainD})).}
\label{Figuruinub8vqtr6}
\end{figure}

Consider a point \(\left(u_{10}, v_{10}\right)\) on \(\Pi^{s}\) close to \(M^{s}\) (e.g. the green point in Figure \ref{Figure00m9n8nb4y2bc}) whose forward orbit goes along the homoclinic loop \(\Gamma\), after a certain time \(\tau\) it crosses \(\Pi^{u}\) at a point \(\left(u_{1\tau}, v_{1\tau}\right)\) (e.g. the blue point in Figure \ref{Figure00m9n8nb4y2bc}), and after a finite time it comes back to \(\Pi^{s}\) at a point \(\left(\overline{u}_{10}, \overline{v}_{10}\right)\) (e.g. the red point in Figure \ref{Figure00m9n8nb4y2bc}). Obviously, \(\tau \rightarrow \infty\) as \(\left(u_{10}, v_{10}\right)\rightarrow M^{s}\). Let \(\mathcal{D} \subset \Pi^{s}\) be the set of all such points \(\left(u_{10}, v_{10}\right)\) that satisfy
\begin{equation}\label{conditionfordomainD}
\left\Vert \left(u_{10}, v_{10}\right)\right\Vert < \epsilon\quad \text{and} \quad\left\Vert \left(u_{1\tau}, v_{1\tau}\right)\right\Vert < \epsilon_{u},
\end{equation}
for some sufficiently small constants \(0 < \epsilon \leq \epsilon_{u} < \delta\) (see Figure \ref{Figuruinub8vqtr6}). It is trivial that \(M^{s}\notin \mathcal{D}\). When \(\mathcal{D}\neq \emptyset\), we define the Poincar\'e map \(T: \mathcal{D} \rightarrow \Pi^{s}\) by \(\left(u_{10}, v_{10}\right) \longmapsto \left(\overline{u}_{10}, \overline{v}_{10}\right)\). The local map \(T^{loc}: \mathcal{D} \rightarrow \Pi^{u}\) is defined by
\begin{equation}\label{localmap}
\left(u_{10}, v_{10}\right) \longmapsto \left(u_{1\tau}, v_{1\tau}\right).
\end{equation}
The global map is defined on the \(\epsilon_{u}\)-ball \(\mathcal{B}_{\epsilon_{u}}\) in \(\Pi^{u}\) centered at \(M^{u}\), i.e.  \(T^{glo}: \mathcal{B}_{\epsilon_{u}}\rightarrow \Pi^{s}\), and its restriction to \(T^{loc}\left(\mathcal{D}\right)\subset \mathcal{B}_{\epsilon_{u}}\) is 
\begin{equation}\label{eq789o0o1ok2nnyxuj}
\left(u_{1\tau}, v_{1\tau}\right) \longmapsto \left(\overline{u}_{10}, \overline{v}_{10}\right).
\end{equation}
Obviously, \(T = T^{glo} \circ T^{loc}\).

\begin{figure}
\centering
\includegraphics[scale=.22]{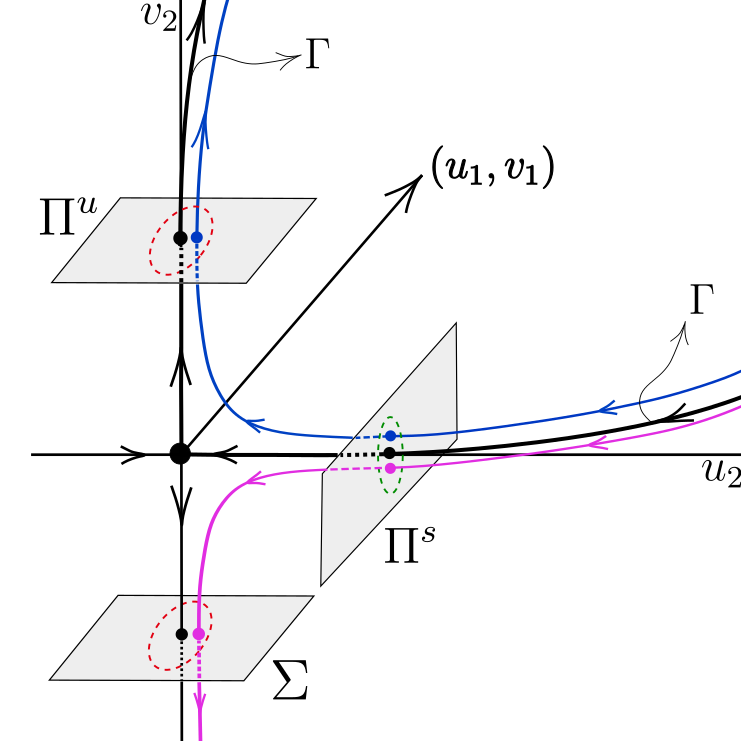}
\caption{\small The domain \(\mathcal{D}\) of the Poincar\'e map is defined as the set of the points on \(\Pi^{s}\) close to \(M^{s}\) that go along the homoclinic loop \(\Gamma\) and intersect \(\Pi^{u} = \{v_{2} = \delta\} \cap \{H=0\}\) at points close to \(M^{u}\). For instance, the blue point on \(\Pi^{s}\) belongs to \(\mathcal{D}\). Similarly, we define \(\mathbb{D}\) as the set of the points on \(\Pi^{s}\) close to \(M^{s}\) that go along \(\Gamma\) until they get close to \(O\) and then go along the negative side of \(v_{2}\)-axis and intersect the cross-section \(\Sigma = \{v_{2} = -\delta\} \cap \{H=0\}\) at points close to the point of the intersection of \(\Sigma\) and \(v_{2}\)-axis (e.g. the pink point on \(\Pi^{s}\) belongs to \(\mathbb{D}\)).}
\label{Figure8767kutvtyqq36cqq}
\end{figure}

Not every orbit starting from \(\Pi^{s}\) goes along \(\Gamma\) and intersects the cross-section \(\Pi^{u}\). Trivial examples are the orbits that start at \(W_{\text{loc}}^{s}\left(O\right) \cap \Pi^{s}\). Other examples are the orbits that go along the other branch of \(W^{u}_{\text{loc}}\left(O\right)\) (negative side of \(v_{2}\)-axis). Consider a cross-section \(\Sigma = \lbrace v_{2} = -\delta\rbrace \cap \lbrace H = 0\rbrace\) to the negative branch of \(W^{u}_{\text{loc}}\left(O\right)\). It will be shown that \(\left(u_{1}, v_{1}\right)\)-coordinates can be chosen on this cross-section. Then
\begin{mydefn}\label{Defniopinouib97v76c5c}
We denote by \(\mathbb{D}\) the set of the points \(\left(u_{10}, v_{10}\right)\) on \(\Pi^{s}\) close to \(M^{s}\) whose forward orbits go along the negative branch of \(W^{u}_{\text{loc}}\left(O\right)\), and after a certain time \(\tau\) they cross \(\Sigma\) at \(\left(u_{1\tau}, v_{1\tau}\right)\) such that (\ref{conditionfordomainD}) holds (see Figure \ref{Figure8767kutvtyqq36cqq}).
\end{mydefn}

For the case of homoclinic figure-eight, we define the domains \(\mathcal{D}\) and \(\mathbb{D}\) for each loop:
\begin{mynotation}\label{Notation76v6uc5x54x4wez}
For \(i=1, 2\), we denote by \(\mathcal{D}^{i}\) and \(\mathbb{D}^{i}\) the corresponding domains \(\mathcal{D}\subset \Pi^{s}_{i}\) and \(\mathbb{D}\subset \Pi^{s}_{i}\) of the loop \(\Gamma_{i}\), respectively.
\end{mynotation}

An orbit starting from \(\mathcal{D}^{1}\subset \Pi^{s}_{1}\) (resp. \(\mathcal{D}^{2}\subset \Pi^{s}_{2}\)) goes along \(\Gamma_{1}\) (resp. \(\Gamma_{2}\)) and intersects \(\Pi^{u}_{1}\) (resp. \(\Pi^{u}_{2}\)), while an orbit which starts from \(\mathbb{D}^{1}\subset \Pi^{s}_{1}\) (resp. \(\mathbb{D}^{2}\subset \Pi^{s}_{2}\)) goes along the negative (resp. positive) side of \(v_{2}\)-axis and intersects \(\Pi^{u}_{2}\) (resp. \(\Pi^{u}_{1}\)).

\begin{figure}
\centering
\begin{subfigure}{0.4\textwidth}
\centering
\includegraphics[scale=0.21]{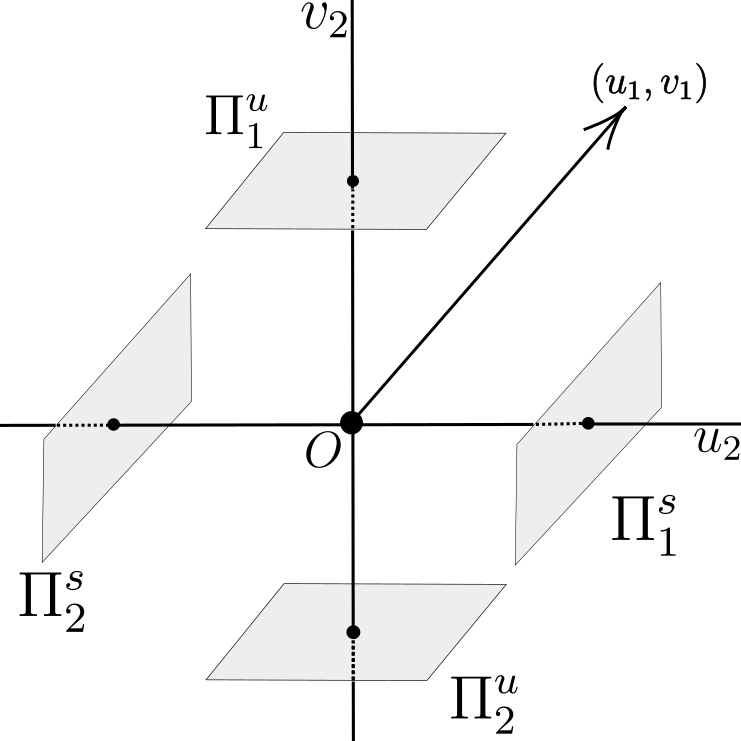}
\end{subfigure}
\begin{subfigure}{0.4\textwidth}
\centering
\includegraphics[scale=0.21]{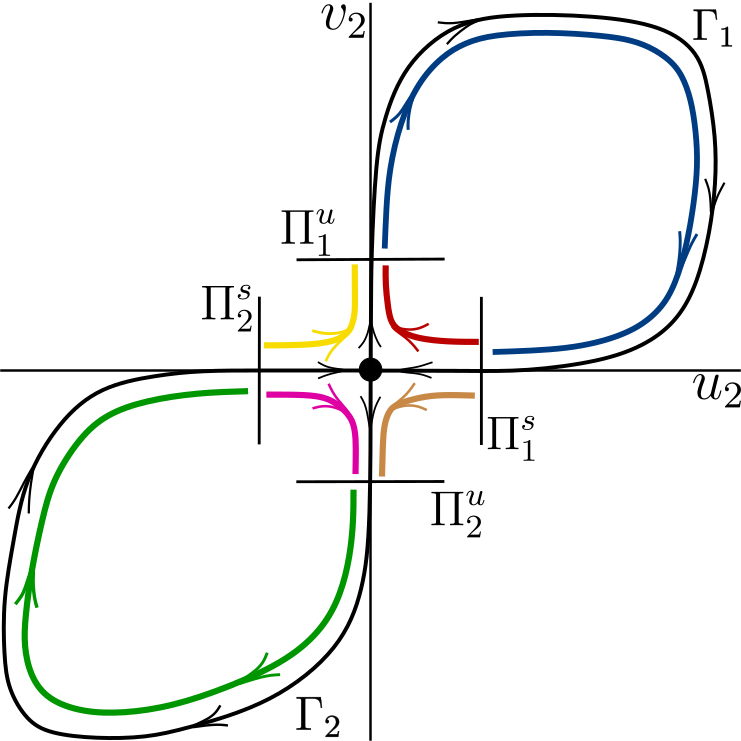}
\end{subfigure}
\caption{\small (left) The positions of the cross-sections \(\Pi^{s}_{1}\), \(\Pi^{u}_{1}\), \(\Pi^{s}_{2}\) and \(\Pi^{u}_{2}\) are shown. (right) \(\Gamma_{1}\) and \(\Gamma_{2}\) are homoclinic orbits. The blue, brown, green, yellow, red and pink curves correspond to the maps \(T_{1}^{\text{glo}}\), \(T_{12}^{\text{loc}}\), \(T_{2}^{\text{glo}}\), \(T_{21}^{\text{loc}}\), \(T_{1}^{\text{loc}}\) and \(T_{2}^{\text{loc}}\), respectively. The Poincar\'e maps \(T_{1}\) (along \(\Gamma_{1}\)) and \(T_{2}\) (along \(\Gamma_{2}\)) are defined by \(T_{1} = T_{1}^{\text{glo}} \circ T_{1}^{\text{loc}}\) and \(T_{2} = T_{2}^{\text{glo}} \circ T_{2}^{\text{loc}}\), respectively.}
\label{Figure673bob8i7vrq7cvraa}
\end{figure}

We introduced the Poincar\'e, local and global maps along a single homoclinic loop above. For the case of homoclinic figure eight, we also define these maps for each loop:
\begin{mynotation}\label{Notation787v76ecx65c3e1ddss}
We denote by \(T_{i}\), \(T^{\text{loc}}_{i}\) and \(T^{\text{glo}}_{i}\) the Poincar\'e, local and global maps along \(\Gamma_{i}\) (\(i=1, 2\)), respectively (see Figure \ref{Figure673bob8i7vrq7cvraa}).
\end{mynotation}
To study the case of a homoclinic figure-eight, we consider two extra local maps:
\begin{mydefn}\label{Defn01ddsdutyc11t1c1asaq}
We define the map \(T^{\text{loc}}_{12}: \mathbb{D}^{1}\subset \Pi^{s}_{1} \rightarrow \Pi^{u}_{2}\) (\(T^{\text{loc}}_{21}: \mathbb{D}^{2}\subset \Pi^{s}_{2} \rightarrow \Pi^{u}_{1}\)) by \(\left(u_{10}, v_{10}\right) \mapsto \left(u_{1\tau}, v_{1\tau}\right)\) where \(\left(u_{10}, v_{10}\right)\in \mathbb{D}^{1}\) (\(\in \mathbb{D}^{2}\)) and \(\left(u_{1\tau}, v_{1\tau}\right)\in\Pi_{2}^{u}\) (\(\in\Pi_{1}^{u}\)) (see Figure \ref{Figure673bob8i7vrq7cvraa}).
\end{mydefn}

\subsection{Choice of coordinates near the equilibrium state \texorpdfstring{\(O\)}{Lg}}\label{choiceofcoordinatesnearO}

This section is dedicated to finding suitable coordinate systems near the equilibrium state \(O\). As it was mentioned above, we consider three different cases of \(\lambda_{1} = \lambda_{2}\), \(\lambda_{1} < \lambda_{2} < 2\lambda_{1}\) and \(2\lambda_{1} < \lambda_{2}\), and for each case we bring system (\ref{eq300}) into a particular normal form. The proofs of the results stated below are postponed to Appendix \ref{Proofs of normal form Lemmas}. We start with the following:

\begin{mylem}\label{Nftheorem0}
Consider system (\ref{eq300}) and first integral (\ref{eq400}). There exists a \({\mathcal{C}}^{\infty}\)-smooth change of coordinates which brings system (\ref{eq300}) to the form
\begin{equation}\label{eq74wt366cuw6gt5uw6v01092}
\begin{aligned}
\dot{u}_{1} &= -\lambda_{1}u_{1} + f_{11}(u_{1}, u_{2}, v_{1}, v_{2}) u_{1} + f_{12}(u_{1}, u_{2}, v_{1}, v_{2}) u_{2},\\
\dot{u}_{2} &= -\lambda_{2}u_{2} + f_{21}(u_{1}, u_{2}, v_{1}, v_{2}) u_{1} + f_{22}(u_{1}, u_{2}, v_{1}, v_{2}) u_{2},\\
\dot{v}_{1} &= +\lambda_{1}v_{1} + g_{11}(u_{1}, u_{2}, v_{1}, v_{2}) v_{1} + g_{12}(u_{1}, u_{2}, v_{1}, v_{2}) v_{2},\\
\dot{v}_{2} &= +\lambda_{2}v_{2} + g_{21}(u_{1}, u_{2}, v_{1}, v_{2}) v_{1} + g_{22}(u_{1}, u_{2}, v_{1}, v_{2}) v_{2},
\end{aligned}
\end{equation}
where the functions \(f_{ij}\), \(g_{ij}\) are \({\mathcal{C}}^{\infty}\)-smooth and vanish at the origin, i.e.
\begin{equation}\label{equnsfuhnrx9uwn984urb348}
f_{ij}\left(0,0,0,0\right) = g_{ij}\left(0,0,0,0\right) = 0,
\end{equation}
and transforms first integral (\ref{eq400}) to
\begin{equation}\label{eq98n9uo8wbr8bau4b819ub92q1u8}
H = \lambda_{1}u_{1}v_{1} - \lambda_{2}u_{2}v_{2}.
\end{equation}
Moreover, system (\ref{eq74wt366cuw6gt5uw6v01092}) remains invariant with respect to symmetry (\ref{eq425}). In particular,
\begin{equation}\label{eq24500}
f_{12}(0, u_{2}, 0, v_{2}) \equiv 0,\quad g_{12}(0, u_{2}, 0, v_{2}) \equiv 0.
\end{equation}
\end{mylem}

The statement of Lemma \ref{Nftheorem0} holds for arbitrary \(0 < \lambda_{1}\leq\lambda_{2}\). However, we will particularly use this normal form for analyzing the case \(\lambda_{1} = \lambda_{2}\).

\begin{mylem}\label{Nftheorem}
Consider system (\ref{eq300}) and first integral (\ref{eq400}), and assume \(\lambda_{1} <\lambda_{2}\). There exists a \({\mathcal{C}}^{\infty}\)-smooth change of coordinates which brings system (\ref{eq300}) to the form
\begin{equation}\label{eq19000}
\begin{aligned}
\dot{u}_{1} &= -\lambda_{1}u_{1} + f_{11}\left(u_{1}, v\right) u_{1} + f_{12}\left(u_{1}, u_{2}, v\right) u_{2},\\
\dot{u}_{2} &= -\lambda_{2}u_{2} + f_{21}\left(u_{1}, v\right)u_{1} + f_{22}\left(u_{1}, u_{2}, v\right) u_{2},\\
\dot{v}_{1} &= +\lambda_{1}v_{1} + g_{11}\left(u, v_{1}\right) v_{1} + g_{12}\left(u, v_{1}, v_{2}\right) v_{2},\\
\dot{v}_{2} &= +\lambda_{2}v_{2} + g_{21}\left(u, v_{1}\right)v_{1} +  g_{22}\left(u, v_{1}, v_{2}\right) v_{2},
\end{aligned}
\end{equation}
where the functions \(f_{ij}\), \(g_{ij}\) are \({\mathcal{C}}^{\infty}\)-smooth and satisfy the identities
\begin{equation}\label{eq20000}
\begin{gathered}
f_{11}(0,v) \equiv 0,\quad f_{11}(u_{1},0) \equiv 0,\quad f_{12}(u,0) \equiv 0,\quad f_{21}(0,v) \equiv 0,\quad f_{22}(0, v) \equiv 0,\\
g_{11}(u,0) \equiv 0,\quad g_{11}(0,v_{1}) \equiv 0,\quad g_{12}(0,v) \equiv 0,\quad g_{21}(u,0) \equiv 0,\quad g_{22}\left(u,0\right) \equiv 0.
\end{gathered}
\end{equation}
This change of coordinates transforms first integral (\ref{eq400}) to
\begin{equation}\label{eq22000}
H = \lambda_{1}u_{1}v_{1} \left[1 + H_{1}\left(u, v\right)\right]-\lambda_{2}u_{2}v_{2} \left[1 + H_{2}\left(u, v\right)\right],
\end{equation}
where \(H_{1}\) and \(H_{2}\) are \(\mathcal{C}^{\infty}\) functions that vanish at \(O\). We can write (\ref{eq22000}) as
\begin{equation}\label{eq22000yuio}
H = \lambda_{1}u_{1}v_{1} \left[1 + o\left(1\right)\right] - \lambda_{2}u_{2}v_{2} \left[1 + o\left(1\right)\right].
\end{equation}
Moreover, normal form (\ref{eq19000}) and first integral (\ref{eq22000}) remain invariant with respect to symmetry (\ref{eq425}). In particular, (\ref{eq24500}) holds.
\end{mylem}

The statement of Lemma \ref{Nftheorem} holds for arbitrary \(\lambda_{1} <\lambda_{2}\). However, we will particularly use this normal form to analyze the local dynamics near \(O\) when \(\lambda_{1} <\lambda_{2} < 2\lambda_{1}\). The normal form that is used for analyzing the case \(2\lambda_{1} < \lambda_{2}\) is given by the following:
\begin{mylem}\label{Nftheorem2}
Consider system (\ref{eq300}) and first integral (\ref{eq400}) and assume \(2\lambda_{1} <\lambda_{2}\). Let \(q\) be the largest integer such that \(q\lambda_{1} < \lambda_{2}\). There exists a \({\mathcal{C}}^{q}\)-smooth change of coordinates which brings system (\ref{eq300}) to the form
\begin{equation}\label{eq23000}
\begin{aligned}
\dot{u}_{1} &= -\lambda_{1}u_{1} + f_{11}\left(u_{1}, v\right) u_{1} + f_{12}\left(u_{1}, u_{2}, v\right) u_{2},\\
\dot{u}_{2} &= -\lambda_{2}u_{2} + f_{22}\left(u_{1}, u_{2}, v\right) u_{2},\\
\dot{v}_{1} &= +\lambda_{1}v_{1} + g_{11}\left(u, v_{1}\right) v_{1} + g_{12}\left(u, v_{1}, v_{2}\right) v_{2},\\
\dot{v}_{2} &= +\lambda_{2}v_{2} + g_{22}\left(u, v_{1}, v_{2}\right) v_{2},
\end{aligned}
\end{equation}
where \(f_{ij}\) and \(g_{ij}\) are \(\mathcal{C}^{q-1}\)-smooth and satisfy  identities (\ref{eq20000}). This change of coordinates transforms first integral (\ref{eq400}) to
\begin{equation}\label{eqhhuhut5989b8b8y4qlz37y7t}
H = \lambda_{1}u_{1}v_{1} \left[1 + H_{1}\left(u, v\right)\right] -\lambda_{2} u_{2}v_{2} \left[1 + H_{2}\left(u, v\right)\right] + u_{2}v_{1}^{2} H_{3}\left(u, v\right) + v_{2}u_{1}^{2} H_{4}\left(u, v\right),
\end{equation}
where \(H\) is \(\mathcal{C}^{q}\), and \(H_{1}\), \(H_{2}\), \(H_{3}\) and \(H_{4}\) are some \(\mathcal{C}^{q-1}\), \(\mathcal{C}^{q}\), \(\mathcal{C}^{q-2}\) and \(\mathcal{C}^{q-2}\) functions, respectively, such that \(H_{1}(O) = H_{2}(O) = 0\). Moreover, system (\ref{eq23000}) and first integral (\ref{eqhhuhut5989b8b8y4qlz37y7t}) remain invariant with respect to symmetry (\ref{eq425}). In particular, (\ref{eq24500}) holds.
\end{mylem}

\begin{myrem}
For simplicity, we can write (\ref{eqhhuhut5989b8b8y4qlz37y7t}) as
\begin{equation*}
H = \lambda_{1}u_{1}v_{1} \left[1 + o\left(1\right)\right] - \lambda_{2}u_{2}v_{2} \left[1 + o\left(1\right)\right] + u_{2}v_{1}^{2} O\left(1\right) + v_{2}u_{1}^{2} O\left(1\right).
\end{equation*}
\end{myrem}

A common structure of all of normal forms (\ref{eq74wt366cuw6gt5uw6v01092}), (\ref{eq19000}) and (\ref{eq23000}) is that the local stable and unstable as well as the local strong stable and strong unstable invariant manifolds of the equilibrium \(O\) are straightened, i.e. \(W^{s}_{\text{loc}} = \{v = 0\}\), \(W^{u}_{\text{loc}} = \{u = 0\}\), \(W^{ss}_{\text{loc}} = \{u_{1} = v_{1} = v_{2} = 0\}\) and \(W^{uu}_{\text{loc}} = \{u_{1} = u_{2} = v_{1} = 0\}\). For the particular case of normal form (\ref{eq23000}), the local extended stable and extended unstable invariant manifolds of \(O\) are straightened too, i.e. \(W^{sE}_{\text{loc}} = \{v_{2} = 0\}\) and \(W^{uE}_{\text{loc}} = \{u_{2} = 0\}\).

\subsection{Trajectories near the equilibrium state \texorpdfstring{\(O\)}{Lg}}\label{trajectoriesnearO}

In this section, we estimate the solutions of systems (\ref{eq74wt366cuw6gt5uw6v01092}), (\ref{eq19000}) and (\ref{eq23000}) near the equilibrium state \(O\) by using the technique of successive approximations.

Consider the system
\begin{equation}\label{eqbwiebro4ibyw43yyyiw0290}
\begin{aligned}
\dot{u}_{i} &= -\lambda_{i}u_{i} + F_{i}\left(u_{1}, u_{2}, v_{1}, v_{2}\right)\\
\dot{v}_{i} &= +\lambda_{i}v_{i} + G_{i}\left(u_{1}, u_{2}, v_{1}, v_{2}\right)
\end{aligned}
,\qquad (i=1,2)
\end{equation}
where \(F_{1}\), \(F_{2}\), \(G_{1}\) and \(G_{2}\) and their first derivatives vanish at the origin. By \cite{Dimabook} (Theorem 2.9), for given \(\tau \geq 0\) and sufficiently small \(u_{10}\), \(u_{20}\), \(v_{1\tau}\) and \(v_{2\tau}\) there exists a unique solution \(\left(u_{1}\left(t\right), u_{2}\left(t\right), v_{1}\left(t\right), v_{2}\left(t\right)\right)\) of system (\ref{eqbwiebro4ibyw43yyyiw0290}) such that
\begin{equation}\label{eq800020}
u_{1}\left(0\right) = u_{10}, \quad u_{2}\left(0\right) = u_{20}, \quad v_{1}\left(\tau\right) = v_{1\tau}, \quad v_{2}\left(\tau\right) = v_{2\tau}.
\end{equation}
The dependence of this solution on each of the variables \(\tau\), \(u_{10}\), \(u_{20}\), \(v_{1\tau}\) and \(v_{2\tau}\) is as smooth as the original system (\ref{eqbwiebro4ibyw43yyyiw0290}).

The following lemmas estimate the solutions of systems (\ref{eq74wt366cuw6gt5uw6v01092}), (\ref{eq19000}) and (\ref{eq23000}) that satisfy boundary condition (\ref{eq800020}). We prove these lemmas in Appendix \ref{ProofsOfBoundaryValueProblems}.

\begin{mylem}\label{flowlemmaresonant}
Let \(\lambda = \lambda_{1} = \lambda_{2}\). There exists \(M > 0\) such that for any sufficiently small \(\delta > 0\), and any \(u_{10}\), \(u_{20}\), \(v_{1\tau}\) and \(v_{2\tau}\), where \(\max \lbrace \lvert u_{10}\rvert, \lvert u_{20}\rvert, \lvert v_{1\tau}\rvert, \lvert v_{2\tau}\rvert\rbrace \leq \delta\), the solution \(\left(u\left(t\right), v\left(t\right)\right)\) of system (\ref{eq74wt366cuw6gt5uw6v01092}) that satisfies boundary condition (\ref{eq800020}) can be written as
\begin{equation}\label{eq7bwwkauser6518092hz9ub}
\begin{aligned}
u_{1}(t) =& e^{-\lambda t}u_{10} + \xi_{1}\left(x\right),\qquad
&& u_{2}(t) = e^{-\lambda t}u_{20} + \xi_{2}\left(x\right),\\
v_{1}(t) =& e^{-\lambda\left(\tau - t\right)} v_{1\tau} + \zeta_{1}\left(x\right),\qquad
&& v_{2}(t) = e^{-\lambda\left(\tau-t\right)} v_{2\tau} +\zeta_{2}\left(x\right),
\end{aligned}
\end{equation}
where \(x = \left(t, \tau, u_{10}, u_{20}, v_{1\tau}, v_{2\tau}\right)\), \(t\in\left[0, \tau\right]\), \(\max\lbrace\lvert \xi_{1}\rvert, \lvert \xi_{2}\rvert\rbrace \leq M e^{-\lambda t}\delta^{2}\) and\\ \(\max\lbrace\lvert \zeta_{1}\rvert, \lvert \zeta_{2}\rvert\rbrace \leq M e^{-\lambda\left(\tau - t\right)}\delta^{2}\). We can also write
\begin{equation}\label{eq89jjnmljdu09o3mjd874njsh}
\begin{aligned}
u_{1}(t) =& e^{-\lambda t}\left[u_{10} + O\left(\delta^{2}\right)\right],\qquad
&& u_{2}(t) = e^{-\lambda t}\left[u_{20} + O\left(\delta^{2}\right)\right],\\
v_{1}(t) =& e^{-\lambda\left(\tau - t\right)} \left[v_{1\tau} + O\left(\delta^{2}\right)\right],\qquad
&& v_{2}(t) = e^{-\lambda\left(\tau-t\right)}\left[v_{2\tau} + O\left(\delta^{2}\right)\right].
\end{aligned}
\end{equation}
\end{mylem}

\begin{mylem}\label{flowlemma}
There exists \(M > 0\) such that for any sufficiently small \(\delta > 0\), and any \(u_{10}\), \(u_{20}\), \(v_{1\tau}\) and \(v_{2\tau}\), where \(\max \lbrace \lvert u_{10}\rvert, \lvert u_{20}\rvert, \lvert v_{1\tau}\rvert, \lvert v_{2\tau}\rvert\rbrace \leq \delta\), the solution \(\left(u\left(t\right), v\left(t\right)\right)\) of system (\ref{eq19000}) that satisfies boundary condition (\ref{eq800020}) can be written as
\begin{equation}\label{eq67900}
\begin{aligned}
u_{1}(t) =& e^{-\lambda_{1} t}u_{10} + \xi_{1}\left(x\right),\qquad
&& u_{2}(t) = e^{-\lambda_{2} t}u_{20} + \xi_{2}\left(x\right),\\
v_{1}(t) =& e^{-\lambda_{1}\left(\tau - t\right)} v_{1\tau} + \zeta_{1}\left(x\right),\qquad
&& v_{2}(t) = e^{-\lambda_{2}\left(\tau-t\right)} v_{2\tau} +\zeta_{2}\left(x\right),
\end{aligned}
\end{equation}
where \(x = \left(t, \tau, u_{10}, u_{20}, v_{1\tau}, v_{2\tau}\right)\), \(t\in\left[0, \tau\right]\), and
\begin{equation*}
\begin{aligned}
\lvert \xi_{1} \rvert \leq & M\left[e^{-\lambda_{1}t}\delta \lvert u_{10}\rvert + e^{-\lambda_{1}\left(\tau - t\right) - \lambda_{2}t} \delta\lvert v_{1\tau}\rvert\right],
&& \lvert \xi_{2} \rvert \leq M e^{-\lambda_{2}t}\delta^{2},\\
\lvert \zeta_{1} \rvert \leq & M \left[e^{-\lambda_{1}\left(\tau - t\right)} \delta \lvert v_{1\tau}\rvert + e^{-\lambda_{2}\left(\tau - t\right) - \lambda_{1}t} \delta\lvert u_{10}\rvert\right],\qquad
&& \lvert \zeta_{2} \rvert \leq M e^{-\lambda_{2}\left(\tau-t\right)}\delta^{2}.
\end{aligned}
\end{equation*}
We can also write
\begin{equation}\label{eqiomljd784987qpqpla874}
\begin{aligned}
u_{1}(t) =& e^{-\lambda_{1}t}u_{10}\left[1 + O\left(\delta\right)\right] + e^{-\lambda_{1}\left(\tau - t\right) - \lambda_{2}t} O\left(\delta v_{1\tau}\right),\quad
&& u_{2}(t) = e^{-\lambda_{2}t}\left[u_{20} + O\left(\delta^{2}\right)\right],\\
v_{1}(t) =& e^{-\lambda_{1}\left(\tau - t\right)} v_{1\tau} \left[1 + O\left(\delta\right)\right] + e^{-\lambda_{2}\left(\tau - t\right) - \lambda_{1}t} O\left(\delta u_{10}\right),\quad
&& v_{2}(t) = e^{-\lambda_{2}\left(\tau-t\right)}\left[v_{2\tau} + O\left(\delta^{2}\right)\right].
\end{aligned}
\end{equation}
\end{mylem}

\begin{mylem}\label{flowlemma0}
There exists \(M > 0\) such that for any sufficiently small \(\delta > 0\), and any \(u_{10}\), \(u_{20}\), \(v_{1\tau}\) and \(v_{2\tau}\), where \(\max \lbrace \lvert u_{10}\rvert, \lvert u_{20}\rvert, \lvert v_{1\tau}\rvert, \lvert v_{2\tau}\rvert\rbrace \leq \delta\), the solution \(\left(u\left(t\right), v\left(t\right)\right)\) of system (\ref{eq23000}) that satisfies boundary condition (\ref{eq800020}) can be written in the form (\ref{eq67900}), where \(t\in\left[0, \tau\right]\) and
\begin{equation*}
\begin{aligned}
\lvert \xi_{1} \rvert \leq & M\left[e^{-\lambda_{1}t}\delta \lvert u_{10}\rvert + e^{-\lambda_{1}\left(\tau + t\right)} \delta\lvert v_{1\tau}\rvert\right],
&& \lvert \xi_{2} \rvert \leq M e^{-\lambda_{2}t}\delta^{2},\\
\lvert \zeta_{1} \rvert \leq & M \left[e^{-\lambda_{1}\left(\tau - t\right)} \delta \lvert v_{1\tau}\rvert + e^{-\lambda_{1}\left(2\tau + t\right)} \delta\lvert u_{10}\rvert\right],\qquad
&& \lvert \zeta_{2} \rvert \leq M e^{-\lambda_{2}\left(\tau-t\right)}\delta^{2}.
\end{aligned}
\end{equation*}
\end{mylem}

\begin{myrem}
For simplicity, we can write the solution given by Lemma \ref{flowlemma0} as
\begin{equation}\label{equos83km3odnin8b83xbv}
\begin{aligned}
u_{1}\left(t\right) &= e^{-\lambda_{1}t}u_{10}\left[1 + O\left(\delta\right)\right] + e^{-\lambda_{1}\left(\tau + t\right)} O\left(\delta v_{1\tau}\right),\quad
&& u_{2}\left(t\right) = e^{-\lambda_{2}t}\left[u_{20} + O\left(\delta^{2}\right)\right],\\
v_{1}\left(t\right) &= e^{-\lambda_{1}\left(\tau-t\right)}v_{1\tau}\left[1 + O\left(\delta\right)\right] + e^{-\lambda_{1}\left(2\tau - t\right)}O\left(\delta u_{10}\right),\quad
&& v_{2}\left(t\right) = e^{-\lambda_{2}\left(\tau-t\right)}\left[v_{2\tau} + O\left(\delta^{2}\right)\right].
\end{aligned}
\end{equation}
\end{myrem}

\subsection{Local maps and their properties}\label{localmapsection}

In this section, we use the results of the previous two sections to study the local maps for each of systems (\ref{eq74wt366cuw6gt5uw6v01092}), (\ref{eq19000}) and (\ref{eq23000}). Recall (\ref{localmap}) and write
\begin{equation}\label{eq68960}
T^{\text{loc}}\left(u_{10}, v_{10}\right) = \left(u_{1\tau}, v_{1\tau}\right) = \big(\eta_{1}\left(u_{10},v_{10}\right), \eta_{2}\left(u_{10},v_{10}\right)\big),
\end{equation}
where \(\eta_{1}\) and \(\eta_{2}\) are some functions. In the previous section, for each of systems (\ref{eq74wt366cuw6gt5uw6v01092}), (\ref{eq19000}) and (\ref{eq23000}), we have approximated the unique solution \(\left(u^{*}, v^{*}\right)\) which satisfies boundary conditions (\ref{eq800020}) (see Lemmas \ref{flowlemmaresonant}, \ref{flowlemma} and \ref{flowlemma0}). We write this solution as
\begin{equation}\label{eq90ioi48iskju389874jksh}
\begin{aligned}
u_{1}^{*}\left(t\right) = u_{1}^{*}\left(t, \tau, u_{10}, u_{20}, v_{1\tau}, v_{2\tau}\right),\\
u_{2}^{*}\left(t\right) = u_{2}^{*}\left(t, \tau, u_{10}, u_{20}, v_{1\tau}, v_{2\tau}\right),\\
v_{1}^{*}\left(t\right) = v_{1}^{*}\left(t, \tau, u_{10}, u_{20}, v_{1\tau}, v_{2\tau}\right),\\
v_{2}^{*}\left(t\right) = v_{2}^{*}\left(t, \tau, u_{10}, u_{20}, v_{1\tau}, v_{2\tau}\right),
\end{aligned}
\end{equation}
to emphasize that it explicitly depends on \(t\), \(\tau\), \(u_{10}\), \(u_{20}\), \(v_{1\tau}\) and \(v_{2\tau}\). This solution represents an orbit which at \(t=0\) is at the point \(\left(u_{10}, u_{20}, v_{10}, v_{20}\right)\) and at \(t=\tau\) is at the point \(\left(u_{1\tau}, u_{2\tau}, v_{1\tau}, v_{2\tau}\right)\).

To study the map \(T^{\text{loc}}\), we consider the case in which \(u_{20} = v_{2\tau} = \delta\), i.e. the points \(\left(u_{10}, u_{20}, v_{10}, v_{20}\right)\)  and \(\left(u_{1\tau}, u_{2\tau}, v_{1\tau}, v_{2\tau}\right)\) belong to \(\Pi^{s}\) and \(\Pi^{u}\), respectively. Evaluating the first equation of  (\ref{eq90ioi48iskju389874jksh}) at \(t = \tau\) and the last two equations of (\ref{eq90ioi48iskju389874jksh}) at \(t=0\) gives
\begin{equation}\label{eqin98nwkndoib112}
\begin{aligned}
u_{1\tau} &= u_{1}^{*}\left(\tau, \tau, u_{10}, \delta, v_{1\tau}, \delta\right),\\
v_{10} &= v_{1}^{*}\left(0, \tau, u_{10}, \delta, v_{1\tau}, \delta\right),\\
v_{20} &= v_{2}^{*}\left(0, \tau, u_{10}, \delta, v_{1\tau}, \delta\right),
\end{aligned}
\end{equation}
which is an implicit relation between \(u_{10}\), \(v_{10}\), \(v_{20}\), \(u_{1\tau} = \eta_{1}\left(u_{10},v_{10}\right)\), \(v_{1\tau} = \eta_{2}\left(u_{10},v_{10}\right)\) and \(\tau\). On the other hand, \(\tau\) and \(v_{20}\) can be expressed as functions of \(\left(u_{10}, v_{10}\right)\). This allows us to approximate the functions \(\eta_{1}\) and \(\eta_{2}\).

\begin{mynotation}
Hereafter, we use the following notation: \(\gamma = \frac{\lambda_{1}}{\lambda_{2}}\).
\end{mynotation}

\subsubsection{Choice of coordinates on the cross-sections}\label{Coordinates-on-cross-sections}

We point out here that we can always choose \(\left(u_{1}, v_{1}\right)\)- coordinate system on each of our cross-sections:

\begin{mylem}\label{Cor518b8b89or48v10010}
In each of the cases \(\lambda_{1} = \lambda_{2}\), \(\lambda_{2} < 2\lambda_{1}\) and \(2\lambda_{1} < \lambda_{2}\), for any arbitrary point \(\left(u_{1}, u_{2}, v_{1}, v_{2}\right)\) of each of the cross-sections \(\Pi^{s}\), \(\Pi^{u}\), \(\Pi^{s}_{1}\), \(\Pi^{u}_{1}\), \(\Pi^{s}_{2}\), \(\Pi^{u}_{2}\) and \(\Sigma\), the variables \(u_{2}\) and \(v_{2}\) are uniquely determined by \(\left(u_{1}, v_{1}\right)\).
\end{mylem}
\begin{proof}
We only prove the statement for \(\Pi^{s}\) and \(\Pi^{u}\). The proof for the other cross-sections is the same.

Consider system (\ref{eq74wt366cuw6gt5uw6v01092}) and suppose \(\lambda_{1} = \lambda_{2}\). Let \(\left(u_{10}, \delta, v_{10}, v_{20}\right)\) and \(\left(u_{1\tau}, u_{2\tau}, v_{1\tau}, \delta\right)\) be two points on \(\Pi^{s}\) and \(\Pi^{u}\), respectively. By virtue of the relation \(\{H = 0\}\), where \(H\) is as in (\ref{eq98n9uo8wbr8bau4b819ub92q1u8}), we have
\begin{equation}\label{eq7901oojjib9iur8b}
v_{20} = {\delta}^{-1} u_{10}v_{10}\quad \text{and} \quad u_{2\tau} = {\delta}^{-1} u_{1\tau}v_{1\tau}.
\end{equation}
This proves the lemma for the case \(\lambda_{1} = \lambda_{2}\).

A straightforward calculation (see \cite{BakraniPhDthesis}) shows that for the cases of systems (\ref{eq19000}) and (\ref{eq23000}), and their corresponding first integrals, we have \(H_{v_{2}} \left(0, \delta, 0, 0\right) \neq 0\) and \(H_{u_{2}} \left(0, 0, 0, \delta\right) \neq 0\).
Then, the proof for the cases \(\lambda_{2} < 2\lambda_{1}\) and \(2\lambda_{1} < \lambda_{2}\) follows from the implicit function theorem.
\end{proof}

\subsubsection{Local maps: case \(\lambda_{1} = \lambda_{2}\)}

We prove \(\mathcal{D} = \emptyset\) by showing that (\ref{conditionfordomainD}) never holds. This implies that the Poincar\'e map along \(\Gamma\) cannot be defined when \(\lambda_{1} = \lambda_{2}\). This also proves Theorem \ref{thmkuyvuyuy090988} for the particular case of \(\lambda_{1} = \lambda_{2}\).

Let \(\lambda = \lambda_{1} = \lambda_{2}\) and consider the case \(u_{20} =  v_{2\tau} = \delta\). Evaluating the first two equations of (\ref{eq89jjnmljdu09o3mjd874njsh}) at \(t=\tau\) and the last two equations at \(t=0\) gives
\begin{equation}\label{equislmkdhyuromns09j39n}
\begin{aligned}
u_{1\tau} =& e^{-\lambda \tau}\left[u_{10} + O\left(\delta^{2}\right)\right],\qquad
&& u_{2\tau} = e^{-\lambda \tau}\delta\left[1 + O\left(\delta\right)\right],\\
v_{10} =& e^{-\lambda\tau} \left[v_{1\tau} + O\left(\delta^{2}\right)\right],\qquad
&& v_{20} = e^{-\lambda\tau}\delta\left[1 + O\left(\delta\right)\right].
\end{aligned}
\end{equation}
Substituting (\ref{eq7901oojjib9iur8b}) into this relation gives \(e^{-\lambda\tau} = \frac{u_{10}v_{10}}{\delta^{2}}\left[1 + O\left(\delta\right)\right]\). Therefore,
\begin{equation*}
v_{1\tau} = e^{\lambda\tau} v_{10} + O\left(\delta^{2}\right) = \delta^{2}u_{10}^{-1}\left[1 + O\left(\delta\right)\right] + O\left(\delta^{2}\right) = u_{10}^{-1}\delta^{2}\left[1 + O\left(\delta\right)\right].
\end{equation*}
For a given sufficiently small \(\delta\), we have
\begin{equation*}
\lim_{u_{10} \rightarrow 0}\lVert \left(u_{1 \tau}, v_{1\tau}\right)\rVert \geq \lim_{u_{10} \rightarrow 0}\lvert v_{1\tau}\rvert = \lim_{u_{10} \rightarrow 0} \lvert u_{10}\rvert^{-1}\delta^{2}\left[1 + O\left(\delta\right)\right] = \infty.
\end{equation*}
This means that (\ref{conditionfordomainD}) does not hold when \(\epsilon\) and \(\epsilon_{u}\) are chosen sufficiently small. On the other hand, it is easily seen that the same happens for the points \(\left(u_{10}, v_{10}\right)\) in \(\mathbb{D}\). The same also holds for the case of homoclinic figure-eight. Therefore,
\begin{myprop}\label{Cor7uijnqnjk3o49inh44}
When \(\lambda_{1} = \lambda_{2}\), we have \(\mathcal{D} = \mathbb{D} = \mathcal{D}^{1} = \mathbb{D}^{1} = \mathcal{D}^{2} = \mathbb{D}^{2} = \emptyset\).
\end{myprop}

\subsubsection{Local maps: case \(\lambda_{1} < \lambda_{2} < 2\lambda_{1}\)}

Let \(\lambda < \lambda_{2} < 2\lambda_{1}\) and consider the case \(u_{20} =  v_{2\tau} = \delta\). Evaluating the first two equations of (\ref{eqiomljd784987qpqpla874}) at \(t=\tau\) and the last two equations at \(t=0\) gives
\begin{equation}\label{eq6739736789927600657}
\begin{aligned}
u_{1\tau} =& e^{-\lambda_{1}\tau}u_{10}\left[1 + O\left(\delta\right)\right] + e^{- \lambda_{2}\tau} O\left(\delta v_{1\tau}\right),\qquad
&& u_{2\tau} = e^{-\lambda_{2}\tau}\delta\left[1 + O\left(\delta\right)\right],\\
v_{10} =& e^{-\lambda_{1}\tau} v_{1\tau} \left[1 + O\left(\delta\right)\right] + e^{-\lambda_{2}\tau} O\left(\delta u_{10}\right),\qquad
&& v_{20} = e^{-\lambda_{2}\tau}\delta\left[1 + O\left(\delta\right)\right].
\end{aligned}
\end{equation}
This, in particular, implies
\begin{equation}\label{eqyusolkuoplkmh23ui8hd7837}
v_{1\tau} = e^{\lambda_{1}\tau} v_{10} \left[1 + O\left(\delta\right)\right] + e^{\left(\lambda_{1}-\lambda_{2}\right)\tau} O\left(\delta u_{10}\right).
\end{equation}
First integral (\ref{eq22000}) vanishes at \(\left(u_{10}, \delta, v_{10}, v_{20}\right)\in\Pi^{s}\). Thus \(v_{20} = \frac{\gamma}{\delta}\cdot u_{10}v_{10}\left[1+o(1)\right]\). Therefore, (\ref{eq6739736789927600657}) implies
\begin{equation}\label{eq703okklcmnajue55mcjlkdi093kdi}
e^{-\lambda_{2}\tau} = \gamma\delta^{-2} u_{10}v_{10}\left[1 + O\left(\delta\right)\right],
\end{equation}
and therefore
\begin{equation}\label{eq8isolpaloue74jndu}
e^{-\lambda_{1}\tau} = \left(\gamma\delta^{-2} u_{10}v_{10}\right)^{\gamma}\left[1 + O\left(\delta\right)\right].
\end{equation}
By these relations, we rewrite (\ref{eqyusolkuoplkmh23ui8hd7837}) as
\begin{equation}\label{eq6kxnjaalkapqp3891000qa}
\eta_{2}\left(u_{10}, v_{10}\right) = v_{1\tau} = e^{\lambda_{1}\tau} v_{10} \left[1 + O\left(\delta\right)\right].
\end{equation}
Substituting this into the equation of \(u_{1\tau}\) in (\ref{eq6739736789927600657}) gives
\begin{equation}\label{eq888u0aplanaki}
\eta_{1}\left(u_{10}, v_{10}\right) = u_{1\tau} = e^{-\lambda_{1}\tau}u_{10}\left[1 + O\left(\delta\right)\right] + e^{\left(\lambda_{1} - \lambda_{2}\right)\tau} O\left(\delta v_{10}\right).
\end{equation}

\begin{figure}
\centering
\includegraphics[scale=.19]{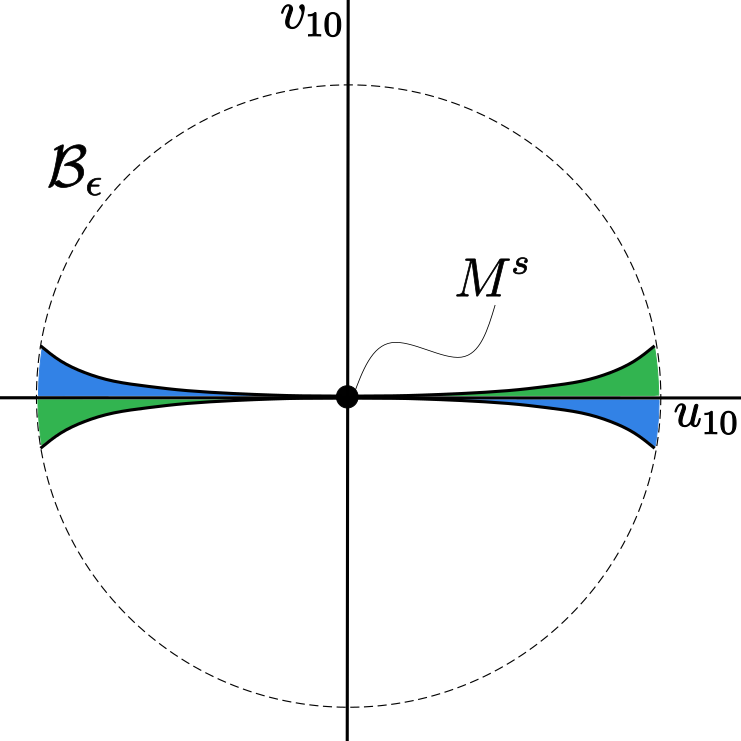}
\caption{\small The regions \(\mathcal{D}\) and \(\mathbb{D}\) for the case \(\lambda_{1} < \lambda_{2} < 2\lambda_{1}\) are shown in green and blue, respectively. They are surrounded by horizontal axis, \(\epsilon\)-ball \(\mathcal{B}_{\epsilon}\) and the curves characterized by (\ref{equios093mljoiudmkl28plso}). Since \(\gamma = \lambda_{1} {\lambda_{2}}^{-1} > 0.5\), these curves are tangent to the horizontal axis at \(M^s = \left(0,0\right)\).}
\label{Figure99nubdqvetc6aaa}
\end{figure}

Let us now explore the domain \(\mathcal{D}\) of the map \(T^{\text{loc}}\).
By choosing \(\delta\) sufficiently small such that \(\lvert O\left(\delta\right) \rvert \leq 1\), we have
\begin{equation*}
\Big\vert \frac{u_{1\tau}}{v_{1\tau}}\Big\vert \leq 2 e^{-\lambda_{2}\tau}\frac{\lvert u_{10}\rvert}{\lvert v_{10}\rvert} + e^{- \lambda_{2}\tau}
\leq \frac{4\gamma}{\delta^{2}}\left(u_{10}^{2} + \lvert u_{10} v_{10}\rvert\right) \leq \frac{8\gamma}{\delta^{2}} \epsilon^{2}.
\end{equation*}
Therefore, for any given (fixed) sufficiently small \(\delta\), we have \(u_{1\tau} = v_{1\tau}O\left(\epsilon^{2}\right)\). Thus,
\begin{equation*}
\begin{aligned}
\left\|\left(u_{1\tau}, v_{1\tau}\right)\right\| = \sqrt{u_{1\tau}^{2} + v_{1\tau}^{2}} = \lvert v_{1\tau}\rvert \left[1 + O\left(\epsilon^{2}\right)\right] = \left(\gamma\delta^{-2}\right)^{-\gamma} \lvert u_{10}\rvert^{-\gamma} \lvert v_{10}\rvert^{1-\gamma} \left[1 + O\left(\delta\right)\right].
\end{aligned}
\end{equation*}
Therefore, \(\lVert \left(u_{1\tau}, v_{1\tau}\right)\rVert < \epsilon_{u}\) if and only if
\begin{equation}\label{equios093mljoiudmkl28plso}
\lvert v_{10}\rvert < {\epsilon_{u}}^{\frac{1}{1-\gamma}} \left(\gamma\delta^{-2}\right)^{\frac{\gamma}{1-\gamma}} \lvert u_{10}\rvert^{\frac{\gamma}{1-\gamma}} \left[1 + O\left(\delta\right)\right],\qquad (\gamma = \lambda_{1}\lambda_{2}^{-1} > 0.5).
\end{equation}
By virtue of (\ref{eq703okklcmnajue55mcjlkdi093kdi}), we see that if \(\left(u_{10}, v_{10}\right) \in \mathcal{D}\), then \(u_{10}v_{10}>0\). It is also easy to see that analogous statements hold for the points in \(\mathbb{D}\). This gives:
\begin{myprop}\label{Cor8i12sehu8901fdg}
Let \(\lambda < \lambda_{2} < 2\lambda_{1}\). For a given sufficiently small \(\delta\), we can choose \(\epsilon\) and \(\epsilon_{u}\) so that the domain \(\mathcal{D}\) (resp. \(\mathbb{D}\)) becomes the set of all points \(\left(u_{10}, v_{10}\right)\) in \(\Pi^{s}\) such that \(0 < u_{10}v_{10}\) (resp. \(u_{10}v_{10} < 0\)), \(\lVert \left(u_{10}, v_{10}\right)\rVert < \epsilon\) and (\ref{equios093mljoiudmkl28plso}) holds (see Figure \ref{Figure99nubdqvetc6aaa}).
\end{myprop}

\begin{myrem}\label{Rem89i9ub8q7v76ce365}
The case of homoclinic figure-eight is the same. Relation \(u_{1\tau} = v_{1\tau}O\left(\epsilon^{2}\right)\) holds for any \(\left(u_{10}, v_{10}\right)\) on \(\Pi^{s}_{i}\), and the domains \(\mathcal{D}^{i}\) and \(\mathbb{D}^{i}\) are given by Proposition \ref{Cor8i12sehu8901fdg} (\(i=1,2\)).
\end{myrem}

\subsubsection{Local maps: case \(2\lambda_{1} < \lambda_{2}\)}\label{Domainjyagvxjycefjthc2tc65}

Let \(2\lambda_{1} < \lambda_{2}\). Evaluating the first two equations of (\ref{equos83km3odnin8b83xbv}) at \(t=\tau\) and the last two equations at \(t=0\) gives
\begin{equation}\label{equilod782plsmk9id9ij}
\begin{aligned}
u_{1\tau} &= e^{-\lambda_{1}\tau}u_{10}\left[1 + O\left(\delta\right)\right] + e^{-2\lambda_{1}\tau} O\left(\delta v_{1\tau}\right),\quad
&& u_{2\tau} = e^{-\lambda_{2}\tau}\left[u_{20} + O\left(\delta^{2}\right)\right],\\
v_{10} &= e^{-\lambda_{1}\tau}v_{1\tau}\left[1 + O\left(\delta\right)\right] + e^{-2\lambda_{1}\tau}O\left(\delta u_{10}\right),\quad
&& v_{20} = e^{-\lambda_{2}\tau}\left[v_{2\tau} + O\left(\delta^{2}\right)\right].
\end{aligned}
\end{equation}
For the particular case of \(u_{20} =  v_{2\tau} = \delta\), we have
\begin{equation}\label{eqi903lokdjkflk}
\begin{aligned}
u_{1\tau} =& e^{-\lambda_{1}\tau}u_{10}\left[1 + O\left(\delta\right)\right] + e^{- 2\lambda_{1}\tau} O\left(\delta\lvert v_{1\tau}\rvert\right),\qquad
&& u_{2\tau} = e^{-\lambda_{2}\tau}\delta\left[1 + O\left(\delta\right)\right],\\
v_{10} =& e^{-\lambda_{1}\tau} v_{1\tau} \left[1 + O\left(\delta\right)\right] + e^{-2\lambda_{1}\tau} O\left(\delta\lvert u_{10}\rvert\right),\qquad
&& v_{20} = e^{-\lambda_{2}\tau}\delta\left[1 + O\left(\delta\right)\right].
\end{aligned}
\end{equation}
This, in particular, implies
\begin{equation}\label{eq09o02oknpqlamncccfg}
\eta_{2}\left(u_{10}, v_{10}\right) = v_{1\tau} = e^{\lambda_{1}\tau} v_{10} \left[1 + O\left(\delta\right)\right] + e^{-\lambda_{1}\tau} O\left(\delta\lvert u_{10}\rvert\right).
\end{equation}
Substituting this into the equation of \(u_{1\tau}\) in (\ref{eqi903lokdjkflk}) gives
\begin{equation}\label{equ9sjhe9opqk39kxnj}
\eta_{1}\left(u_{10}, v_{10}\right) = u_{1\tau} = e^{-\lambda_{1}\tau}u_{10}\left[1 + O\left(\delta\right)\right] + e^{-\lambda_{1}\tau} O\left(\delta\lvert v_{10}\rvert\right).
\end{equation}
Then, local map (\ref{eq68960}) maps \(\left(u_{10}, v_{10}\right)\) to \(\left(u_{1\tau}, v_{1\tau}\right)\), where \(u_{1\tau}\) and \(v_{1\tau}\) are as in (\ref{equ9sjhe9opqk39kxnj}) and (\ref{eq09o02oknpqlamncccfg}), respectively, and \(\tau\) is a function of \(\left(u_{10}, v_{10}\right)\). It is not as straightforward as the previous two cases to express \(\tau\) as a function of \(\left(u_{10}, v_{10}\right)\). This is not straightforward either to find the domain \(\mathcal{D}\) of \(T^{\text{loc}}\). Below, we divide \(\mathcal{D}\) into three regions (it is shown that \(\mathcal{D} \neq \emptyset\)) and study each region separately.

\begin{figure}
\centering
\begin{subfigure}{0.4\textwidth}
\centering
\includegraphics[scale=0.18]{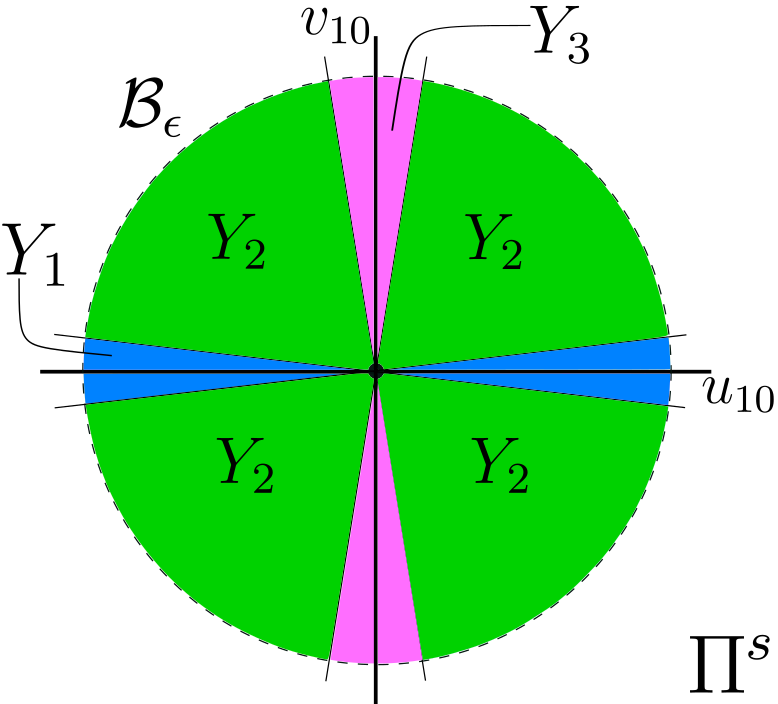}
\end{subfigure}
%\hfill
\begin{subfigure}{0.4\textwidth}
\centering
\includegraphics[scale=0.18]{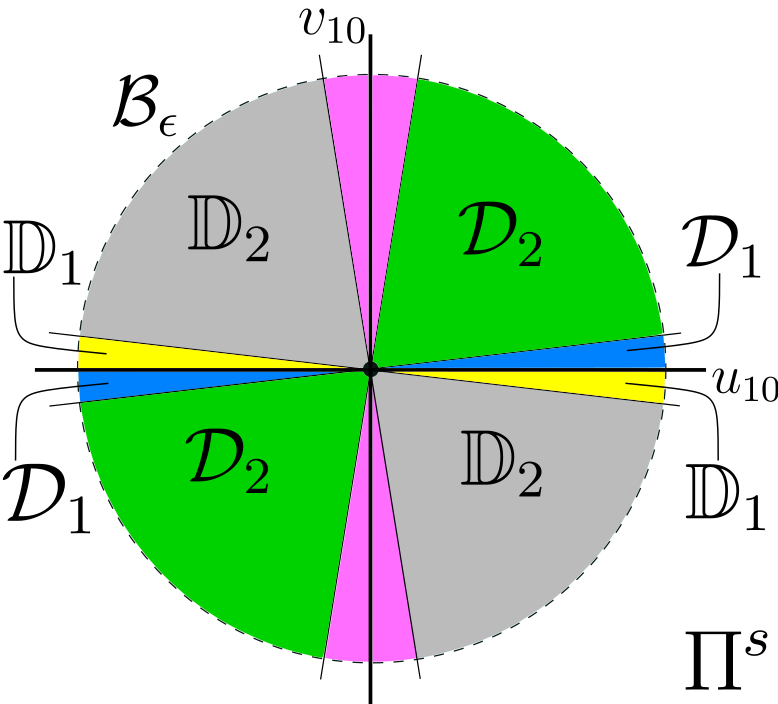}
\end{subfigure}
\caption{\small We divide the \(\epsilon\)-ball in \(\Pi^{s}\) centered at \(M^{s}\) into three disjoint regions: \(Y_{1}\) (blue), \(Y_{2}\) (green) and \(Y_{3}\) (pink), shown in the left figure. To investigate the sets \(\mathcal{D}\) and \(\mathbb{D}\) when \(\lambda_{2} > 2\lambda_{1}\), we consider the intersection of each of these sets with the regions \(Y_{1}\), \(Y_{2}\) and \(Y_{3}\). We then define \(\mathcal{D}_{i} = \mathcal{D} \cap Y_{i}\) and \(\mathbb{D}_{i} = \mathbb{D} \cap Y_{i}\). The regions \(\mathcal{D}_{1}\), \(\mathcal{D}_{2}\), \(\mathbb{D}_{1}\) and \(\mathbb{D}_{2}\) are shown by blue, green, yellow and gray, respectively in the right figure. The sets \(\mathcal{D}_{3}\) and \(\mathbb{D}_{3}\) are subsets of the pink region.}
\label{Fig7uib9uoubyvdaetv}
\end{figure}

Let \(\mathcal{B}_{\epsilon}\) be the \(\epsilon\)-ball in \(\Pi^{s}\) centered at \(M^{s}\). For a given \(m>1\) define
\begin{equation}\label{eqyujnj68999jhhhgg34}
\begin{gathered}
Y^{m}_{1} = \left\{\left(u_{10}, v_{10}\right)\in \mathcal{B}_{\epsilon}: \lvert v_{10}\rvert < m^{-1}\lvert u_{10}\rvert\right\},\\
Y^{m}_{2} = \left\{\left(u_{10}, v_{10}\right)\in \mathcal{B}_{\epsilon}: m^{-1}\lvert u_{10}\rvert \leq \lvert v_{10}\rvert \leq m\lvert u_{10}\rvert\right\}\\
Y^{m}_{3} = \left\{\left(u_{10}, v_{10}\right)\in \mathcal{B}_{\epsilon}: m\lvert u_{10}\rvert < \lvert v_{10}\rvert \right\}
\end{gathered}
\end{equation}
(see Figure \ref{Fig7uib9uoubyvdaetv}). Obviously, \(\mathcal{B}_{\epsilon} = Y^{m}_{1} \cup Y^{m}_{2} \cup Y^{m}_{3}\). We define
\begin{equation}\label{eqiuvyseyt74vqi8l3viul181}
\mathcal{D}_{1}^{\epsilon}:= \mathcal{D}\cap Y^{m}_{1},\qquad \mathcal{D}_{2}^{\epsilon}:= \mathcal{D}\cap Y^{m}_{2}, \qquad \mathcal{D}_{3}^{\epsilon}:= \mathcal{D}\cap Y^{m}_{3}.
\end{equation}
Analogously, we define \(\mathbb{D}_{1}^{\epsilon}:= \mathbb{D}\cap Y^{m}_{1}\), \(\mathbb{D}_{2}^{\epsilon}:= \mathbb{D}\cap Y^{m}_{2}\) and \(\mathbb{D}_{3}^{\epsilon}:= \mathbb{D}\cap Y^{m}_{3}\). We may drop the subscript \(\epsilon\) and \(m\), when no confusion arises.

For \(\left(u_{10}, v_{10}\right)\in Y^{m}_{1} \cup Y^{m}_{2}\), we have \(\lvert v_{10}\rvert \leq m\lvert u_{10}\rvert\) and therefore \(v_{10} = O\left(u_{10}\right)\). By virtue of this relation and taking into account that first integral (\ref{eqhhuhut5989b8b8y4qlz37y7t}) vanishes at \(\left(u_{10}, \delta, v_{10}, v_{20}\right)\in \Pi^{s}\), we derive
\begin{equation}\label{eq09uis4kkn0cusbgdj}
v_{20} = \gamma\delta^{-1} u_{10}v_{10}\left[1 + O\left(\delta\right)\right].
\end{equation}
This relation together with (\ref{eqi903lokdjkflk}) implies that any point \(\left(u_{10}, v_{10}\right)\in Y^{m}_{1} \cup Y^{m}_{2}\) reaches \(\Pi^{u}\) if \(u_{10}v_{10} > 0\), and reaches \(\Sigma\) if \(u_{10}v_{10} < 0\). Therefore, to find \(\mathcal{D}_{1}\cup \mathcal{D}_{2}\) (\(\mathbb{D}_{1}\cup \mathbb{D}_{2}\)), it is sufficient to find the points in \(Y^{m}_{1} \cup Y^{m}_{2}\) for which \(\lVert \left(u_{1\tau}, v_{1\tau}\right)\rVert < \epsilon_{u}\).

Like the preceding two cases, relation (\ref{eq09uis4kkn0cusbgdj}) yields (\ref{eq703okklcmnajue55mcjlkdi093kdi}) and (\ref{eq8isolpaloue74jndu}). Let \(\delta\) be sufficiently small. With (\ref{eq09o02oknpqlamncccfg}), (\ref{equ9sjhe9opqk39kxnj}) and some straightforward calculation, we derive
\begin{equation*}
\lvert u_{1\tau}\rvert \leq \left(2 +  \delta\right) e^{-\lambda_{1}\tau}\epsilon,\qquad \lvert v_{1\tau}\rvert \leq \left[4 m^{1-\gamma} \left(\delta^{2}\gamma^{-1}\right)^{\gamma} + 1\right] \epsilon^{1-2\gamma}.
\end{equation*}
This gives the following:
\begin{myprop}
For given \(m\), sufficiently small \(\delta\) and sufficiently small \(\epsilon_{u}\), we can choose \(\epsilon\) sufficiently small such that for \(i= 1, 2\) we have
\begin{equation*}
\begin{aligned}
&\mathcal{D}_{i} = \lbrace \left(u_{10}, v_{10}\right)\in Y^{m}_{i},\quad u_{10}v_{10} > 0\rbrace,\quad
&&\mathbb{D}_{i} = \lbrace \left(u_{10}, v_{10}\right)\in Y^{m}_{i},\quad u_{10}v_{10} < 0\rbrace.
\end{aligned}
\end{equation*}
\end{myprop}

Now, consider \(\left(u_{10}, v_{10}\right)\in Y^{m}_{2}\cup Y^{m}_{3}\). We have \(\lvert u_{10}\rvert \leq m\lvert v_{10}\rvert\) and hence \(u_{10} = O\left(v_{10}\right)\). By virtue of this relation and relation (\ref{equilod782plsmk9id9ij}), we obtain
\begin{equation}\label{eqrtyuwikknj67234}
v_{10} = e^{-\lambda_{1}\tau} v_{1\tau} \left[1 + O\left(\delta\right)\right] \qquad \text{and}\qquad v_{1\tau} = e^{\lambda_{1}\tau} v_{10} \left[1 + O\left(\delta\right)\right].
\end{equation}
This relation together with (\ref{equilod782plsmk9id9ij}) gives
\begin{equation}\label{eq679oo89lkikmn44hujj}
\frac{u_{1\tau}}{v_{1\tau}} = \frac{e^{-\lambda_{1}\tau}u_{10}\left[1 + O\left(\delta\right)\right] + e^{-2\lambda_{1}\tau} O\left(\delta v_{1\tau}\right)}{v_{1\tau}} = e^{-2\lambda_{1}\tau}\left[\frac{u_{10}}{v_{10}} + O\left(\delta\right)\right] = o\left(1\right),
\end{equation}
which implies \(u_{1\tau} = o(v_{1\tau})\). Thus, when \(\left(u_{10}, v_{10}\right)\in Y^{m}_{2}\cup Y^{m}_{3}\), we have
\begin{equation*}
\big\Vert \left(u_{1\tau}, v_{1\tau}\right)\big\Vert = \sqrt{u_{1\tau}^{2} + v_{1\tau}^{2}} = \lvert v_{1\tau}\rvert \left[1 + o\left(1\right)\right].
\end{equation*}

Note that it was relation (\ref{eq09uis4kkn0cusbgdj}) that enabled us to, first, identify the points in \(Y^{m}_{i}\) for which \(\lVert \left(u_{1\tau}, v_{1\tau}\right)\rVert < \epsilon_{u}\) holds, and second, distinguish \(\mathcal{D}_{i}\) from \(\mathbb{D}_{i}\) for \(i=1,2\). For the case of \(i=3\), we cannot deduce such a relation from first integral (\ref{eqhhuhut5989b8b8y4qlz37y7t}). However, as we see later, the dynamics on \(Y^{m}_{3}\) is quite simple and can be analyzed without knowing \(\mathcal{D}_{3}\) and \(\mathbb{D}_{3}\) precisely.

Meanwhile, we have shown the following

\begin{myprop}\label{Cor67ujokjsuijwww}
If \(\left(u_{10}, v_{10}\right)\in Y^{m}_{2}\), then (\ref{eq09uis4kkn0cusbgdj}), (\ref{eq703okklcmnajue55mcjlkdi093kdi}), (\ref{eq8isolpaloue74jndu}), (\ref{eqrtyuwikknj67234}), \(u_{10} = O\left(v_{10}\right)\) and \(v_{10} = O\left(u_{10}\right)\) hold. If \((u_{10}, v_{10})\in Y^{m}_{3}\), then (\ref{eq679oo89lkikmn44hujj}) and \(u_{1\tau} = o(v_{1\tau})\) hold.
\end{myprop}

\section{Analysis near homoclinics and super-homoclinics}
\label{Dynamics nearhomoclinic}

The purpose of this section is to study the dynamics near (single and figure-eight) homoclinic and super-homoclinic orbits. In particular, we prove in this section, all the theorems stated in the Introduction. In the first section below, we introduce some concepts and notations. The second and the third sections are dedicated to study the dynamics near a single homoclinic orbit. We prove Theorems \ref{Th782ju2iu22}, \ref{thmkuyvuyuy090988} and \ref{Invariantmanifoldthm} in these two sections. The ideas and techniques which are used to prove these theorems are also used in the later sections. In the fourth section, we extend the results obtained for a single homoclinic to the case of the homoclinic figure-eight. The proofs of Theorems \ref{Thmyipoinoib53}, \ref{thm89bqyvrtyvtv} and \ref{Thm6892jjdibbea} are provided in this section. Finally, we study the case of a super-homoclinic and prove Theorems \ref{superhomoclinicthm} and \ref{superhomoclinicthmforfigure8} in the fifth (and the last) section.

\subsection{Set-up and notations}\label{Somenotationsandsetting}

Choose a sufficiently small \(\delta > 0\) such that all the statements of the previous sections hold. Fix this \(\delta\). According to (\ref{eq789o0o1ok2nnyxuj}) and (\ref{eq63000}), for \((u_{10}, v_{10}) \in \mathcal{D}\subset\Pi^{s}\), we have
\begin{equation}\label{eq78uikolp23456fga}
\left(
\begin{array}{c}
\overline{u}_{10}\\ \overline{v}_{10}
\end{array} \right) = T\left(
\begin{array}{c}
u_{10}\\ v_{10}
\end{array} \right) = \left(
\begin{array}{c}
\left[a + o\left(1\right)\right]\,u_{1\tau} + \left[b + o\left(1\right)\right]\,v_{1\tau}\\
\left[c + o\left(1\right)\right]\,u_{1\tau} + \left[d + o\left(1\right)\right]\,v_{1\tau}
\end{array} \right),
\end{equation}
where \(a\), \(b\), \(c\) and \(d\) are real constants (in fact, these coefficients are functions of \(\delta\) but since \(\delta\) is assumed to be fixed, we treat these coefficients as constants). Our job is to analyze this map for different values of \(a\), \(b\), \(c\) and \(d\), and for each of the cases \(\lambda_{2} < 2\lambda_{1}\) and \(2\lambda_{1} < \lambda_{2}\). In this strand, we first introduce some notations:
\begin{mynotation}\label{Notation6735oilboiyvkiutc}
Let \(\mathcal{N} \subset \mathcal{M}\) be two arbitrary sets and \(f: \mathcal{N} \rightarrow \mathcal{M}\) be an injective map. We denote the set of the points in \(\mathcal{N}\) whose forward orbits lie entirely in \(\mathcal{N}\) by \(\Lambda_{\mathcal{N}, f}^{s}\) or \(\Lambda_{\mathcal{N}}^{s}\), when no confusion arises. Indeed,
\begin{equation*}
\Lambda_{\mathcal{N}, f}^{s} = \Lambda_{\mathcal{N}}^{s} = \lbrace x\in \mathcal{N}: f^{n}\left(x\right) \in \mathcal{N}, \quad\forall n \geq 0\rbrace.
\end{equation*}
We denote the set of the points in \(\mathcal{N}\) whose backward orbits lie entirely in \(\mathcal{N}\) by \(\Lambda_{\mathcal{N}, f}^{u}\) or \(\Lambda_{\mathcal{N}}^{u}\), when no confusion arises. Indeed,
\begin{equation*}
\Lambda_{\mathcal{N}, f}^{u} = \Lambda_{\mathcal{N}}^{u} = \lbrace x\in \mathcal{N}: \text{ for all } n \geq 0,\quad f^{-n}\left(x\right) \text{ exists and belongs to } \mathcal{N}\rbrace.
\end{equation*}
\end{mynotation}
\begin{myrem}\label{Remiuboyboweouberuvg}
Recall \(W^{s}_{\mathcal{U}}\left(O\right)\) and \(W^{u}_{\mathcal{U}}\left(O\right)\) from the Introduction. Taking into account that \(W^{s}_{\mathcal{U}}\left(O\right) \cap \Lambda^{s}_{\mathcal{D}, T} = \emptyset\) and \(W^{u}_{\mathcal{U}}\left(O\right) \cap \Lambda^{u}_{\mathcal{D}, T} = \emptyset\), we can reformulate Theorem \ref{Th782ju2iu22} as follows: The forward (resp. backward) orbit of any point on \(\Lambda^{s}_{\mathcal{D}, T}\) (resp. \(\Lambda^{u}_{\mathcal{D}, T}\)) converges to the homoclinic orbit \(\Gamma\).
\end{myrem}
\begin{mynotation}\label{Not89kookkjweer}
Given a point \((u_{10}, v_{10})\) on a given cross-section, we denote the quantity \(\frac{v_{10}}{u_{10}}\) (when \(u_{10}\neq 0\)) by \(w(u_{10},v_{10})\) or \(w\). Consider the case \((u_{10}, v_{10}) \in \mathcal{D}\) and let \((\overline{u}_{10}, \overline{v}_{10}) \in \Pi^{s}\) be its image under the Poincar\'e map \(T\). We denote the quantity \(\frac{\overline{v}_{10}}{\overline{u}_{10}}\) (when \(\overline{u}_{10}\neq 0\)) by \(\overline{w}\left(u_{10},v_{10}\right)\) or \(\overline{w}\).
\end{mynotation}

\begin{mynotation}\label{Notiou7v86c54x3}
We denote the straight line \(\lbrace v_{10} = \frac{d}{b} u_{10}\rbrace\) in \(\Pi^{s}\) by \(\ell^{*}\).
\end{mynotation}

\subsection{Dynamics near the homoclinic orbit \(\Gamma\): case \(\lambda_{2} < 2\lambda_{1}\)}\label{Dynamicsnearasinglehomoclinictrivialcase}

\begin{figure}
\centering
\begin{subfigure}{0.4\textwidth}
\centering
\includegraphics[scale=0.19]{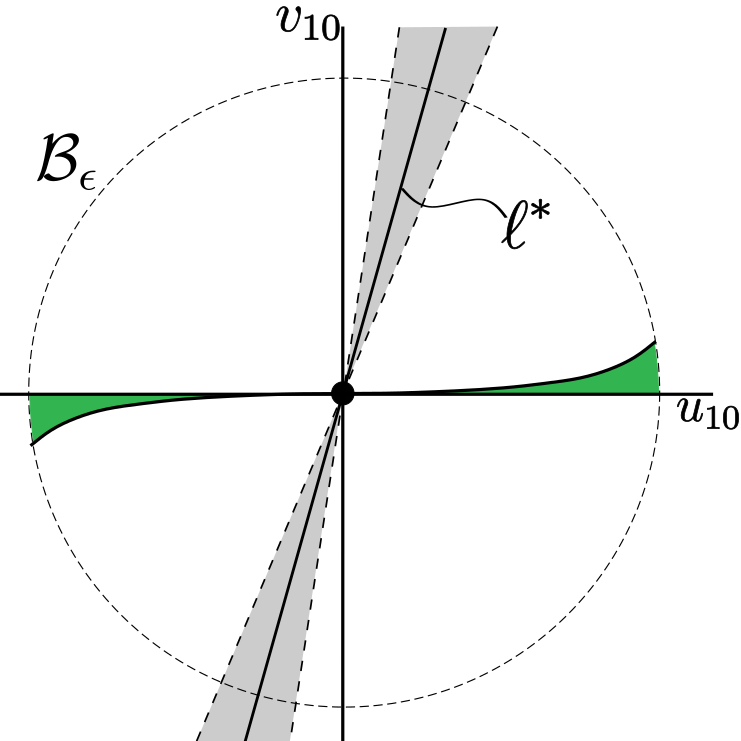}
\end{subfigure}
\begin{subfigure}{0.4\textwidth}
\centering
\includegraphics[scale=0.19]{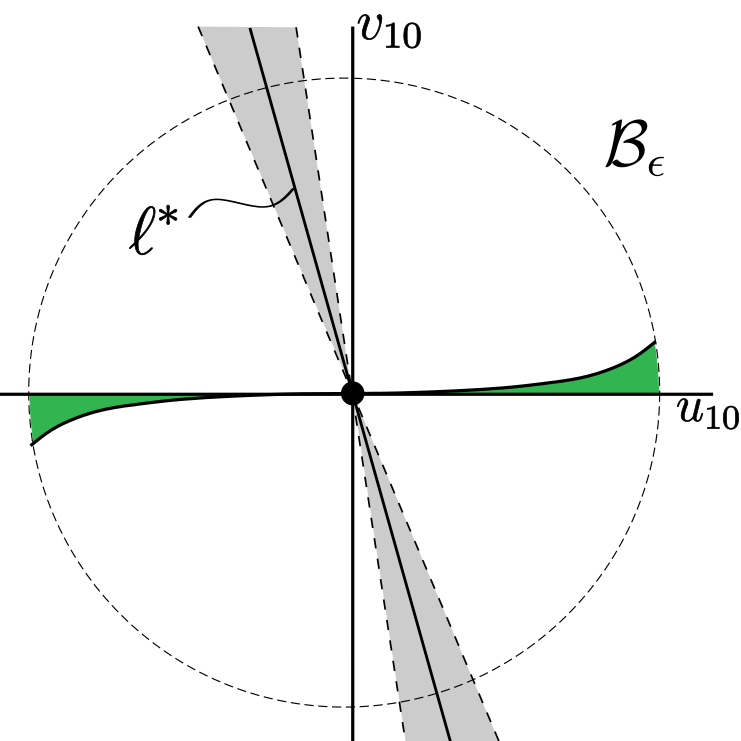}
\end{subfigure}
\caption{\small Case \(\lambda_{1} < \lambda_{2} < 2\lambda_{1}\): the domain \(\mathcal{D}\) of the Poincar\'e map \(T\) is shown in green. The images of the points in \(\mathcal{D}\) under the Poincar\'e map \(T\) accumulate near the straight line \(\ell^{*}\) (the line whose slope is \(\frac{d}{b}\)) in the gray region. As it is shown, the green and the gray regions have no intersection which means \(\mathcal{D} \cap T\left(\mathcal{D}\right) = \emptyset\). This implies that the backward and forward orbits of any point of the domain \(\mathcal{D}\) leaves \(\mathcal{D}\). The left and right figures correspond to the cases \(bd > 0\) and \(bd < 0\), respectively.}
\label{Fig1uim3jj4hk4ji4j4}
\end{figure}

Here, we show that when \(\lambda_{2} < 2\lambda_{1}\), any point in the domain \(\mathcal{D}\) of the Poincar\'e map \(T\) leaves \(\mathcal{D}\) by both forward and backward iterations of the Poincar\'e map. The proof of the case \(\lambda_{1} = \lambda_{2}\) directly follows from Proposition \ref{Cor7uijnqnjk3o49inh44} in which we have shown that the domain \(\mathcal{D}\) of the Poincar\'e map is empty. For the case of \(\lambda_{1} < \lambda_{2} < 2\lambda_{1}\), we prove that the image of the domain \(\mathcal{D}\) under the Poincar\'e map \(T\) has no intersection with \(\mathcal{D}\) (see Figure \ref{Fig1uim3jj4hk4ji4j4}). We formalize this discussion in the following lemma:

\begin{mylem}\label{Lem74njk2po2pn}
When \(\lambda_{2} < 2\lambda_{1}\), we have \(\Lambda^{s}_{\mathcal{D}, T} = \Lambda^{u}_{\mathcal{D}, T} = \emptyset\).
\end{mylem}

\begin{proof}
When \(\lambda_{1} = \lambda_{2}\), the statement follows from Proposition \ref{Cor7uijnqnjk3o49inh44}.

Suppose \(\lambda_{1} < \lambda_{2} < 2\lambda_{1}\). By Proposition \ref{Cor8i12sehu8901fdg}, the domain \(\mathcal{D}\) of the Poincar\'e map is
\begin{equation*}
\{ \left(u_{10}, v_{10}\right) \in \Pi^{s}: u_{10}v_{10} > 0,\quad\lVert \left(u_{10}, v_{10}\right)\rVert < \epsilon, \quad \lvert v_{10}\rvert < K_{\epsilon_{u}} \lvert u_{10}\rvert^{\frac{\gamma}{1-\gamma}} \left[1 + O\left(\delta\right)\right]\},
\end{equation*}
where \(\frac{1}{2} < \gamma =\frac{\lambda_{1}}{\lambda_{2}} < 1\) and \(K_{\epsilon_{u}}\) is some constant (see (\ref{equios093mljoiudmkl28plso})). Since \(u_{1\tau} = v_{1\tau}O\left(\epsilon^{2}\right)\), Poincar\'e map (\ref{eq78uikolp23456fga}) can be written as
\begin{equation*}
\left(\overline{u}_{10}, \overline{v}_{10}\right) = \Big(\left[b + O\left(\epsilon^{2}\right)\right] v_{1\tau}, \left[d + O\left(\epsilon^{2}\right)\right] v_{1\tau}\Big),
\end{equation*}
which implies \(\overline{w} = \frac{d}{b} + O\left(\epsilon^{2}\right)\). This means that the images of the points in the domain \(\mathcal{D}\) under the Poincar\'e map \(T\) accumulate near \(\ell^{*}\). However, for a fixed \(\delta\) and a sufficiently small \(\epsilon\), this line has no intersection with the domain \(\mathcal{D}\) (see Figure \ref{Fig1uim3jj4hk4ji4j4}). This implies \(\Lambda^{s}_{\mathcal{D}, T} = \Lambda^{u}_{\mathcal{D}, T} = \emptyset\), as desired.
\end{proof}

We can now prove Theorems \ref{Th782ju2iu22} (case \(\lambda_{2} < 2\lambda_{1}\)) and \ref{thmkuyvuyuy090988}. The proof of Theorem \ref{Th782ju2iu22} for the case \(2\lambda_{1} < \lambda_{2}\) is provided in the next section.
\begin{proof}[Proof of Theorem \ref{Th782ju2iu22}: case \(\lambda_{2} < 2\lambda_{1}\)]
The proof is an immediate consequence of Remark \ref{Remiuboyboweouberuvg} and Lemma \ref{Lem74njk2po2pn}.
\end{proof}
\begin{proof}[Proof of Theorem \ref{thmkuyvuyuy090988}]
Any orbit in \(W^{s}_{loc}\left(\Gamma\right)\) other than \(\Gamma\) must intersect \(\Pi^{s}\) at \(\Lambda_{\mathcal{D}, T}^{s}\). However, by Lemma \ref{Lem74njk2po2pn}, we have \(\Lambda_{\mathcal{D}, T}^{s} = \emptyset\). This implies \(W^{s}_{loc}\left(\Gamma\right) = \Gamma\). The proof of \(W^{u}_{loc}\left(\Gamma\right) = \Gamma\) is the same.
\end{proof}

\subsection{Dynamics near the homoclinic orbit \(\Gamma\): case \(2\lambda_{1} < \lambda_{2}\)}\label{kjhliub812ftf878g3x7c}

In this section, we study the dynamics near the homoclinic orbit \(\Gamma\) for the case \(2\lambda_{1} < \lambda_{2}\), and prove Theorems \ref{Th782ju2iu22} (case \(2\lambda_{1} < \lambda_{2}\)) and \ref{Invariantmanifoldthm}.

\begin{figure}
\centering
\includegraphics[scale=.22]{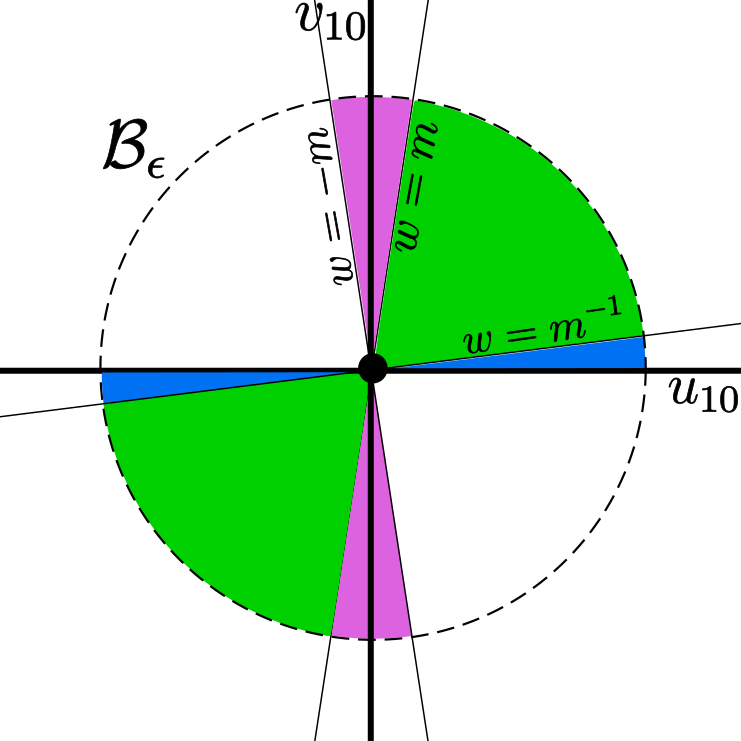}
\caption{\small When \(2\lambda_{1} < \lambda_{2}\), we write the domain \(\mathcal{D}\) of the Poincar\'e map \(T\) as the disjoint union of three subsets \(\mathcal{D}_{1}\), \(\mathcal{D}_{2}\) and \(\mathcal{D}_{3}\), i.e. \(\mathcal{D} = \mathcal{D}_{1} \cup\mathcal{D}_{2}\cup \mathcal{D}_{3}\). The subset \(\mathcal{D}_{1}\) is shown in blue and \(\mathcal{D}_{2}\) is shown in green. The set \(\mathcal{D}_{3}\) is a subset of the purple region.}
\label{FigDomainsD1D2D3aaaaaaa7777uytvyt}
\end{figure}

Recall from Section \ref{Domainjyagvxjycefjthc2tc65} that when \(2\lambda_{1} < \lambda_{2}\), we divide the domain \(\mathcal{D}\) of the Poincar\'e map \(T\) into three subsets \(\mathcal{D}_{1}\), \(\mathcal{D}_{2}\) and \(\mathcal{D}_{3}\), i.e. \(\mathcal{D} = \mathcal{D}_{1} \cup\mathcal{D}_{2}\cup \mathcal{D}_{3}\) (see Figure \ref{FigDomainsD1D2D3aaaaaaa7777uytvyt}). In order to understand the dynamics near the homoclinic loop \(\Gamma\), we need to investigate the set of the points on the domain \(\mathcal{D}\) whose forward or backward orbits (under the iterations of the Poincar\'e map \(T\)) lie in \(\mathcal{D}\), i.e. the sets \(\Lambda^{s}_{\mathcal{D}, T}\) and \(\Lambda^{u}_{\mathcal{D}, T}\) (see Notation \ref{Notation6735oilboiyvkiutc}). To this end, we take the following three steps:
\begin{itemize}
\item {\it Step 1}: Investigating the set of the points in \(\mathcal{D}_{2} \cup \mathcal{D}_{3}\) whose forward or backward orbits lie entirely in \(\mathcal{D}_{2} \cup \mathcal{D}_{3}\), i.e. the sets \(\Lambda^{s}_{\mathcal{D}_{2} \cup \mathcal{D}_{3}, T}\) and \(\Lambda^{u}_{\mathcal{D}_{2} \cup \mathcal{D}_{3}, T}\).
\item {\it Step 2}: Investigating the set of the points in \(\mathcal{D}_{1}\) whose forward or backward orbits lie entirely in \(\mathcal{D}_{1}\), i.e. the sets \(\Lambda^{s}_{\mathcal{D}_{1}, T}\) and \(\Lambda^{u}_{\mathcal{D}_{1}, T}\).
\end{itemize}
Obviously, \(\Lambda^{s}_{\mathcal{D}_{1}, T}\) and \(\Lambda^{s}_{\mathcal{D}_{2} \cup \mathcal{D}_{3}, T}\) are subsets of \(\Lambda^{s}_{\mathcal{D}, T}\). In addition, \(\Lambda^{u}_{\mathcal{D}_{1}, T}\) and \(\Lambda^{u}_{\mathcal{D}_{2} \cup \mathcal{D}_{3}, T}\) are subsets of \(\Lambda^{u}_{\mathcal{D}, T}\). In the third step, we show that the reverse directions also hold: \(\Lambda^{s}_{\mathcal{D}, T} \subset \Lambda^{s}_{\mathcal{D}_{1}, T} \cup \Lambda^{s}_{\mathcal{D}_{2} \cup \mathcal{D}_{3}, T}\) and \(\Lambda^{u}_{\mathcal{D}, T}  \subset \Lambda^{u}_{\mathcal{D}_{1}, T} \cup  \Lambda^{u}_{\mathcal{D}_{2} \cup \mathcal{D}_{3}, T}\). Equivalently,
\begin{itemize}
\item {\it Step 3}: We show \(\Lambda^{s}_{\mathcal{D}, T} = \Lambda^{s}_{\mathcal{D}_{1}, T} \cup \Lambda^{s}_{\mathcal{D}_{2} \cup \mathcal{D}_{3}, T}\) and \(\Lambda^{u}_{\mathcal{D}, T} = \Lambda^{u}_{\mathcal{D}_{1}, T} \cup \Lambda^{u}_{\mathcal{D}_{2} \cup \mathcal{D}_{3}, T}\).
\end{itemize}

Notice that the statement of Step 3 is not trivial. In fact, at the first stage, one can consider the possibility of the existence of a point \(x\in\mathcal{D}\) such that its forward orbit lies entirely in \(\mathcal{D}\), i.e. \(x\in \Lambda^{s}_{\mathcal{D}, T}\), but it does not lie entirely in only one of the sets \(\mathcal{D}_{1}\) or \(\mathcal{D}_{2} \cup \mathcal{D}_{3}\), i.e. \(x\notin\Lambda^{s}_{\mathcal{D}_{1}, T}\) and \(x\notin \Lambda^{s}_{\mathcal{D}_{2} \cup \mathcal{D}_{3}, T}\). In other words, the forward orbit of \(x\) stays in \(\mathcal{D}\) but switches between \(\mathcal{D}_{1}\) and \(\mathcal{D}_{2} \cup \mathcal{D}_{3}\). In Step 3, we indeed show that this scenario does not happen.

We take Step 1 in the following lemma. This lemma helps us to understand the dynamics of the Poincar\'e map \(T\) on the set \(\mathcal{D}_{2}\cup \mathcal{D}_{3}\). We explore in this lemma how \(T\) behaves on this set, with which rate the orbits of this set grow, and how \(\Lambda^{s}_{\mathcal{D}_{2} \cup \mathcal{D}_{3}, T}\) and \(\Lambda^{u}_{\mathcal{D}_{2} \cup \mathcal{D}_{3}, T}\) look like. From a technical point of view, part (\ref{Item892ib3bu8yi4y006}) of this lemma which shows the existence of the unstable manifold of the Poincar\'e map \(T\) is the main result of this section. The techniques which are used in the proof of this part are also used in Section \ref{Dynamicsnearfigure8homoclinic} for the proof of the existence of the unstable manifold of the homoclinic figure-eight. We prove Lemma \ref{Lem899i3n9n9in3o2ubqc} in Section \ref{ProofofLemma11111liuo8784v}.

\begin{mylem}\label{Lem899i3n9n9in3o2ubqc}
Let \(w\) and \(\ell^{*}\) be as in Notations \ref{Not89kookkjweer} and \ref{Notiou7v86c54x3}, respectively. Assume \(2\lambda_{1} < \lambda_{2}\) and consider \(\left(u_{10}, v_{10}\right) \in \mathcal{D}_{2} \cup \mathcal{D}_{3}\). Then
\begin{enumerate}[(i)]
\item\label{Item892ib3bu8yi4y001} \(\overline{w} = w\left(T\left(u_{10}, v_{10}\right)\right) = \frac{d}{b} + o(1)\), where \(o\left(1\right)\) stands for a function of \((u_{10}, v_{10})\) that converges to zero as \((u_{10}, v_{10})\rightarrow (0,0)\).
\item\label{Item892ib3bu8yi4y002} There exists a constant \(C>0\) such that \(\left\Vert \left(u_{10}, v_{10}\right)\right\Vert^{1-2\gamma} < C \left\Vert T\left(u_{10}, v_{10}\right)\right\Vert\) holds for arbitrary \(\left(u_{10}, v_{10}\right)\), where \(\gamma = \lambda_{1}{\lambda_{2}}^{-1} < 0.5\).
\item\label{Item892ib3bu8yi4y003} if \(bd > 0\), then \(T\left(u_{10}, v_{10}\right)\) lies in \(\mathcal{D}_{2}\) unless it leaves \(\mathcal{B}_{\epsilon}\).
\item\label{Item892ib3bu8yi4y004} \(\Lambda^{s}_{\mathcal{D}_{2} \cup \mathcal{D}_{3}, T} = \emptyset\).
\item\label{Item0oubuv18v81716c71c} \(\Lambda^{u}_{\mathcal{D}_{2} \cup \mathcal{D}_{3}, T} = \Lambda^{u}_{\mathcal{D}_{2}, T}\)
\item\label{Item892ib3bu8yi4y005} when \(bd < 0\) we have \(\Lambda^{u}_{\mathcal{D}_{2} \cup \mathcal{D}_{3}, T} = \emptyset\).
\item\label{Item892ib3bu8yi4y006} when \(bd > 0\), the set \(\lbrace M^{s}\rbrace \cup\Lambda^{u}_{\mathcal{D}_{2} \cup \mathcal{D}_{3}, T}\) is a one-dimensional \(\mathcal{C}^{1}\)-manifold which is tangent to \(\ell^{*}\) at \(M^{s}\).
\end{enumerate}
\end{mylem}

\begin{figure}
\centering
\begin{subfigure}{0.4\textwidth}
\centering
\includegraphics[scale=0.23]{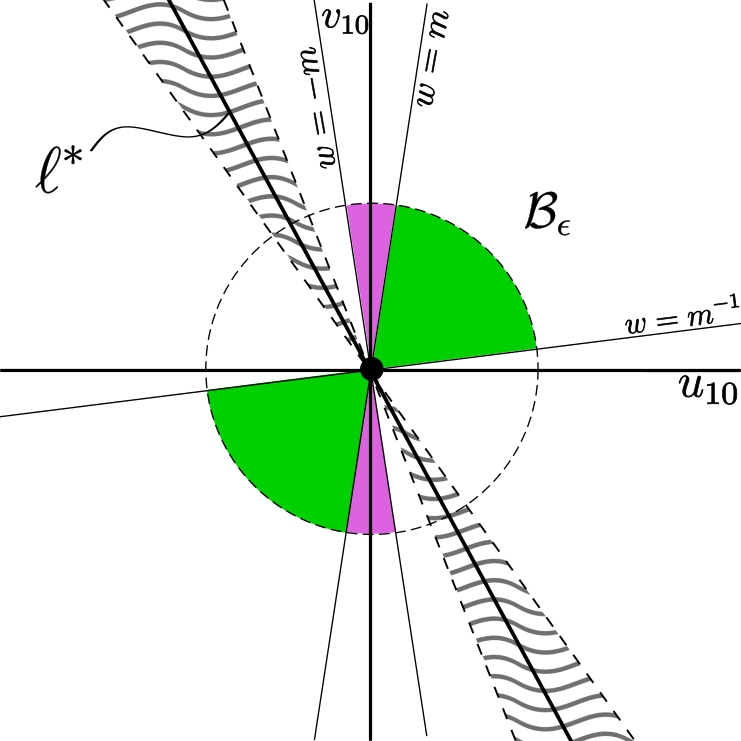}
\end{subfigure}
\begin{subfigure}{0.4\textwidth}
\centering
\includegraphics[scale=0.23]{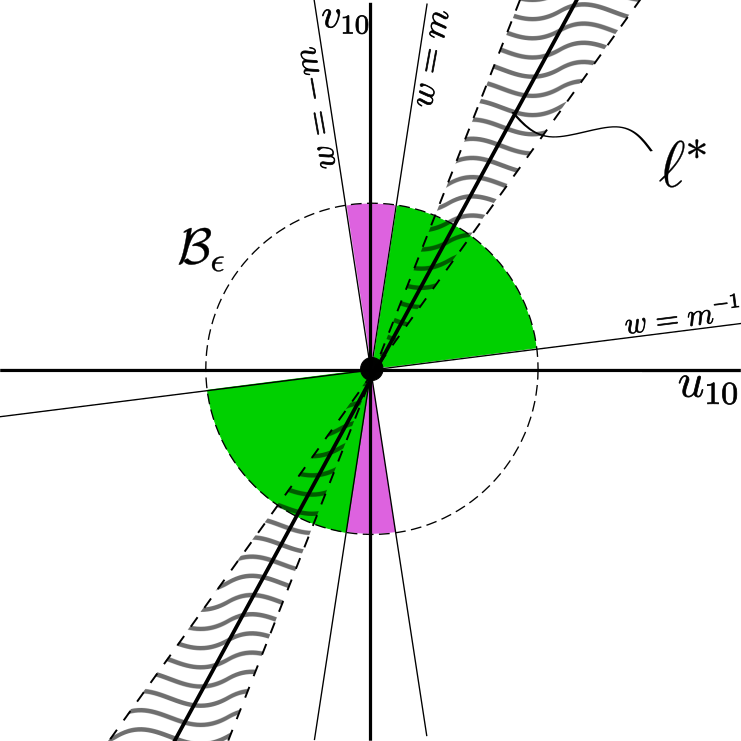}
\end{subfigure}
\caption{\small The straight line whose slope is \(\frac{d}{b}\) is denoted by \(\ell^{*}\). The left figure corresponds to the case \(bd<0\) and the right one corresponds to the case \(bd>0\). The set \(\mathcal{D}_{2}\) is shown in green. The set \(\mathcal{D}_{3}\) is a subset of the purple region. the image of \(\mathcal{D}_{2} \cup \mathcal{D}_{3}\) under the Poincar\'e map, i.e. \(T\left(\mathcal{D}_{2} \cup \mathcal{D}_{3}\right)\), is a subset of the wavy region. Informally speaking, the Poincar\'e map \(T\) {\it preserves} and {\it expands} the region \(\mathcal{D}_{2} \cup \mathcal{D}_{3}\). We show in Lemma \ref{Lem899i3n9n9in3o2ubqc} that when \(bd >0\), there exists an unstable invariant manifold for the Poincar\'e map \(T\) in the wavy region, tangent to \(\ell^*\) at \(M^{s}\).}
\label{Figure10989be378bbbbbw}
\end{figure}

It follows from this lemma that the image of \(\mathcal{D}_{2} \cup \mathcal{D}_{3}\) under the Poincar\'e map \(T\) lies near \(\ell^{*}\), and the Poincar\'e map increases the norm of any point of this set. Informally speaking, for the particular case of \(bd >0\), this means that the Poincar\'e map \(T\) {\it preserves} and {\it expands} the region \(\mathcal{D}_{2} \cup \mathcal{D}_{3}\). A geometrical picture of this behavior is illustrated in Figure \ref{Figure10989be378bbbbbw}.

We now take the second step in the next lemma. In this lemma, we study the dynamics of \(T^{-1}\) on the set \(\mathcal{D}_{1}\). Most of the statements of the following lemma are analogous to the statements of the preceding lemma. This is not a coincidence. In fact, we see later in the proof of Lemma \ref{Lem6647ububquybdu1tyqbuq} that the dynamics of \(T^{-1}\) on \(\mathcal{D}_{1}\) can be obtained from the dynamics of \(T\) on \(\mathcal{D}_{2}\cup\mathcal{D}_{3}\) by a permutation and reversion of time. The proof of this lemma is postponed to Section \ref{ProofofLemma22222liuo8784v}.

\begin{mylem}\label{Lem6647ububquybdu1tyqbuq}
Let \(2\lambda_{1} < \lambda_{2}\) and \(\left(u_{10}, v_{10}\right) \in \mathcal{D}_{1}\). 
\begin{enumerate}[(i)]
\item\label{Item345gysiibwiuwy1} if \(cd > 0\), then \(T\left(\mathcal{D}\right) \cap \mathcal{D}_{1} = \emptyset\).
\item\label{Item345gysiibwiuwy2} if \(cd < 0\), then \(w\left(T^{-1}\left(u_{10}, v_{10}\right)\right) = o\left(1\right)\), where \(o\left(1\right)\) stands for a function of \((u_{10}, v_{10})\) that converges to zero as \((u_{10}, v_{10})\rightarrow (0,0)\). In other words, \(T^{-1}\left(\mathcal{D}_{1}\right)\) accumulates near the horizontal axis.
\item\label{Item345gysiibwiuwy3} if \(cd < 0\), then \(\left\Vert \left(u_{10}, v_{10}\right)\right\Vert^{1-2\gamma} < C \left\Vert T^{-1}\left(u_{10}, v_{10}\right)\right\Vert\) for some constant \(C >0\).
\item\label{Item345gysiibwiuwy4} if \(cd < 0\), then \(T^{-1}\left(u_{10}, v_{10}\right)\) remains in \(\mathcal{D}_{1}\) unless it leaves \(\mathcal{B}_{\epsilon}\).
\item \(\Lambda^{s}_{\mathcal{D}_{1}, T^{-1}} = \emptyset\). Equivalently, \(\Lambda^{u}_{\mathcal{D}_{1}, T} = \emptyset\).
\item\label{Item540o0kju2h2j2n2d} if \(cd < 0\), then the set \(\lbrace M^{s}\rbrace \cup\Lambda^{u}_{\mathcal{D}_{1}, T^{-1}}\) (equivalently, the set \(\lbrace M^{s}\rbrace \cup\Lambda^{s}_{\mathcal{D}_{1}, T}\)) is a one-dimensional \(\mathcal{C}^{1}\)-manifold which is tangent to the horizontal axis at \(M^{s}\).
\end{enumerate}
\end{mylem}

In the preceding two lemmas, we have shown that the sets \(\Lambda^{s}_{\mathcal{D}_{2} \cup \mathcal{D}_{3}, T}\) and \(\Lambda^{u}_{\mathcal{D}_{1}, T}\) are always empty. It was also shown that \(\Lambda^{u}_{\mathcal{D}_{2} \cup \mathcal{D}_{3}, T} = \Lambda^{u}_{\mathcal{D}_{2}, T}\). This allows us to reformulate Step 3 as in the following lemma:

\begin{mylem}\label{Lem12wo0ij4u4hhsza0}
\begin{inparaenum}[(i)]
\item \label{Item46710kmnryt} \(\Lambda^{s}_{\mathcal{D}, T} = \Lambda^{s}_{\mathcal{D}_{1}, T}\). \hspace{1cm}
\item \label{Item546uhuuwhu4hugtei}\(\Lambda^{u}_{\mathcal{D}, T} = \Lambda^{u}_{\mathcal{D}_{2}, T}\).
\end{inparaenum}
\end{mylem}

\begin{proof}
Let \(\mathtt{x} \in \mathcal{D}_{2} \cup \mathcal{D}_{3}\). It follows from parts (\ref{Item892ib3bu8yi4y001}), (\ref{Item892ib3bu8yi4y002}) and (\ref{Item892ib3bu8yi4y003}) of Lemma \ref{Lem899i3n9n9in3o2ubqc} that if \(bd < 0\), then \(T\left(\mathtt{x}\right)\notin \mathcal{D}\), and if \(bd > 0\), then for some \(k\), \(T^{k}\left(\mathtt{x}\right)\notin\mathcal{B}_{\epsilon}\). Thus, any point in \(\Lambda^{s}_{\mathcal{D}, T}\) must belong to \(\mathcal{D}_{1}\). This proves part (\ref{Item46710kmnryt}).

To prove part (\ref{Item546uhuuwhu4hugtei}), notice that if \(\Lambda^{u}_{\mathcal{D}, T} = \emptyset\), then \(\Lambda^{u}_{\mathcal{D}_{2}, T} = \emptyset\) and therefore \(\Lambda^{u}_{\mathcal{D}, T} = \Lambda^{u}_{\mathcal{D}_{2}, T}\). So we assume that \(\Lambda^{u}_{\mathcal{D}, T}\) is non-empty. Let \(\mathtt{x} \in \Lambda^{u}_{\mathcal{D}, T}\). We need to show \(\mathtt{x}\in \mathcal{D}_{2}\). To do this, we first prove \(\mathtt{x}\notin \mathcal{D}_{1}\). Assume the contrary, i.e. \(\mathtt{x}\in \mathcal{D}_{1}\). It follows from parts (\ref{Item892ib3bu8yi4y001}) and (\ref{Item892ib3bu8yi4y003}) of Lemma \ref{Lem899i3n9n9in3o2ubqc} that if \(T^{-1}\left(\mathtt{x}\right)\in \mathcal{D}_{2}\cup \mathcal{D}_{3}\), then \(\mathtt{x} = T\left(T^{-1}\left(\mathtt{x}\right)\right)\) either belongs to \(\mathcal{D}_{2}\) or lies outside the domain \(\mathcal{D}\) which contradicts the assumption \(\mathtt{x}\in\mathcal{D}_{1}\). Therefore, \(T^{-1}\left(\mathtt{x}\right)\notin \mathcal{D}_{2}\cup \mathcal{D}_{3}\), and so \(T^{-1}\left(\mathtt{x}\right)\in \mathcal{D}_{1}\). By virtue of part (\ref{Item345gysiibwiuwy1}) of Lemma \ref{Lem6647ububquybdu1tyqbuq}, this relation implies \(cd<0\). On the other hand, when \(cd<0\), it follows from parts (\ref{Item345gysiibwiuwy3}) and (\ref{Item345gysiibwiuwy4}) of Lemma \ref{Lem6647ububquybdu1tyqbuq} that there exists a \(k>0\) such that \(T^{-k}\left(\mathtt{x}\right) \notin \mathcal{B}_{\epsilon}\) and hence \(T^{-k}\left(\mathtt{x}\right) \notin \mathcal{D}\). This contradicts the preliminary assumption \(\mathtt{x} \in \Lambda^{u}_{\mathcal{D}, T}\). Therefore, if \(\mathtt{x} \in \Lambda^{u}_{\mathcal{D}, T}\), then \(\mathtt{x}\notin\mathcal{D}_{1}\), or equivalently, \(T^{-n}\left(\mathtt{x}\right)\notin\mathcal{D}_{1}\) for all \(n\geq 0\).

To finish the proof, it is sufficient to show that \(\mathtt{x}\notin\mathcal{D}_{3}\). Assume the contrary, i.e. \(\mathtt{x}\in\mathcal{D}_{3}\). Since \(\mathtt{x} \in \Lambda^{u}_{\mathcal{D}, T}\) implies \(T^{-n}\left(\mathtt{x}\right)\notin\mathcal{D}_{1}\) for all \(n\geq 0\), we have \(T^{-1}\left(\mathtt{x}\right) \notin\mathcal{D}_{1}\). On the other hand, parts (\ref{Item892ib3bu8yi4y001}) and (\ref{Item892ib3bu8yi4y003}) of Lemma \ref{Lem899i3n9n9in3o2ubqc} imply that if \(T^{-1}\left(\mathtt{x}\right)\in \mathcal{D}_{2}\cup \mathcal{D}_{3}\), then \(\mathtt{x} = T\left(T^{-1}\left(\mathtt{x}\right)\right)\) either belongs to \(\mathcal{D}_{2}\) or lies outside the domain \(\mathcal{D}\) which contradicts the assumption \(\mathtt{x}\in\mathcal{D}_{3}\). Therefore, \(\mathtt{x}\notin\mathcal{D}_{3}\), as desired.
\end{proof}

Recall that the local stable (unstable) set of the homoclinic loop \(\Gamma\), denoted by \(W^{s}_{\text{loc}}(\Gamma)\) (\(W^{u}_{\text{loc}}(\Gamma)\)), is the union of \(\Gamma\) itself and the set of the points in a sufficiently small neighborhood \(\mathcal{U}\) of \(\Gamma\) whose forward (backward) orbits lie in \(\mathcal{U}\) and their \(\omega\)-limit sets (\(\alpha\)-limit sets) coincide with \(\Gamma \cup \lbrace O\rbrace\). By this definition, the intersection of \(W^{s}_{\text{loc}}(\Gamma)\) and \(\Pi^{s}\) must belong to \(\{M^{s}\}\cup\Lambda^{s}_{\mathcal{D}}\), and the intersection of \(W^{u}_{\text{loc}}(\Gamma)\) and \(\Pi^{s}\) must belong to \(\{M^{s}\}\cup\Lambda^{u}_{\mathcal{D}}\). On the other hand, we have shown in the above lemmas that when \(\Lambda^{s}_{\mathcal{D}}\) (\(\Lambda^{u}_{\mathcal{D}}\)) is non-empty, any point on this set converges to \(M^{s}\) by the forward (backward) iterations of the Poincar\'e map \(T\). This leads to the following:
\begin{myprop}\label{Cor12byeijheuije}
Let \(\phi_{t}\) be the flow of system (\ref{eq23000}). Then
\begin{equation*}
W^{s}_{loc}\left(\Gamma\right) = \Gamma \cup \phi_{t}\left(\Lambda^{s}_{\mathcal{D}, T}\right)\,\text{ for } t\geq 0,\qquad \text{and}\qquad 
W^{u}_{loc}\left(\Gamma\right) = \Gamma \cup \phi_{t}\left(\Lambda^{u}_{\mathcal{D}, T}\right)\,\text{ for } t\leq 0.
\end{equation*}
\end{myprop}

\begin{figure}
\centering
\begin{subfigure}{1.0\textwidth}
\centering
\includegraphics[scale=0.15]{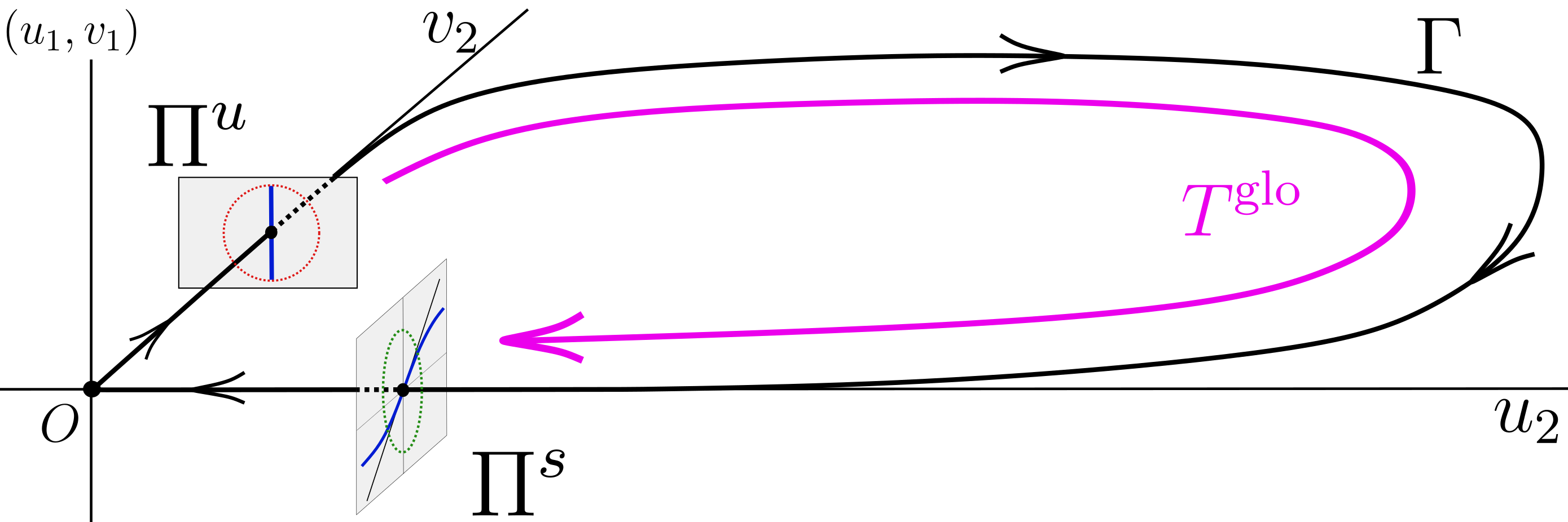}
\caption{}
\end{subfigure}
\newline
\begin{subfigure}{1.0\textwidth}
\centering
\includegraphics[scale=0.12]{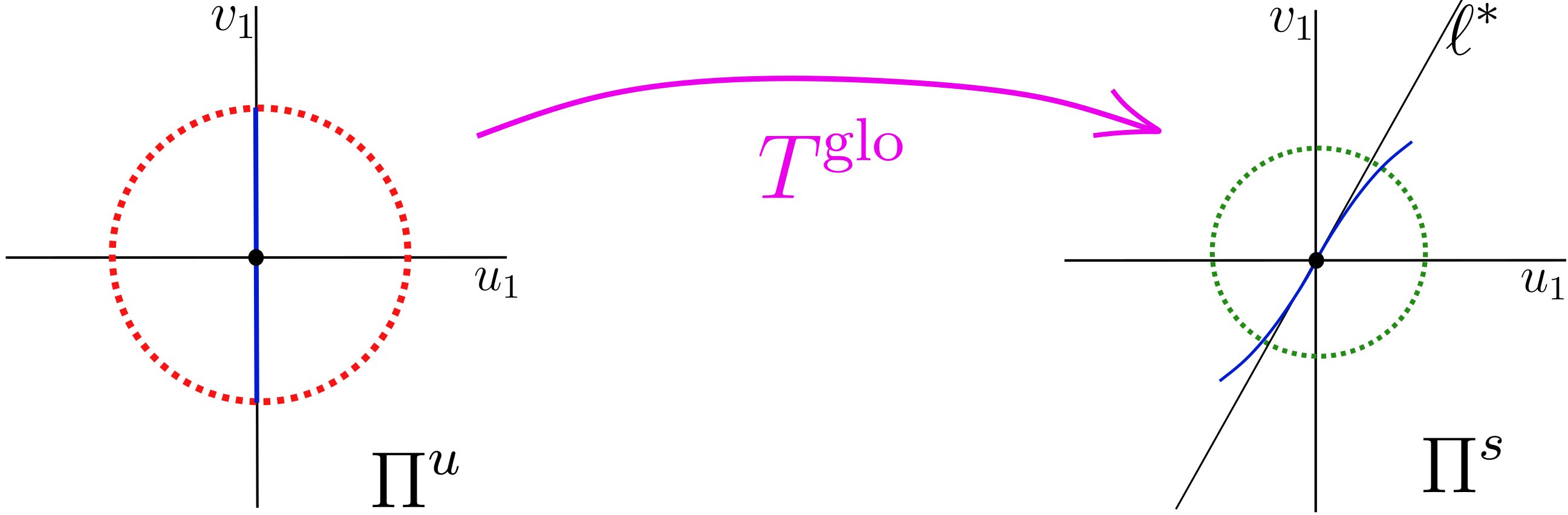}
\caption{}
\end{subfigure}
\caption{\small The local unstable invariant manifold of the equilibrium \(O\) intersects \(\Pi^{u}\) at the \(v_{1}\)-axis. Thus, the blue curve (\(v_{1}\)-axis restricted to a small neighborhood of \(M^{u}\) in \(\Pi^{u}\)) lies at the intersection of the local unstable invariant manifold of \(O\) and the cross-section \(\Pi^{u}\). This curve is mapped to the blue curve on \(\Pi^{s}\) by \(T^{\text{glo}}\) which means that the blue curve on \(\Pi^{s}\) lies in \(W^{u}_{\text{glo}}\left(O\right) \cap \Pi^{s}\). Since \(v_{1}\)-axis on \(\Pi^{u}\) is mapped to \(\ell^{*}\) on \(\Pi^{s}\) by \(d T^{\text{glo}}\), the straight line \(\ell^{*}\) is tangent to the blue curve on \(\Pi^{s}\) at \(M^{s}\).}
\label{Figurejytavxu7c99nc6aa7vx7ca}
\end{figure}

In system (\ref{eq23000}), the local unstable invariant manifold of the equilibrium \(O\) is straightened, i.e. \(W^{u}_{loc}\left(O\right) = \lbrace u = 0\rbrace\). Thus, the intersection of this manifold and the cross-section \(\Pi^{u} = \{v_{2} = \delta\} \cap\{H=0\}\) is the straight line \(\lbrace u_{1} = 0\rbrace\), i.e. \(v_{1}\)-axis. Consider the restriction of this line to a small neighborhood of \(M^{u}\) (in Figure \ref{Figurejytavxu7c99nc6aa7vx7ca}, it is shown by blue color on \(\Pi^{u}\)). The global map \(T^{\text{glo}}\) maps this restricted piece to a curve, denote it by \(\gamma^{u}\), on \(\Pi^{s}\) (shown by blue color on  \(\Pi^{s}\) in Figure \ref{Figurejytavxu7c99nc6aa7vx7ca}). This curve is in fact at the intersection of the global unstable invariant manifold of \(O\) and the cross-section \(\Pi^{s}\). Since \(T^{\text{glo}}\) is a diffeomorphism and the vector \(\left(\begin{smallmatrix}
0 \\ 1
\end{smallmatrix}\right)\) is tangent to \(v_{1}\)-axis at \(M^{u}\), the vector \(d T^{\text{glo}}\left(\begin{smallmatrix}
0 \\ 1
\end{smallmatrix}\right) = \left(\begin{smallmatrix}
b \\ d
\end{smallmatrix}\right)\) is tangent to \(\gamma^{u}\) at \(M^{s}\) in \(\Pi^{s}\), i.e. \(\gamma^{u}\) is tangent to \({\ell}^{*}\) at \(M^{s}\) (recall that \(\ell^{*}\) is the line in \(\Pi^{s}\) whose slope is \(\frac{d}{b}\)). Therefore, it follows from Lemma \ref{Lem12wo0ij4u4hhsza0} and part (\ref{Item892ib3bu8yi4y006}) of Lemma \ref{Lem899i3n9n9in3o2ubqc} that when \(bd > 0\), \(W^{u}_{glo}\left(O\right) \cap \Pi^{s}\) and \(\lbrace M^{s}\rbrace \cup\Lambda^{u}_{\mathcal{D}, T}\) are tangent at \(M^{s}\). On the other hand, it follows from Lemma \ref{Lem12wo0ij4u4hhsza0} and part (\ref{Item540o0kju2h2j2n2d}) of Lemma \ref{Lem6647ububquybdu1tyqbuq} that when \(cd < 0\) the intersection of the local stable manifold of \(O\) and the cross-section \(\Pi^{s}\), i.e. the horizontal axis, is tangent to \(\lbrace M^{s}\rbrace \cup\Lambda^{s}_{\mathcal{D}, T}\) at \(M^{s}\).  Moreover, by Assumption \ref{assumption50}, the homoclinic orbit \(\Gamma\) is at the transverse intersection of the global stable and unstable invariant manifolds of the equilibrium \(O\). Therefore, the intersection of these two manifolds with the cross-section \(\Pi^{s}\), i.e. the horizontal axis and the curve \(\gamma^{u}\), intersect transversely at \(M^{s}\). Since \(\gamma^{u}\) is tangent to \(\ell^{*}\) at \(M^{s}\), we have that the intersection of \(W_{\text{glo}}^{s}\left(O\right)\) and \(W_{\text{glo}}^{u}\left(O\right)\) at \(\Gamma\) is transverse if and only if the horizontal axis on \(\Pi^{s}\) and the straight line \(\ell^{*}\) are distinct. These statements give
\begin{myprop}\label{Prop453ityuuuihbbdhjdndaa}
\begin{enumerate}[(i)]
\item When \(bd > 0\), the 2-dimensional \(\mathcal{C}^1\)-smooth invariant manifold \(W^{u}_{\text{loc}}\left(\Gamma\right)\) is tangent to \(W_{\text{glo}}^{u}\left(O\right)\) at every point of \(\Gamma\).
\item When \(cd < 0\), the 2-dimensional \(\mathcal{C}^1\)-smooth invariant manifold \(W^{s}_{\text{loc}}\left(\Gamma\right)\) is tangent to \(W_{\text{glo}}^{s}\left(O\right)\) at every point of \(\Gamma\).
\item The intersection of \(W_{\text{glo}}^{s}\left(O\right)\) and \(W_{\text{glo}}^{u}\left(O\right)\) at \(\Gamma\) is transverse if and only if \(d \neq 0\).
\end{enumerate}
\end{myprop}

By virtue of the above results, we can prove Theorems \ref{Th782ju2iu22} (case \(2\lambda_{1}< \lambda_{2}\)) and \ref{Invariantmanifoldthm}:
\begin{proof}[Proof of Theorem \ref{Th782ju2iu22}: case \(2\lambda_{1}< \lambda_{2}\)]
The proof is an immediate consequence of Remark \ref{Remiuboyboweouberuvg} and Lemmas \ref{Lem899i3n9n9in3o2ubqc}, \ref{Lem6647ububquybdu1tyqbuq} and \ref{Lem12wo0ij4u4hhsza0}.
\end{proof}
\begin{proof}[Proof of Theorem \ref{Invariantmanifoldthm}]
By Proposition \ref{Cor12byeijheuije} and the preceding Lemmas we have that \(W^{u}_{\text{loc}}\left(\Gamma\right) = \Gamma\) when \(bd < 0\), and \(W^{s}_{\text{loc}}\left(\Gamma\right) = \Gamma\) when \(cd > 0\). The rest of the theorem is already proved (see Proposition \ref{Prop453ityuuuihbbdhjdndaa}).
\end{proof}

The following remark suggests an alternative formulation of Theorem \ref{Invariantmanifoldthm}:
\begin{myrem}\label{Remou8auvuyfvuflugfvugfjyg}
Consider the global stable and unstable invariant manifolds of the equilibrium \(O\) of system (\ref{eq23000}). Let \(\gamma^{u}\) (resp. \(\gamma^{s}\)) be a curve at the intersection of the global unstable (resp. stable) invariant manifold of \(O\) and the cross-section \(\Pi^{s}\) (resp. \(\Pi^{u}\)) which passes through \(M^{s}\) (resp. \(M^{u}\)). Following the discussion above, the slope of the tangent line to the curve \(\gamma^{u}\) at \(M^{s}\) is \(\frac{d}{b}\). Moreover, the slope of the tangent line to the curve \(\gamma^{s}\) at \(M^{u}\) is \(\frac{-c}{d}\). This suggests an alternative way to detect \(\text{sgn}\left(bd\right)\) and \(\text{sgn}\left(cd\right)\) which are required in the statement of Theorem \ref{Invariantmanifoldthm}. Indeed, instead of computing the coefficients \(a\), \(b\), \(c\) and \(d\) in Theorem \ref{Invariantmanifoldthm}, one can look at the slopes of the intersection curves of the global stable and unstable invariant manifolds of the equilibrium \(O\) with the cross-sections \(\Pi^{u}\) and \(\Pi^{s}\) at the points \(M^{u}\) and \(M^{s}\).
\end{myrem}

\subsubsection{Proof of Lemma \ref{Lem899i3n9n9in3o2ubqc}}\label{ProofofLemma11111liuo8784v}

\begin{proof}[Proof of part (\ref{Item892ib3bu8yi4y001})]
By (\ref{eq679oo89lkikmn44hujj}) and Proposition \ref{Cor67ujokjsuijwww}, \((u_{10}, v_{10}) \in \mathcal{D}_{2} \cup \mathcal{D}_{3}\) implies \(u_{1\tau} = o(v_{1\tau})\). Thus, Poincar\'e map (\ref{eq78uikolp23456fga}) takes the form
\begin{equation}\label{eq7890kiowert}
\left(\overline{u}_{10}, \overline{v}_{10}\right) = \Big(\left[b + o\left(1\right)\right] v_{1\tau}, \left[d + o\left(1\right)\right] v_{1\tau}\Big),
\end{equation}
which implies \(\overline{w} = \frac{d}{b} + o\left(1\right)\).
\end{proof}

\begin{proof}[Proof of part (\ref{Item892ib3bu8yi4y002})]
For \(\left(u_{10}, v_{10}\right)\in \mathcal{D}_{2}\), relations (\ref{eq8isolpaloue74jndu}) and (\ref{eq7890kiowert}) imply
\begin{equation*}
\left\Vert T\left(u_{10}, v_{10}\right)\right\Vert = \left\Vert \left(\overline{u}_{10}, \overline{v}_{10}\right)\right\Vert = \left[b^{2} + d^{2} + o\left(1\right)\right]^{\frac{1}{2}}\left\vert v_{1\tau}\right\vert = K \left\vert u_{10}\right\vert^{-\gamma} \left\vert v_{10}\right\vert^{1-\gamma},
\end{equation*}
where \(K = K\left(u_{10}, v_{10}\right) = \gamma^{-\gamma}\delta^{2\gamma}\sqrt{b^{2} + d^{2} + o\left(1\right)}\). For \(C > K^{-1} m^{\gamma} \left(1 + m\right)^{\frac{1}{2}}\), we have
\begin{equation*}
\frac{\left\Vert \left(u_{10}, v_{10}\right)\right\Vert}{\left\Vert T\left(u_{10}, v_{10}\right)\right\Vert} = \frac{\lvert v_{10}\rvert \sqrt{1 + \lvert\frac{u_{10}}{v_{10}}\rvert}}{K \left\vert u_{10}\right\vert^{-\gamma} \left\vert v_{10}\right\vert^{1-\gamma}} \leq K^{-1} m^{\gamma} \left(1 + m\right)^{\frac{1}{2}} \lvert v_{10}\rvert^{2\gamma} < C \left\Vert \left(u_{10}, v_{10}\right)\right\Vert^{2\gamma},
\end{equation*}
as desired.
\end{proof}

\begin{proof}[Proof of part (\ref{Item892ib3bu8yi4y003})]
By part (\ref{Item892ib3bu8yi4y001}) of Lemma \ref{Lem899i3n9n9in3o2ubqc}, \(T\left(u_{10}, v_{10}\right)\) is somewhere close to the line \(\ell^{*}\) and since, for \(bd>0\), the restriction of \(\ell^{*}\setminus\lbrace M^{s}\rbrace\) to \(\mathcal{B}_{\epsilon}\) lies in \(\mathcal{D}_{2}\) we have that if \(T\left(u_{10}, v_{10}\right)\) lies in \(\mathcal{B}_{\epsilon}\), then it must belong to \(\mathcal{D}_{2}\).
\end{proof}

\begin{proof}[Proofs of parts (\ref{Item892ib3bu8yi4y004}), (\ref{Item0oubuv18v81716c71c}) and (\ref{Item892ib3bu8yi4y005})]
All are easy consequences of the previous parts.
\end{proof}

\begin{proof}[Proof of part (\ref{Item892ib3bu8yi4y006})]
The proof is based on a theory of invariant manifolds for cross-maps (see \cite{Dimabook} and \cite{GonchenkoCrossmap2010}). An introduction to this theory is provided in Appendix \ref{Invariantmanifoldscrossmaps}.

Consider \(\mathcal{D}^{\epsilon_{1}}_{2}\) for a sufficiently small \(\epsilon_{1} > 0\) (see (\ref{eqiuvyseyt74vqi8l3viul181})). Choose \(\epsilon_{2} < \epsilon_{1}\) such that \(\mathcal{X} \subset \mathcal{D}_{2}^{\epsilon_{1}}\), where \(\mathcal{X} = \lbrace \left(u_{10}, v_{10}\right) \in \Pi^{s}: m^{-1} \leq  \frac{v_{10}}{u_{10}} \leq m,\,\, u_{10}\neq 0,\,\, \lvert v_{10}\rvert \leq \epsilon_{2}\rbrace\) and \(m\) is as in (\ref{eqyujnj68999jhhhgg34}) (see Figure \ref{FigXinu10-v10coordinates}). Recall \(w\) in Notation \ref{Not89kookkjweer} and define the new variable \(z\) by
\begin{equation}\label{P800980}
z = z\left(u_{10}, v_{10}\right) = \text{sgn}\left(v_{10}\right) \lvert v_{10}\rvert^{\alpha}, \qquad (0 < \alpha \,\,\text{will be specified later)}.
\end{equation}

\begin{figure}
\centering
\includegraphics[scale=.25]{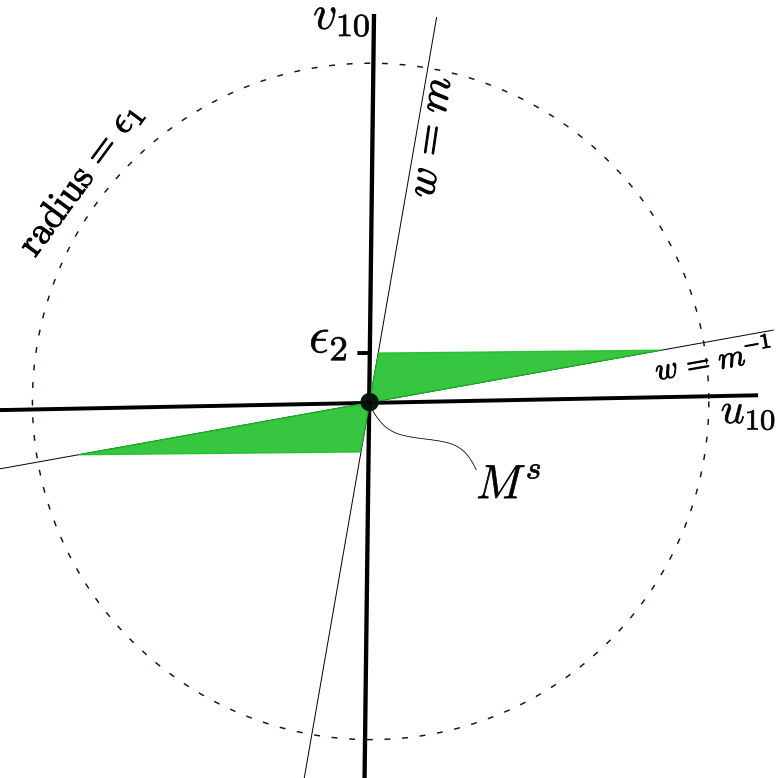}
\caption{\small The set \(\mathcal{X} \subset \mathcal{D}_{2}^{\epsilon_{1}}\) in \(\left(u_{10}, v_{10}\right)\)-plane is shown by green color.}
\label{FigXinu10-v10coordinates}
\end{figure}

Let \(\mathcal{Y}\) be the set \(\mathcal{X}\) equipped with \(\left(w, z\right)\)-coordinates. Thus, \(\mathcal{Y} = \left[m^{-1}, m\right]\times \left(\left[-{\epsilon_{2}}^{\alpha}, {\epsilon_{2}}^{\alpha}\right]\setminus \lbrace 0\rbrace\right)\) (see Figure \ref{FigY}).  Consider the restriction of the Poincar\'e map \(T\) to the set \(\mathcal{X}\), i.e. \(T\vert_{\mathcal{X}}\), and denote the representation of this map in \(\left(w, z\right)\)-coordinates by \(\mathcal{T}\). We write
\begin{equation}\label{eq78i9o0ij3hujb3gh121}
\mathcal{T}: \left(w, z\right)\mapsto \left(\overline{w}, \overline{z}\right) = \left(f\left(w, z\right), g\left(w, z\right)\right),
\end{equation}
for some smooth functions \(f\) and \(g\) defined on \(\mathcal{Y}\). Note that by (\ref{eq8isolpaloue74jndu}) and the relation \(\overline{z} = g\left(w, z\right) = \text{sgn}\left(\overline{v}_{10}\right) \lvert \overline{v}_{10}\rvert^{\alpha}\), we can derive
\begin{equation}\label{eq89oknehujrnhjnwj}
\overline{z} = g\left(w, z\right) = \text{sgn}\left(dz\right) \lvert d\rvert^{\alpha}\left(\frac{\gamma}{\delta^{2}}\right)^{-\gamma\alpha} w^{\gamma\alpha} \lvert z\rvert^{1 - 2\gamma} \left[1 + O\left(\delta\right)\right]
= O\left(\lvert z\rvert^{1 - 2\gamma}\right)
\end{equation}
and
\begin{equation}\label{P809000}
z = O\left(\lvert \overline{z}\rvert^{\frac{1}{{1 - 2\gamma}}}\right).
\end{equation}
We now make a statement which is proved in Appendix \ref{Appendix989897g86g7656u}:
\begin{mylem}\label{State920jj3nu3uhwu}
\(g_{z}\left(w, z\right)\) is non-zero for any \((w, z)\in\mathcal{Y}\).
\end{mylem}

\begin{figure}
\centering
\includegraphics[scale=.20]{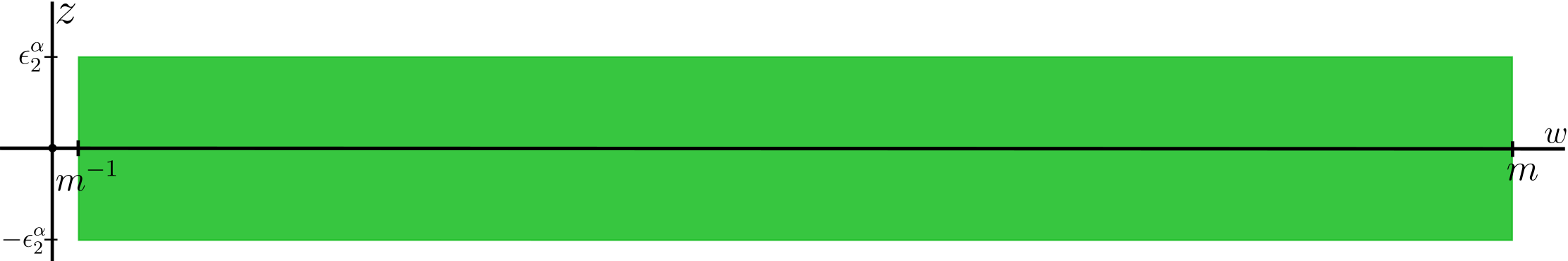}
\caption{\small The set \(\mathcal{Y}\) (the set \(\mathcal{X}\) equipped with \(\left(w, z\right)\)-coordinates) is shown by green color. It contains two connected components (below and above the horizontal axis).}
\label{FigY}
\end{figure}

According to this lemma and the implicit function theorem, the variable \(z\) is a \(\mathcal{C}^{q}\)-smooth (\(q\) is as in Lemma \ref{Nftheorem2}) function of \(\left(w,\overline{z}\right)\) for \(w\in \left[m^{-1}, m\right]\) and \(\overline{z} \in g\left(\mathcal{Y}\right)\). Denote this function by \(G\). Regarding the domain of this function, note that not every \(\left(w, \overline{z}\right)\) necessarily belongs to the domain of \(G\). In other words, for an arbitrary \(\left(w, \overline{z}\right)\), there might not exist \(z\in\left[-{\epsilon_{2}}^{\alpha}, {\epsilon_{2}}^{\alpha}\right]\setminus \lbrace 0\rbrace\) such that \(z = G\left(w, \overline{z}\right)\). However, by (\ref{eq89oknehujrnhjnwj}), this relation holds if \(\overline{z}\) is chosen sufficiently small, i.e. for a sufficiently small \(\theta > 0\) we have
\begin{equation*}
\left[m^{-1}, m\right]\times\left(\left[-\theta, \theta\right]\setminus \lbrace 0\rbrace\right) \subset \text{domain}\left(G\right).
\end{equation*}
Denote this set by \(\mathcal{R}\), i.e. \(\mathcal{R} = \left[m^{-1}, m\right]\times \left(\left[-\theta, \theta\right]\setminus\lbrace 0\rbrace\right)\). Without loss of generality, assume \(\theta < {\epsilon_{2}}^{\alpha}\). Having the function \(G\) in hand means that we can write the Poincar\'e map \(\mathcal{T}\) in cross-form: we define the cross-map \(\mathcal{T}^{\times}: \left(w, \overline{z}\right) \mapsto \left(\overline{w}, z\right)\) by
\begin{equation}\label{P808000}
\left(\overline{w}, z\right) = \left(F\left(w, \overline{z}\right), G\left(w, \overline{z}\right)\right),\qquad \text{where \,\,} F\left(w, \overline{z}\right) = f\left(w, G\left(w, \overline{z}\right)\right),
\end{equation}
and \(\left(w, \overline{z}\right)\in \text{domain}\left(G\right)\). It follows from part (\ref{Item892ib3bu8yi4y001}) of Lemma \ref{Lem899i3n9n9in3o2ubqc} (proved earlier), relation (\ref{P809000}) and the fact that \(z = 0\) if and only if \(\overline{z} = 0\) (follows from (\ref{eq89oknehujrnhjnwj})) that \(\mathcal{T}^{\times}\left(\mathcal{R}\right) \subset \mathcal{R}\). Hereafter, we focus on the restriction of \(\mathcal{T}^{\times}\) on \(\mathcal{R}\). Our approach to prove the existence of the desired invariant manifold for the Poincar\'e map \(\mathcal{T}\) is to apply Theorem \ref{thm3000} (see Appendix \ref{Invariantmanifoldscrossmaps}) on the cross-map \(\mathcal{T}^{\times}\). However, to do this, there are two issues that we need to take care of. The first is that the domain \(\mathcal{R}\) does not satisfy the assumption of Theorem \ref{thm3000} (in that proposition, the domain must be written as a Cartesian product of two convex closed sets but \(\mathcal{R}\) is not of this form since it does not contain the line \(\overline{z} = 0\)). Second, we need to compute the partial derivatives of the cross-map \(\mathcal{T}^{\times}\). The second issue is resolved by the following lemma:
\begin{mylem}\label{Stat4u74i4jm1ml2mo3p}
Let \(\beta = \alpha^{-1} \min\lbrace 4\gamma, 1-2\gamma\rbrace\). We have
\begin{equation*}
\begin{aligned}
F_{w}\left(w, \overline{z}\right) &= O\left(\lvert \overline{z}\rvert^{\frac{\beta}{{1 - 2\gamma}}}\right),\qquad
&& F_{\overline{z}}\left(w, \overline{z}\right) = O\left(\lvert \overline{z}\rvert^{\frac{\beta - 1 + 2\gamma}{{1 - 2\gamma}}}\right),\\
G_{w}\left(w, \overline{z}\right) &= O\left(\lvert \overline{z}\rvert^{\frac{1}{{1 - 2\gamma}}}\right),\qquad
&& G_{\overline{z}}\left(w,\overline{z}\right) = O\left(\lvert \overline{z}\rvert^{\frac{2\gamma}{{1 - 2\gamma}}}\right).
\end{aligned}
\end{equation*}
\end{mylem}
This lemma is proved in Appendix \ref{Appendix989897g86g7656u}. We now extend the domain \(\mathcal{R}\) to \(\widetilde{\mathcal{R}}\), where \(\widetilde{\mathcal{R}} = \left[m^{-1}, m\right]\times \left[-\theta, \theta\right]\). We also extend the map \(\mathcal{T}^{\times}\) to the map \({\widetilde{\mathcal{T}}}^{\times}\) defined on \(\widetilde{\mathcal{R}}\) by
\begin{equation*}
{\widetilde{\mathcal{T}}}^{\times}\left(w, \overline{z}\right) := \left\{\begin{array}{ll}
\mathcal{T}^{\times}\left(w, \overline{z}\right) = \left(F\left(w, \overline{z}\right), G\left(w, \overline{z}\right)\right)\quad & \overline{z} \neq 0, \vspace*{2mm}\\
\left(\frac{d}{b}, 0\right)\quad & \overline{z} = 0,
\end{array}\right.
\end{equation*}
Lemma \ref{Stat4u74i4jm1ml2mo3p} implies that for a fixed sufficiently small \(\alpha\), the map \(\widetilde{\mathcal{T}}^{\times}: \widetilde{\mathcal{R}} \rightarrow \widetilde{\mathcal{R}}\) is a \(\mathcal{C}^{1}\)-smooth extension of \(\mathcal{T}^{\times}\) to \(\widetilde{\mathcal{R}}\).

Now, let us come back to the Poincar\'e map \(\mathcal{T}\) defined on \(\mathcal{Y}\). We extend this map to
\begin{equation*}
\widetilde{\mathcal{T}}\left(w, z\right) := \left\{\begin{array}{ll}
\mathcal{T}\left(w, z\right)\quad & \left(w, z\right)\in \mathcal{Y}, \vspace*{2mm}\\
\left(\frac{d}{b}, 0\right)\quad & z = 0.
\end{array}\right.
\end{equation*}
It is clear that the map \(\widetilde{\mathcal{T}}^{\times}\) is in fact the cross-map of \(\widetilde{\mathcal{T}}\) on \(\widetilde{\mathcal{R}}\). Note that since \(\theta < {\epsilon_{2}}^{\alpha}\), we have \(\widetilde{\mathcal{R}} \subset \mathcal{Y}\). Thus, both of the maps \(\widetilde{\mathcal{T}}\) and \(\widetilde{\mathcal{T}}^{\times}\) are defined on \(\widetilde{\mathcal{R}}\). Therefore, for a sufficiently small \(\theta\), the map \(\widetilde{\mathcal{T}}^{\times}\) satisfies the assumptions of Theorem \ref{thm3000} and Proposition \ref{Prop7jnhiojhui22bh}. This implies that the map \(\widetilde{\mathcal{T}}\) possesses a \(\mathcal{C}^{1}\)-smooth invariant manifold
\begin{equation*}
M^{*} = \big\{\left(w, z\right): w = h^{*}\left(z\right)\big\} \subset \widetilde{\mathcal{R}},
\end{equation*}
where \(h^{*}\) is some \(\mathcal{C}^{1}\)-smooth function defined on \(\left[-\theta, \theta\right]\). Moreover, by Proposition \ref{Prop7jnhiojhui22bh}, if the backward orbit of a point in \(\widetilde{\mathcal{R}}\) remains in \(\widetilde{\mathcal{R}}\), then it must belong to \(M^{*}\). Therefore, \(\Lambda^{u}_{\widetilde{\mathcal{R}}, \widetilde{\mathcal{T}}} \subset M^{*}\). Removing the point \(\left(\frac{d}{b}, 0\right)\) from \(M^{*}\), we obtain a set which is invariant under the map \(\mathcal{T}\). Moreover, we have \(\Lambda^{u}_{\mathcal{R}, \mathcal{T}} \subset M^{*}\setminus \lbrace \left(\frac{d}{b}, 0\right)\rbrace\).

Let us now come back to \(\left(u_{10}, v_{10}\right)\)-coordinates and the Poincar\'e map \(T\). Equip \(\mathcal{R}\) with \(\left(u_{10}, v_{10}\right)\)-coordinates and choose \(0 < \epsilon < \theta\). Thus, \(\mathcal{D}_{2}^{\epsilon} \subset \mathcal{R}\). Consider the manifold \(M^{*}\) in \(\left(u_{10}, v_{10}\right)\)-coordinates and restrict it to \(\mathcal{D}_{2}^{\epsilon}\). Denote this restriction by \(\mathcal{M}^{*}\). We have that \(\mathcal{M}^{*}\setminus \lbrace M^{s}\rbrace\) is invariant under \(T\), and \(\Lambda^{u}_{\mathcal{D}_{2}^{\epsilon}, T} \subset \mathcal{M}^{*}\setminus \lbrace M^{s}\rbrace\). Choosing a sufficiently small \(\epsilon\) also guarantees that \(\mathcal{M}^{*}\) is a connected piece of \(M^{*}\) and hence is a \(\mathcal{C}^{1}\)-manifold.

The manifold \(\mathcal{M}^{*}\) is our desired manifold if we show \(\Lambda^{u}_{\mathcal{D}_{2}^{\epsilon}, T} = \mathcal{M}^{*}\setminus \lbrace M^{s}\rbrace\). So far, we have shown that \(\Lambda^{u}_{\mathcal{D}_{2}^{\epsilon}, T} \subset \mathcal{M}^{*}\setminus \lbrace M^{s}\rbrace\) and so it is sufficient to show \(\mathcal{M}^{*}\setminus \lbrace M^{s}\rbrace \subset \Lambda^{u}_{\mathcal{D}_{2}^{\epsilon}, T}\). However, this is just a direct consequence of part (\ref{Item892ib3bu8yi4y002}) of Lemma \ref{Lem899i3n9n9in3o2ubqc}  (proved earlier). The fact that \(\mathcal{M}^{*}\) is tangent to \(\ell^{*}\) at \(M^{s}\) is also a direct consequence of part (\ref{Item892ib3bu8yi4y001}) of this lemma. This ends the proof of part (\ref{Item892ib3bu8yi4y006}).
\end{proof}

\begin{myrem}\label{rem800n09b98qba9gwakl}
Lemma \ref{Lem899i3n9n9in3o2ubqc} states that when \(bd > 0\), the set \(\lbrace M^{s}\rbrace \cup\Lambda^{u}_{\mathcal{D}_{2} \cup \mathcal{D}_{3}, T}\) is a \(\mathcal{C}^{1}\)-smooth curve which is tangent to \(\ell^{*}\) at \(M^{s}\), and any point on this curve converges to \(M^{s}\) by the backward iterations of the Poincar\'e map \(T\). It follows from part (\ref{Item53iniubyv7t65czcr}) of Theorem \ref{thm3000} and the proof of Lemma \ref{Lem899i3n9n9in3o2ubqc} that, when \(bd > 0\), if we take a curve \(\zeta\) in \(\mathcal{D}_{2}\), then \(\{T^{n}\left(\zeta\right)\vert_{\mathcal{D}_{2}}\}_{n=1}^{\infty}\) converges uniformly to the curve  \(\lbrace M^{s}\rbrace \cup\Lambda^{u}_{\mathcal{D}_{2} \cup \mathcal{D}_{3}, T}\).
\end{myrem}

\subsubsection{Proof of Lemma \ref{Lem6647ububquybdu1tyqbuq}}\label{ProofofLemma22222liuo8784v}

Reverse the time direction in system (\ref{eq23000}) (i.e. \(t\rightarrow -t\)) and exchange the stable and unstable components, i.e. apply the linear change of coordinates
\begin{equation}\label{P850060}
\left(\widetilde{u}_{1}, \widetilde{u}_{2}, \widetilde{v}_{1}, \widetilde{v}_{2}\right) = \left(v_{1}, v_{2}, u_{1}, u_{2}\right)
\end{equation}
This gives a system which is of the form of system (\ref{eq23000}), where all the assumptions of Lemma \ref{Nftheorem2} are satisfied. The global map along \(\Gamma\) for this system is \(J \left(T^{\text{glo}}\right)^{-1} J^{-1}\), where \(J=\)\begin{scriptsize}
\(\left(\begin{array}{cc}
0 & 1\\
1 & 0
\end{array}\right)\)\end{scriptsize} and \(T^{\text{glo}}\) is the global map of system (\ref{eq23000}). Thus, the differential of this map at \(M^{s}\) is
\begin{equation*}
J \left(dT^{glo}\left(M^{s}\right)\right)^{-1} J^{-1} = J\cdot \frac{1}{ad-bc}\left(\begin{array}{cc}
d & -b\\
-c & a
\end{array}\right)\cdot J^{-1} = \frac{1}{ad-bc}\left(\begin{array}{cc}
a & -c\\
-b & d
\end{array}\right).
\end{equation*}
This implies that if we replace conditions \(bd>0\) and \(bd < 0\) in Lemma \ref{Lem899i3n9n9in3o2ubqc} by \(cd < 0\) and \(cd> 0\), respectively, and the line \(\ell^{*}\) by the straight line whose slope is \(\frac{-d}{c}\), then all the statements of Lemma \ref{Lem899i3n9n9in3o2ubqc} also hold for this system and the region \(\mathcal{D}_{2} \cup \mathcal{D}_{3} \subset \Pi^{s}\). Consequently, by applying the inverse of change of coordinates (\ref{P850060}), all the statements of Lemma \ref{Lem899i3n9n9in3o2ubqc} also hold for the system which is derived from system (\ref{eq23000}) by a reversion of time and the region \(\{\left(u_{1}, v_{1}\right)\in \mathcal{B}_{\epsilon_{u}}\subset \Pi^{u}:\, 0 < \frac{v_{1}}{u_{1}} \leq m, \, u_{1}\neq 0\} \subset \Pi^{u}\). In this case, the line \(\ell^{*}\) is replaced by the straight line in \(\Pi^{u}\) whose slope is \(\frac{-c}{d}\). The homoclinic loop \(\Gamma\) in this system leaves and enters \(O\) along the positive sides of \(u_{2}\) and \(v_{2}\), respectively, and the corresponding Poincar\'e map, call it \(\widetilde{T}\), is defined on \(\Pi^{u}\). Therefore, the statements of Lemma \ref{Lem899i3n9n9in3o2ubqc} also hold for the map
\begin{equation}\label{eq78001arrrnb09b08bx}
T^{\text{glo}}\circ\widetilde{T}\circ\left(T^{\text{glo}}\right)^{-1}
\end{equation}
and the set
\begin{equation}\label{eq7b8qlbwltwyiu11gpoq}
K = T^{\text{glo}} \left(\{\left(u_{1}, v_{1}\right)\in \mathcal{B}_{\epsilon_{u}}\subset \Pi^{u}:\, 0 < \frac{v_{1}}{u_{1}} \leq m, \, u_{1}\neq 0\}\right),
\end{equation}
where the line \(\ell^{*}\) is replaced by the horizontal axis in \(\Pi^{s}\). The later one is simply because 
\begin{equation*}
d T^{\text{glo}}\left(M^{u}\right)\left(\begin{array}{c}
d\\
-c
\end{array}\right) = \left(\begin{array}{c}
ad-bc\\
0
\end{array}\right).
\end{equation*}

Notice that map (\ref{eq78001arrrnb09b08bx}) is conjugate to the inverse of the Poincar\'e map \(T^{-1}\).

The map (\ref{eq78001arrrnb09b08bx}) coincides with \(T^{-1}\) on \(T\left(\mathcal{D}\right)\). Note that, for sufficiently large \(m\), the set \(T^{\text{glo}} \left(\{\left(u_{1}, v_{1}\right)\in \mathcal{B}_{\epsilon_{u}}\subset \Pi^{u}:\, u_{1} \neq 0, \, m < \frac{v_{1}}{u_{1}}\}\right)\) has no intersection with \(\mathcal{D}_{1}\). Therefore, Lemma \ref{Lem6647ububquybdu1tyqbuq} will be proved once we show that \(\mathcal{D}_{1}\subset K\) if \(cd < 0\), and \(K \cap \mathcal{D}_{1} = \emptyset\) if \(cd > 0\). However, this is an immediate consequence of the discussion above. In fact, it follows from the above discussion that the line \(\ell^{*}\) passes through \(\mathcal{D}_{2} \cup \mathcal{D}_{3}\) if and only if the horizontal axis passes through \(K\). The later case, for sufficiently large \(m\), is equivalent to the condition \(\mathcal{D}_{1}\subset K\) and happens if and only if \(cd < 0\). This ends the proof of Lemma \ref{Lem6647ububquybdu1tyqbuq}.

\begin{myrem}\label{rem78b78v6v2327872b78b}
Lemma \ref{Lem6647ububquybdu1tyqbuq} states that if \(cd < 0\), then the set \(\lbrace M^{s}\rbrace \cup\Lambda^{u}_{\mathcal{D}_{1}, T^{-1}}\) is a \(\mathcal{C}^{1}\)-smooth curve which is tangent to the horizontal axis at \(M^{s}\). Moreover, any point on this curve converges to \(M^{s}\) by the forward iterations of the Poincar\'e map \(T\). It follows from Remark \ref{rem800n09b98qba9gwakl} and the proof of Lemma \ref{Lem6647ububquybdu1tyqbuq} that if we take a curve \(\zeta\) in \(K \cap \mathcal{B}_{\epsilon}\), then \(\{T^{-n}\left(\zeta\right)\vert_{K \cap \mathcal{B}_{\epsilon}}\}_{n=1}^{\infty}\) converges uniformly to the curve \(\lbrace M^{s}\rbrace \cup\Lambda^{u}_{\mathcal{D}_{1}, T^{-1}}\).
\end{myrem}

\subsection{Dynamics near the homoclinic figure-eight}\label{Dynamicsnearfigure8homoclinic}

In this section, we study the dynamics near the homoclinic figure-eight \(\Gamma_{1} \cup \Gamma_{2}\). In particular, we prove Theorems \ref{Thmyipoinoib53}, \ref{thm89bqyvrtyvtv} and \ref{Thm6892jjdibbea} in this section. We start with recalling some definitions and notations from Section \ref{Setupsection76q5y6d45ljbiywbux}.

\begin{figure}[!htb]
\centering
\includegraphics[scale=.25]{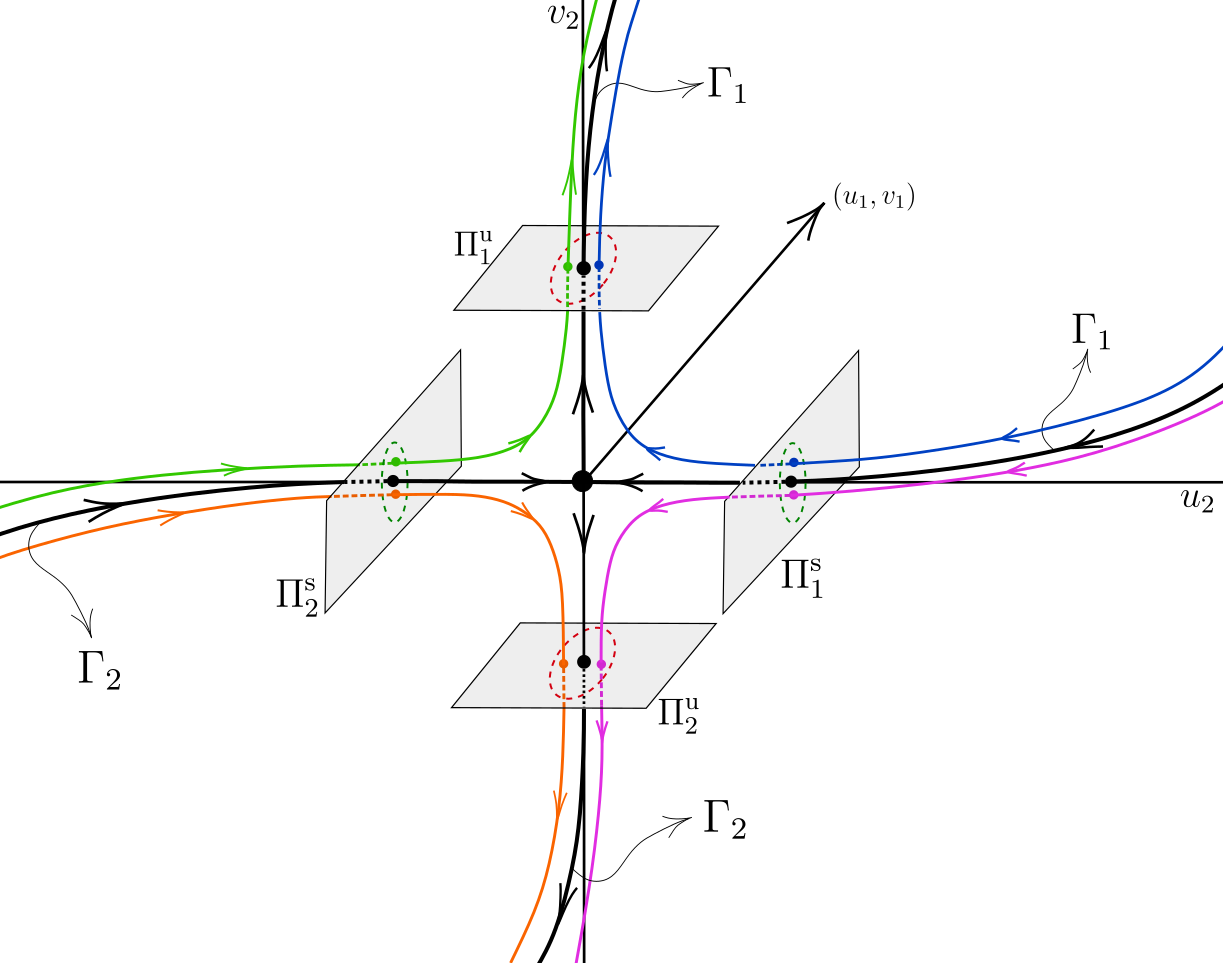}
\caption{\small The homoclinic figure-eight \(\Gamma_{1}\cup\Gamma_{2}\) and the the cross-sections \(\Pi_{1}^{s}\), \(\Pi_{1}^{u}\), \(\Pi_{2}^{s}\) and \(\Pi_{2}^{u}\) are shown. We consider \(\epsilon\)-neighborhoods of \(M_{1}^{s}\) and \(M_{2}^{s}\) (green dashed circles) in \(\Pi^{s}_{1}\) and \(\Pi^{s}_{2}\), respectively, as well as \(\epsilon_{u}\)-neighborhoods of \(M_{1}^{u}\) and \(M_{2}^{u}\) (red dashed circles) in \(\Pi^{u}_{1}\) and \(\Pi^{u}_{2}\), respectively. The set \(\mathcal{D}^{1}\) (resp. \(\mathcal{D}^{2}\)) is the set of the points in the \(\epsilon\)-neighborhood in \(\Pi^{s}_{1}\) (resp. \(\Pi^{s}_{2}\)) whose forward orbits go along \(\Gamma_{1}\) (resp. \(\Gamma_{2}\)) and intersect the \(\epsilon_{u}\)-neighborhood in \(\Pi_{1}^{u}\) (resp. \(\Pi_{2}^{u}\)). The blue point on \(\Pi^{s}_{1}\) and the brown point on \(\Pi^{s}_{2}\) belong to \(\mathcal{D}^{1}\) and \(\mathcal{D}^{2}\), respectively. We denote by \(\mathbb{D}^{1}\) (resp. \(\mathbb{D}^{2}\)) the set of the points in the \(\epsilon\)-neighborhood in \(\Pi_{1}^{s}\) (resp. \(\Pi_{2}^{s}\)) whose forward orbits go along the negative (resp. positive) side of \(v_{2}\)-axis and intersect the \(\epsilon_{u}\)-neighborhood in \(\Pi_{2}^{u}\) (resp. \(\Pi_{1}^{u}\)).}
\label{Figure00oin9n9u23723109012bygsdjn}
\end{figure}

For \(i=1, 2\), we denote by \(\mathcal{D}^{i}\) the set of the points \(\left(u_{10}, v_{10}\right)\) on \(\Pi_{i}^{s}\) whose forward orbits go along the homoclinic orbit \(\Gamma_{i}\) and intersect \(\Pi_{i}^{u}\) at \(\left(u_{1\tau}, v_{1\tau}\right)\) such that
\begin{equation}\label{conditionfordomainDsecondformula}
\left\Vert \left(u_{10}, v_{10}\right)\right\Vert < \epsilon\quad \text{and} \quad\left\Vert \left(u_{1\tau}, v_{1\tau}\right)\right\Vert < \epsilon_{u},
\end{equation}
for some sufficiently small constants \(0 < \epsilon \leq \epsilon_{u} < \delta\). We denote by \(\mathbb{D}^{1}\) (\(\mathbb{D}^{2}\)) the set of the points \(\left(u_{10}, v_{10}\right)\) on \(\Pi_{1}^{s}\) (\(\Pi_{2}^{s}\)) whose forward orbits go along the negative (positive) side of \(v_{2}\)-axis and intersect \(\Pi_{2}^{u}\) (\(\Pi_{1}^{u}\)) at \(\left(u_{1\tau}, v_{1\tau}\right)\) such that (\ref{conditionfordomainDsecondformula}) holds (see Figure \ref{Figure00oin9n9u23723109012bygsdjn}). We also denote by \(T_{i}\), \(T^{\text{loc}}_{i}\) and \(T^{\text{glo}}_{i}\) the Poincar\'e, local and global maps along \(\Gamma_{i}\) (\(i=1, 2\)), respectively (see Figure \ref{Figure673bob8i7vrq7cvraa}). The maps \(T_{1}^{\text{glo}}\) and \(T_{2}^{\text{glo}}\) are defined on the open \(\epsilon_{u}\)-balls around \(M_{1}^{s}\) and \(M_{2}^{s}\), respectively. Regarding the other maps, we have \(\text{domain}\left(T_{1}^{\text{loc}}\right) = \text{domain}\left(T_{1}\right) = \mathcal{D}^{1}\) and \(\text{domain}\left(T_{2}^{\text{loc}}\right) = \text{domain}\left(T_{2}\right) = \mathcal{D}^{2}\). We also define the map \(T^{\text{loc}}_{12}: \mathbb{D}^{1}\subset \Pi^{s}_{1} \rightarrow \Pi^{u}_{2}\) (resp. \(T^{\text{loc}}_{21}: \mathbb{D}^{2}\subset \Pi^{s}_{2} \rightarrow \Pi^{u}_{1}\)) by \(\left(u_{10}, v_{10}\right) \mapsto \left(u_{1\tau}, v_{1\tau}\right)\), where \(\left(u_{10}, v_{10}\right)\in \mathbb{D}^{1}\) (resp. \(\in \mathbb{D}^{2}\)) and \(\left(u_{1\tau}, v_{1\tau}\right)\in\Pi_{2}^{u}\) (resp. \(\in\Pi_{1}^{u}\)) (see Figure \ref{Figure673bob8i7vrq7cvraa}).

Let \(\mathcal{V}\) be a sufficiently small neighborhood of \(\Gamma_{1}\cup \Gamma_{2}\) and define \(\Xi = \mathcal{D}^{1} \cup \mathbb{D}^{1}\cup \mathcal{D}^{2} \cup \mathbb{D}^{2}\). For any \(\mathtt{x}\in \Xi\), we correspond a (finite or infinite) sequence \(\{\mathtt{x}_k\}\) to \(\mathtt{x}\) in the following way: {\it (i)} \(\mathtt{x}_{0} = \mathtt{x}\), {\it (ii)} if \(\mathtt{x}_{k}\in \Xi\) (\(k\geq 0\)), we define \(\mathtt{x}_{k+1}\) to be the first intersection point of the forward orbit of \(\mathtt{x}_{k}\) and \(\Pi^{s}_{1}\cup\Pi^{s}_{2}\). Similarly, if \(\mathtt{x}_{k}\in \Xi\) (\(k\leq 0\)), we define \(\mathtt{x}_{k-1}\) to be the first intersection point of the backward orbit of \(\mathtt{x}_{k}\) and \(\Pi^{s}_{1}\cup\Pi^{s}_{2}\). In order to understand the dynamics in \(\mathcal{V}\), we need to find the set of the points whose forward or backward orbits lie entirely in \(\mathcal{V}\), i.e. the set of the points \(\mathtt{x}\in \Xi\) for which the sequence \(\{\mathtt{x}_k\}\) is well-defined for all \(k\geq 0\) or \(k\leq 0\).

\begin{figure}
\centering
\includegraphics[scale=.19]{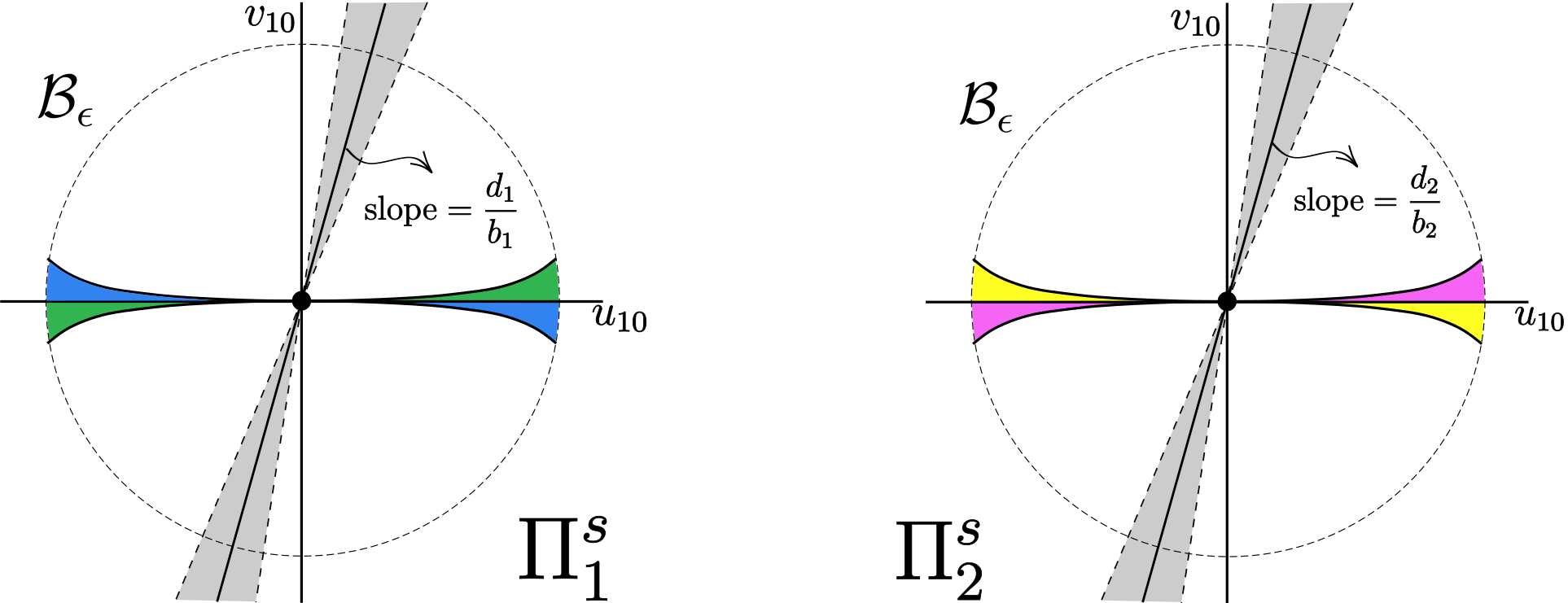}
\caption{\small This figure corresponds to the case \(\lambda_{1} < \lambda_{2} < 2\lambda_{1}\), \(b_{1}d_{1}>0\) and \(b_{2}d_{2}>0\). The regions \(\mathcal{D}^{1}\), \(\mathbb{D}^{1}\), \(\mathcal{D}^{2}\) and \(\mathbb{D}^{2}\) are shown in green, blue, pink and yellow, respectively. We define \(\Xi\) as the union of these four regions. Let \(\mathtt{x}_{1}\) be the first intersection point of the forward orbit of \(\mathtt{x} \in \Xi\) and \(\Pi_{1}^{s} \cup \Pi_{2}^{s}\). It is shown in the proof of Theorem \ref{thm89bqyvrtyvtv} that for any \(\mathtt{x}\in \Xi\) the point \(\mathtt{x}_{1}\) lies in one of the gray regions on \(\Pi_{1}^{s}\) or \(\Pi_{2}^{s}\)}
\label{Figure-yy8678679v98c786c87c}
\end{figure}

When \(\lambda_{2} < 2\lambda_{1}\), the dynamics near the homoclinic figure-eight is quite similar to the case of a single homoclinic loop: the forward and backward orbit of any arbitrary point in \(\mathcal{V}\) leaves \(\mathcal{V}\). When \(\lambda_{1}= \lambda_{2}\), it follows from Proposition \ref{Cor7uijnqnjk3o49inh44} that \(\Xi = \emptyset\) and so there is no dynamics near the homoclinic figure-eight. For the case of \(\lambda_{1} < \lambda_{2} < 2\lambda_{1}\), we show in the next proof that for any \(\mathtt{x}\in\Xi\) whose corresponding \(\mathtt{x}_{1}\) is defined, the point \(\mathtt{x}_{1}\) lies close to the straight lines with slope \(\frac{d_{1}}{b_{1}}\) (if \(\mathtt{x}_{1}\) lies in \(\mathcal{D}^{1} \cup \mathbb{D}^{1}\)) or \(\frac{d_{2}}{b_{2}}\) (if \(\mathtt{x}_{1}\) lies in \(\mathcal{D}^{2} \cup \mathbb{D}^{2}\)), and hence, it lies outside of the set \(\Xi\) (see Figure \ref{Figure-yy8678679v98c786c87c}).

\begin{proof}[Proof of Theorem \ref{thm89bqyvrtyvtv}]
The proof for the case \(\lambda_{1} = \lambda_{2}\) is an immediate consequence of Proposition \ref{Cor7uijnqnjk3o49inh44}.

Suppose \(\lambda_{1} < \lambda_{2} < 2\lambda_{1}\). By Proposition \ref{Cor8i12sehu8901fdg}, we have
\begin{equation*}
\Xi = \{\left(u_{10}, v_{10}\right):\quad \|\left(u_{10}, v_{10}\right)\| < \epsilon, \text{ and } 0 <\lvert v_{10}\rvert < K_{\epsilon_{u}} \lvert u_{10}\rvert^{\frac{\gamma}{1-\gamma}} \left[1 + O\left(\delta\right)\right]\},
\end{equation*}
where \(K_{\epsilon_{u}} > 0\) is some constant and \(\gamma = \lambda_{1}\lambda_{2}^{-1} > 0.5\) (see Figure \ref{Figure-yy8678679v98c786c87c}). Consider \(\left(u_{10}, v_{10}\right) \in \Xi\). Since \(u_{1\tau} = v_{1\tau}O\left(\epsilon^{2}\right)\) (see Remark \ref{Rem89i9ub8q7v76ce365}), the forward orbit of this point intersects one of the cross-sections \(\Pi^{u}_{1}\) or \(\Pi^{u}_{2}\) at a point close to the vertical axis and then it ends up either in the cross-section \(\Pi^{s}_{1}\) close to the straight line with the slope \(\frac{d_{1}}{b_{1}}\) or in the cross-section \(\Pi^{s}_{2}\) close to the straight line with the slope \(\frac{d_{2}}{b_{2}}\). In both cases, this point is outside of the set \(\Xi\) (see Figure \ref{Figure-yy8678679v98c786c87c}). This proves Theorem \ref{thm89bqyvrtyvtv}.
\end{proof}

A point in \(\Pi^{s}_{1} \cup \Pi^{s}_{2}\) whose forward orbit lies entirely in \(\mathcal{V}\) and does not lie on the stable manifold of \(O\) must belong to \(\Xi\). We denote the set of these points by \(\Lambda^{s}\). The same holds for backward orbits. We also define the set \(\Lambda^{u}\) analogously. In order to understand the dynamics in \(\mathcal{V}\), we need to investigate these two sets. For the case of \(\lambda_{2} < 2\lambda_{1}\), Theorem \ref{thm89bqyvrtyvtv} states that both of these sets are empty. Our approach to investigate \(\Lambda^{s}\) and \(\Lambda^{u}\) for the case \(2\lambda_{1} < \lambda_{2}\) is similar to what we have done in the previous section for the case of a single homoclinic loop.

\begin{figure}
\centering
\begin{subfigure}{0.4\textwidth}
\centering
\includegraphics[scale=0.20]{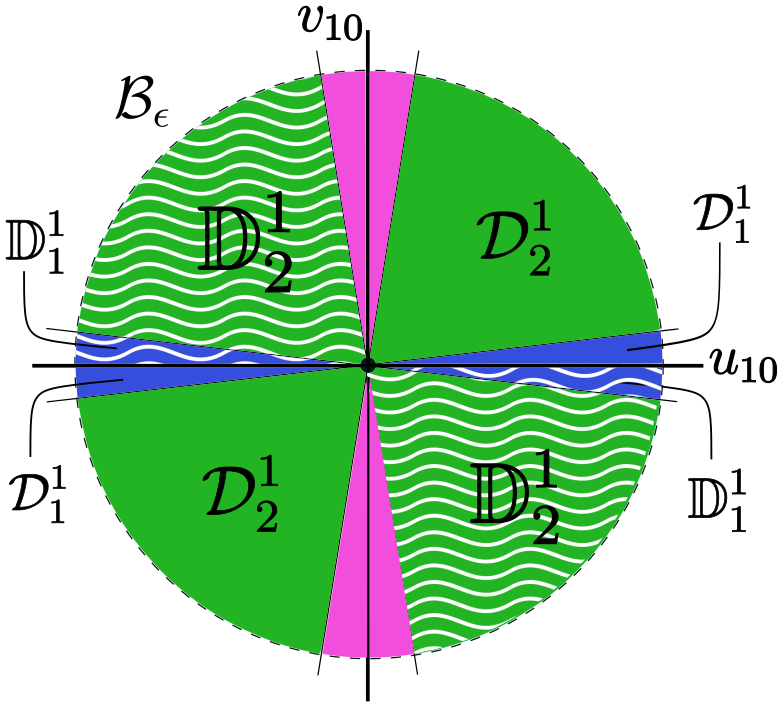}
\end{subfigure}
\begin{subfigure}{0.4\textwidth}
\centering
\includegraphics[scale=0.20]{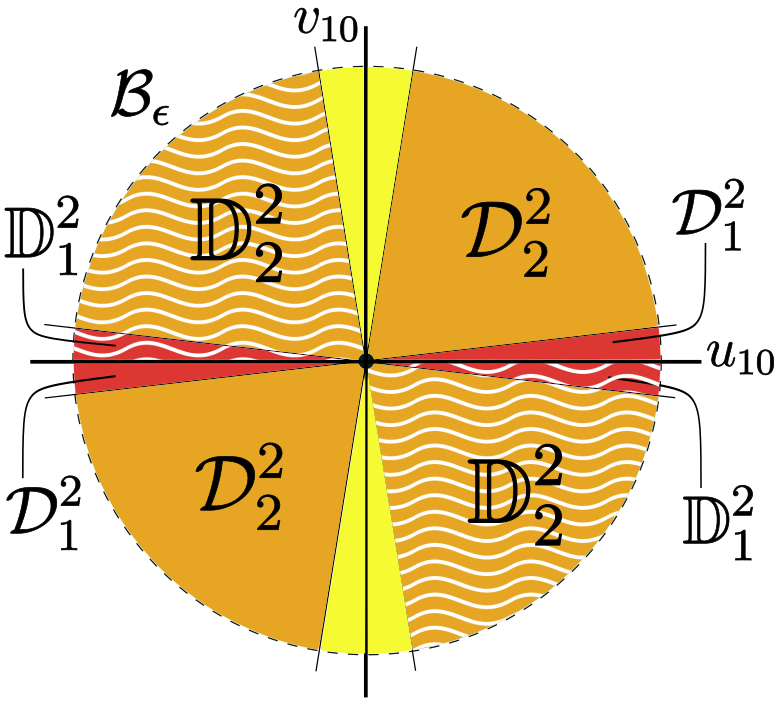}
\end{subfigure}
\caption{\small The left and right figures show \(\Pi_{1}^{s}\) and \(\Pi_{2}^{s}\), respectively. When \(2\lambda_{1} < \lambda_{2}\), we divide \(\mathcal{D}^{i}\) into three subsets \(\mathcal{D}^{i}_{1}\), \(\mathcal{D}^{i}_{2}\) and \(\mathcal{D}^{i}_{3}\) (\(i = 1, 2\)). Similarly, we divide \(\mathbb{D}^{i}\) into three subsets \(\mathbb{D}^{i}_{1}\), \(\mathbb{D}^{i}_{2}\) and \(\mathbb{D}^{i}_{3}\) (\(i = 1, 2\)). The sets \(\mathcal{D}^{1}_{3}\) and \(\mathbb{D}^{1}_{3}\) are subsets of the purple region, and the sets \(\mathcal{D}^{2}_{3}\) and \(\mathbb{D}^{2}_{3}\) are subsets of the yellow region.}
\label{Figi988n1uim3jj4rrj4uybuyv}
\end{figure}

Recall from Section \ref{Domainjyagvxjycefjthc2tc65} that when \(2\lambda_{1} < \lambda_{2}\), we divide each of the sets \(\mathcal{D}^{1}\), \(\mathbb{D}^{1}\), \(\mathcal{D}^{2}\) and \(\mathbb{D}^{2}\) into three regions, i.e. for \(i=1,2\), we write \(\mathcal{D}^{i} = \mathcal{D}^{i}_{1} \cup \mathcal{D}^{i}_{2} \cup \mathcal{D}^{i}_{3}\) and \(\mathbb{D}^{i} = \mathbb{D}^{i}_{1} \cup \mathbb{D}^{i}_{2} \cup \mathbb{D}^{i}_{3}\) (see Figure \ref{Figi988n1uim3jj4rrj4uybuyv}). Write \(\Xi = \mathcal{I} \cup \mathcal{J} \subset \Pi^{s}_{1}\cup \Pi^{s}_{2}\), where
\begin{equation*}
\mathcal{I} := \bigcup_{\substack{i =1, 2}} \left(\mathcal{D}_{1}^{i} \cup\mathbb{D}_{1}^{i}\right) \quad \text{and}
\quad \mathcal{J} := \bigcup_{\substack{i =1, 2\\ j=2,3}} \left(\mathcal{D}_{j}^{i} \cup\mathbb{D}_{j}^{i}\right).
\end{equation*}
\begin{mydefn}
We define \(\Lambda^{s}_{\mathcal{I}}\) (\(\Lambda^{u}_{\mathcal{I}}\)) as the set of the points in \(\mathcal{I}\) whose forward (backward) orbits intersect \(\Xi\) infinitely many times and all the intersection points belong to \(\mathcal{I}\). More precisely, 
\begin{equation*}
\Lambda^{s}_{\mathcal{I}} = \{\mathtt{x} = \mathtt{x}_{0}: \quad \mathtt{x}_{k}\in\mathcal{I} \text{ for all } k\geq 0\} \quad \text{and} \quad \Lambda^{u}_{\mathcal{I}} = \{\mathtt{x} = \mathtt{x}_{0}: \quad \mathtt{x}_{k}\in\mathcal{I} \text{ for all } k\leq 0\}
\end{equation*}
The sets \(\Lambda^{s}_{\mathcal{J}}\) and \(\Lambda^{u}_{\mathcal{J}}\) are defined analogously.
\end{mydefn}

Similar to the case of a single homoclinic, we take three steps to investigate the sets \(\Lambda^{s}\) and \(\Lambda^{u}\). In the first step, we investigate the sets \(\Lambda^{s}_{\mathcal{J}}\) and \(\Lambda^{u}_{\mathcal{J}}\). This is done in Lemma \ref{Lemmaliunoyb7b1878b}. From technical point of view, part (\ref{Item7783b8b7b7qqq4}) of this lemma which proves the existence of an unstable invariant manifold of the homoclinic figure-eight is the main result of this section. The techniques which are used in the proof of this part rely on the proof of part (\ref{Item892ib3bu8yi4y006}) of Lemma \ref{Lem899i3n9n9in3o2ubqc}. In the second step, we investigate the sets \(\Lambda^{s}_{\mathcal{I}}\) and \(\Lambda^{u}_{\mathcal{I}}\). This is also done in Lemma \ref{Lemouioqucrv84yrviq7y11}. Finally, in Lemma \ref{Lemi876c65x54zxyhhfhgxs}, we clarify the relations between the sets \(\Lambda^{s}_{\mathcal{J}}\), \(\Lambda^{u}_{\mathcal{J}}\), \(\Lambda^{s}_{\mathcal{I}}\) and \(\Lambda^{u}_{\mathcal{I}}\), and the sets \(\Lambda^{s}\) and \(\Lambda^{u}\). This enables us to prove Theorem \ref{thm89bqyvrtyvtv}. We start with the following:

\begin{mylem}\label{Lemmaliunoyb7b1878b}
Assume \(2\lambda_{1} < \lambda_{2}\) and let \(w\) be as in Notation \ref{Not89kookkjweer}. For \(\mathtt{x}\in \mathcal{J}\), we have
\begin{enumerate}[(i)]
\item\label{Item001pp2ok77u6iu} if \(\mathtt{x}\in \mathcal{D}^{1}_{2}\cup \mathcal{D}^{1}_{3} \cup \mathbb{D}^{2}_{2} \cup \mathbb{D}^{2}_{3}\) (i.e. \(\mathtt{x}\in \mathcal{J}\cap \Pi^{s}_{1}\)), then \(w\left(\mathtt{x}_{1}\right) = \frac{d_{1}}{b_{1}} + o(1)\). If \(\mathtt{x}\in \mathcal{D}^{2}_{2} \cup \mathcal{D}^{2}_{3} \cup \mathbb{D}^{1}_{2} \cup \mathbb{D}^{1}_{3}\) (i.e. \(\mathtt{x}\in \mathcal{J}\cap \Pi^{s}_{2}\)), then \(w\left(\mathtt{x}_{1}\right) = \frac{d_{2}}{b_{2}} + o(1)\). Here, \(o\left(1\right)\) stands for a function of \(\mathtt{x}\) that converges to zero as \(\mathtt{x}\rightarrow M^{s}_{1,2}\).
\item\label{Item22poiib97v86v37e} There exists a constant \(C>0\) such that \(\|\mathtt{x}\|^{1-2\gamma} < C \|\mathtt{x}_{1}\|\) holds for arbitrary \(\mathtt{x}\) (\(0 < \gamma = \frac{\lambda_{1}}{\lambda_{2}}<0.5\)).
\item\label{Item0991ub8uv8yrc7t} \(\mathtt{x}_{1}\in \mathcal{B}_{\epsilon}\) implies \(\mathtt{x}_{1}\in \mathcal{J}\).
\item\label{Item882786c76c11} \(\Lambda_{\mathcal{J}}^{s} = \emptyset\).
\item\label{Item7783b8b7b7qqq1} if \(b_{1}d_{1} > 0\) and \(b_{2}d_{2} < 0\), then \(\Lambda_{\mathcal{J}}^{u} = W^{u}_{\text{loc}}\left(\Gamma_{1}\right)\cap \mathcal{D}^{1}_{2}\).
\item\label{Item7783b8b7b7qqq2} if \(b_{1}d_{1} < 0\) and \(b_{2}d_{2} > 0\), then \(\Lambda_{\mathcal{J}}^{u} = W^{u}_{\text{loc}}\left(\Gamma_{2}\right)\cap \mathcal{D}^{2}_{2}\).
\item\label{Item7783b8b7b7qqq3} if \(b_{1}d_{1} > 0\) and \(b_{2}d_{2} > 0\), then \(\Lambda_{\mathcal{J}}^{u} = \left[W^{u}_{\text{loc}}\left(\Gamma_{1}\right)\cap \mathcal{D}^{1}_{2}\right] \bigcup \left[W^{u}_{\text{loc}}\left(\Gamma_{2}\right)\cap \mathcal{D}^{2}_{2}\right]\).
\item\label{Item7783b8b7b7qqq4} if \(b_{1}d_{1} < 0\) and \(b_{2}d_{2} < 0\), then \(\Lambda_{\mathcal{J}}^{u}\subset \mathbb{D}^{1}_{2}\cup\mathbb{D}^{2}_{2}\). More precisely, for each \(i = 1,2\), the union of \(M_{i}^{s}\) and \(\Lambda_{\mathcal{J}}^{u}\cap \mathbb{D}^{i}_{2}\) is a one-dimensional \(\mathcal{C}^{1}\)-manifold in \(\Pi^{s}_{i}\) which at \(M_{i}^{s}\) is tangent to the straight line with slope \(\frac{d_{i}}{b_{i}}\). Moreover, the backward orbit of any point in \(\Lambda_{\mathcal{J}}^{u}\) intersects these two manifolds alternately, i.e. for any \(\mathtt{x}\in\Lambda_{\mathcal{J}}^{u}\), all the points \(\mathtt{x}_{k}\) for even and negative \(k\)s belong to only one of the manifolds and all the other \(\mathtt{x}_{k}\) (odd and negative \(k\)s) belong to the other manifold.
\end{enumerate}
\end{mylem}

\begin{proof}
The same techniques that were used in the proof of Lemma \ref{Lem899i3n9n9in3o2ubqc} also prove parts (\ref{Item001pp2ok77u6iu}), (\ref{Item22poiib97v86v37e}) and (\ref{Item0991ub8uv8yrc7t}). 

Part (\ref{Item882786c76c11}) is an immediate consequence of (\ref{Item22poiib97v86v37e}) and (\ref{Item0991ub8uv8yrc7t}).

In the rest of the proof, we assume \(\mathtt{x}\in \Lambda_{\mathcal{J}}^{u}\). Notice that \(\mathtt{x} = \mathtt{x}_{0}\in \Lambda_{\mathcal{J}}^{u}\) implies that \(\mathtt{x}_{k}\) is defined for all \(k\leq 0\) and \(\mathtt{x}_{k}\in \Lambda_{\mathcal{J}}^{u}\). Since \(\Lambda_{\mathcal{J}}^{u}\subset \mathcal{J}\), we have two possibilities for \(\mathtt{x}_{k}\):
\begin{equation*}
\mathtt{x}_{k} \in \mathbb{D}^{2}_{2} \cup \mathbb{D}^{2}_{3} \cup \mathcal{D}^{1}_{2} \cup \mathcal{D}^{1}_{3} \qquad \text{or}\qquad
\mathtt{x}_{k} \in \mathcal{D}^{2}_{2} \cup \mathcal{D}^{2}_{3} \cup \mathbb{D}^{1}_{2} \cup \mathbb{D}^{1}_{3}.
\end{equation*}
Our strategy for proving the rest of this lemma is to consider both of these possibilities and keep track of the sequence \(\mathtt{x}_{k}, \mathtt{x}_{k+1}, \cdots, \mathtt{x}_{-1}, \mathtt{x}_{0}\). We analyze the behaviors and patterns of this sequence for arbitrary \(\mathtt{x}\in \Lambda_{\mathcal{J}}^{u}\).

To prove part (\ref{Item7783b8b7b7qqq1}), suppose \(b_{1}d_{1} > 0\) and \(b_{2}d_{2} < 0\). By part (\ref{Item001pp2ok77u6iu}), for \(\mathtt{x}_{-2}\), we observe
\begin{enumerate}[(1)]
\item \(\mathtt{x}_{-2}\in \mathbb{D}^{2}_{2} \cup \mathbb{D}^{2}_{3} \cup \mathcal{D}^{1}_{2} \cup \mathcal{D}^{1}_{3} \Longrightarrow \mathtt{x}_{-1}\in \mathcal{D}^{1}_{2}\Longrightarrow \mathtt{x}\in \mathcal{D}^{1}_{2}\), and
\vspace{-2mm}
\item \(\mathtt{x}_{-2}\in \mathcal{D}^{2}_{2} \cup \mathcal{D}^{2}_{3} \cup \mathbb{D}^{1}_{2} \cup \mathbb{D}^{1}_{3} \Longrightarrow \mathtt{x}_{-1}\in \mathbb{D}^{2}_{2}\Longrightarrow \mathtt{x}\in \mathcal{D}^{1}_{2}\).
\end{enumerate}
According to this observation, \(\mathtt{x}\in \Lambda_{\mathcal{J}}^{u}\) implies \(\mathtt{x}\in \mathcal{D}^{1}_{2}\). In other words, \(\Lambda_{\mathcal{J}}^{u}\) is in fact the set of all \(\mathtt{x}\in \mathcal{D}^{1}_{2}\) whose backward orbits only intersect \(\Pi_{1}^{s}\) (and not \(\Pi_{2}^{s}\)), and all the intersection points belong to \(\mathcal{D}^{1}_{2}\). It follows from Theorem \ref{Invariantmanifoldthm} that this set is nothing but \(W^{u}_{\text{loc}}\left(\Gamma_{1}\right)\cap \mathcal{D}^{1}_{2}\). This proves part (\ref{Item7783b8b7b7qqq1}).

The proof of part (\ref{Item7783b8b7b7qqq2}) is analogous to the proof of part (\ref{Item7783b8b7b7qqq1}).

To prove part (\ref{Item7783b8b7b7qqq3}), let \(b_{1}d_{1} > 0\) and \(b_{2}d_{2} > 0\). By (\ref{Item001pp2ok77u6iu}), for \(\mathtt{x}_{k-2}\) (\(k \leq 0\)) we observe
\begin{enumerate}[(1)]
\item \(\mathtt{x}_{k-2}\in \mathbb{D}^{2}_{2} \cup \mathbb{D}^{2}_{3} \cup \mathcal{D}^{1}_{2} \cup \mathcal{D}^{1}_{3} \Longrightarrow \mathtt{x}_{k-1}\in \mathcal{D}^{1}_{2}\Longrightarrow \mathtt{x}_{k}\in \mathcal{D}^{1}_{2}
\Longrightarrow \cdots
\Longrightarrow \mathtt{x}\in \mathcal{D}^{1}_{2}\), and
\vspace{-2mm}
\item \(\mathtt{x}_{k-2}\in \mathcal{D}^{2}_{2} \cup \mathcal{D}^{2}_{3} \cup \mathbb{D}^{1}_{2} \cup \mathbb{D}^{1}_{3}
\Longrightarrow \mathtt{x}_{k-1}\in \mathcal{D}^{2}_{2}
\Longrightarrow \mathtt{x}_{k}\in \mathcal{D}^{2}_{2}
\Longrightarrow \cdots
\Longrightarrow \mathtt{x}\in \mathcal{D}^{2}_{2}\).
\end{enumerate}
This observation holds for any arbitrary \(k\leq 0\) which means that the set \(\Lambda_{\mathcal{J}}^{u}\) consists of two disjoint sets: the first is the set of all \(\mathtt{x}\in \mathcal{D}^{1}_{2}\) whose backward orbits intersect \(\Xi\) infinitely many times and every time at \(\mathcal{D}^{1}_{2}\), and the second is the set of all \(\mathtt{x}\in \mathcal{D}^{2}_{2}\) whose backward orbits intersect \(\Xi\) infinitely many times and every time at \(\mathcal{D}^{2}_{2}\). According to Theorem \ref{Invariantmanifoldthm}, the first set is in fact \(W^{u}_{\text{loc}}\left(\Gamma_{1}\right)\cap \mathcal{D}^{1}_{2}\) and the second one is \(W^{u}_{\text{loc}}\left(\Gamma_{2}\right)\cap \mathcal{D}^{2}_{2}\). This proves part (\ref{Item7783b8b7b7qqq3}).

To prove part (\ref{Item7783b8b7b7qqq4}), let \(b_{1}d_{1} < 0\) and \(b_{2}d_{2} < 0\). By (\ref{Item001pp2ok77u6iu}), for \(\mathtt{x}_{k-1}\) (\(k \leq -4\)), we observe
\begin{enumerate}[(1)]
\item \(\mathtt{x}_{k-1}\in \mathbb{D}^{2}_{2} \cup \mathbb{D}^{2}_{3} \cup \mathcal{D}^{1}_{2} \cup \mathcal{D}^{1}_{3} \Longrightarrow \mathtt{x}_{k}\in \mathbb{D}^{1}_{2}\Longrightarrow \mathtt{x}_{k+1}\in \mathbb{D}^{2}_{2}
\Longrightarrow \mathtt{x}_{k+2}\in \mathbb{D}^{1}_{2}\\
\Longrightarrow \mathtt{x}_{k+3}\in \mathbb{D}^{2}_{2}
\Longrightarrow \cdots
\Longrightarrow \mathtt{x}\in \mathbb{D}^{1}_{2} \,\,(\text{if\,\,} \mathtt{x}_{-1}\in \mathbb{D}^{2}_{2}) \,\,\text{or}\,\, \mathtt{x}\in \mathbb{D}^{2}_{2} \,\,(\text{if\,\,} \mathtt{x}_{-1}\in \mathbb{D}^{1}_{2})\), and
\item \(\mathtt{x}_{k-1}\in \mathcal{D}^{2}_{2} \cup \mathcal{D}^{2}_{3} \cup \mathbb{D}^{1}_{2} \cup \mathbb{D}^{1}_{3}
\Longrightarrow \mathtt{x}_{k}\in \mathbb{D}^{2}_{2}\Longrightarrow \mathtt{x}_{k+1}\in \mathbb{D}^{1}_{2}
\Longrightarrow \mathtt{x}_{k+2}\in \mathbb{D}^{2}_{2}\\
\Longrightarrow \mathtt{x}_{k+3}\in \mathbb{D}^{1}_{2}
\Longrightarrow \cdots
\Longrightarrow \mathtt{x}\in \mathbb{D}^{1}_{2} \,\,(\text{if\,\,} \mathtt{x}_{-1}\in \mathbb{D}^{2}_{2}) \,\,\text{or}\,\, \mathtt{x}\in \mathbb{D}^{2}_{2} \,\,(\text{if\,\,} \mathtt{x}_{-1}\in \mathbb{D}^{1}_{2})\).
\end{enumerate}
This observation holds for any arbitrary \(k\leq -4\) and means that the backward orbit of \(\mathtt{x}\) intersects \(\mathcal{J}\) at \(\mathbb{D}_{2}^{1}\) and \(\mathbb{D}_{2}^{2}\) alternately.

Define the maps \(\mathbb{T}_{12}: \mathbb{D}^{1} \rightarrow\Pi^{s}_{2}\) and \(\mathbb{T}_{21}: \mathbb{D}^{2} \rightarrow\Pi^{s}_{1}\) by \(\mathbb{T}_{12} := T_{2} \circ T_{12}\) and \(\mathbb{T}_{21} := T_{1} \circ T_{21}\). We then define \(\mathbb{T}: \mathbb{D}^{1} \rightarrow \Pi^{s}_{1}\) by \(\mathbb{T}:= \mathbb{T}_{21} \circ \mathbb{T}_{12}\). According to the above observation, the set \(\Lambda_{\mathcal{J}}^{u}\) is in fact the set of the points \(x\in\mathbb{D}_{2}^{1}\) such that \(\mathbb{T}^{-n}\left(x\right)\in \mathbb{D}_{2}^{1}\) for all integers \(n>0\).

Recall \(\left(w, z\right)\) coordinate system and the map \(\widetilde{T}\) introduced in the proof of Lemma \ref{Lem899i3n9n9in3o2ubqc}. Similar to that proof, we equip \(\mathbb{D}^{1}_{2}\) and \(\mathbb{D}^{2}_{2}\) with \(\left(w, z\right)\) coordinates and define the maps \(\widetilde{\mathbb{T}}_{12}\) and \(\widetilde{\mathbb{T}}_{21}\) by
\begin{equation*}
\widetilde{\mathbb{T}}_{12}\left(w, z\right) := \left\{\begin{array}{ll}
\left(\overline{w}, \overline{z}\right)\quad & z \neq 0, \vspace*{2mm}\\
\left(\frac{d_{2}}{b_{2}}, 0\right)\quad & z = 0,
\end{array}\right.\qquad \text{for } \left(w, z\right) \in \widetilde{\mathcal{R}}_{1},
\end{equation*}
and
\begin{equation*}
\widetilde{\mathbb{T}}_{21}\left(w, z\right) := \left\{\begin{array}{ll}
\left(\overline{w}, \overline{z}\right)\quad & z \neq 0, \vspace*{2mm}\\
\left(\frac{d_{1}}{b_{1}}, 0\right)\quad & z = 0,
\end{array}\right.\qquad \text{for } \left(w, z\right) \in \widetilde{\mathcal{R}}_{2},
\end{equation*}
where \(\widetilde{\mathcal{R}}_{1}\) and \(\widetilde{\mathcal{R}}_{2}\) are some appropriate rectangles defined analogous to the proof of Lemma \ref{Lem899i3n9n9in3o2ubqc}. According to Remark \ref{rem90876c5xxz5ktc}, the estimates given by Lemma \ref{Derivativeflowlemma} also hold for the local maps \(T_{12}\) and \(T_{21}\). Therefore, with exactly the same proof as the proof of Lemma \ref{Lem899i3n9n9in3o2ubqc}, we see that both of the maps \(\widetilde{\mathbb{T}}_{12}\) and \(\widetilde{\mathbb{T}}_{21}\) can be written in cross-form and the partial derivatives of the cross-map satisfies the estimates given by Lemma \ref{Stat4u74i4jm1ml2mo3p}. Moreover, as it can be seen from the proof of Lemma \ref{Lem899i3n9n9in3o2ubqc}, we can make the estimates in Lemma \ref{Stat4u74i4jm1ml2mo3p} sufficiently small by choosing \(\theta\) small enough. This means that the maps \(\widetilde{\mathbb{T}}_{12}\) and \(\widetilde{\mathbb{T}}_{21}\) satisfy the assumptions of Lemma \ref{Lem22009o0ub8yvy7v21w} for sufficiently small \(K_{1}\) and \(K_{2}\). Thus, Lemma \ref{Lem22009o0ub8yvy7v21w} implies that by choosing an appropriate norm, the map \(\widetilde{\mathbb{T}}:= \widetilde{\mathbb{T}}_{21} \circ \widetilde{\mathbb{T}}_{12}\) (which is in fact the representation of \(\mathbb{T}\) in \(\left(w, z\right)\) coordinates) can be written in cross-form and the cross-map has sufficiently small partial derivatives. Therefore, this cross-map satisfies the assumptions of Theorem \ref{thm3000}. The rest of the proof follows from the proof of Lemma \ref{Lem899i3n9n9in3o2ubqc}.
\end{proof}

The following lemma is analogous to Lemma \ref{Lemmaliunoyb7b1878b}. The proof of this lemma is a simple modification of the proof of Lemma \ref{Lem6647ububquybdu1tyqbuq} for the case of homoclinic figure-eight.

\begin{mylem}\label{Lemouioqucrv84yrviq7y11}
Assume \(2\lambda_{1} < \lambda_{2}\) and let \(w\) be as in Notation \ref{Not89kookkjweer}. For \(\mathtt{x}\in \mathcal{I}\), we have
\begin{enumerate}[(i)]
\item\label{Item001pp2ok77u6iu76767} \(w\left(\mathtt{x}_{-1}\right) = o(1)\), where \(o\left(1\right)\) stands for a function of \(\mathtt{x}\) that converges to zero as \(\mathtt{x}\rightarrow M^{s}_{1,2}\).
\item\label{Item22poiib97v86v37e11} There exists a constant \(C>0\) such that \(\|\mathtt{x}\|^{1-2\gamma} < C \|\mathtt{x}_{-1}\|\) holds for any \(\mathtt{x}\) (\(0 < \gamma = \frac{\lambda_{1}}{\lambda_{2}}<0.5\)).
\item\label{Item0991ub8uv8yrc7t22} \(\mathtt{x}_{-1}\in \mathcal{B}_{\epsilon}\) implies \(\mathtt{x}_{-1}\in \mathcal{I}\).
\item\label{Item882786c76c1133} \(\Lambda_{\mathcal{I}}^{u} = \emptyset\).
\item\label{Item7783b8b7b7qqq144} if \(c_{1}d_{1} < 0\) and \(c_{2}d_{2} > 0\), then \(\Lambda_{\mathcal{I}}^{s} = W^{s}_{\text{loc}}\left(\Gamma_{1}\right)\cap \mathcal{D}^{1}_{1}\).
\item\label{Item7783b8b7b7qqq255} if \(c_{1}d_{1} > 0\) and \(c_{2}d_{2} < 0\), then \(\Lambda_{\mathcal{I}}^{s} = W^{s}_{\text{loc}}\left(\Gamma_{2}\right)\cap \mathcal{D}^{2}_{1}\).
\item\label{Item7783b8b7b7qqq366} if \(c_{1}d_{1} < 0\) and \(c_{2}d_{2} < 0\), then \(\Lambda_{\mathcal{I}}^{s} = \left[W^{u}_{\text{loc}}\left(\Gamma_{1}\right)\cap \mathcal{D}^{1}_{1}\right] \bigcup \left[W^{s}_{\text{loc}}\left(\Gamma_{2}\right)\cap \mathcal{D}^{2}_{1}\right]\).
\item\label{Item7783b8b7b7qqq477} if \(c_{1}d_{1} > 0\) and \(c_{2}d_{2} > 0\), then \(\Lambda_{\mathcal{I}}^{s}\subset \mathbb{D}^{1}_{1}\cup\mathbb{D}^{2}_{1}\). More precisely, for each \(i = 1,2\), the union of \(M_{i}^{s}\) and \(\Lambda_{\mathcal{I}}^{s}\cap \mathbb{D}^{i}_{1}\) is a one-dimensional \(\mathcal{C}^{1}\)-manifold in \(\Pi^{s}_{i}\) which at \(M_{i}^{s}\) is tangent to the horizontal axis. Moreover, the forward orbit of any point in \(\Lambda_{\mathcal{I}}^{s}\) intersects these two manifolds alternately, i.e. for any \(\mathtt{x}\in\Lambda_{\mathcal{I}}^{s}\), all the points \(\mathtt{x}_{k}\) for even and negative \(k\)s belong to only one of the manifolds and all the other \(\mathtt{x}_{k}\) (odd and negative \(k\)s) belong to the other manifold.
\end{enumerate}
\end{mylem}

The following Lemma states that the forward (resp. backward) orbit of a point in \(\mathcal{V}\) lies in \(\mathcal{V}\) if and only if it intersects the cross-sections \(\Pi_{1}^{s}\) and \(\Pi_{2}^{s}\) only at \(\mathcal{I}\) (resp. \(\mathcal{J}\)).

\begin{mylem}\label{Lemi876c65x54zxyhhfhgxs}
We have \(\Lambda^{u} = \Lambda_{\mathcal{J}}^{u}\) and \(\Lambda^{s} = \Lambda^{s}_{\mathcal{I}}\).
\end{mylem}
\begin{proof}
It follows from parts (\ref{Item22poiib97v86v37e11}) and (\ref{Item0991ub8uv8yrc7t22}) of Lemma \ref{Lemouioqucrv84yrviq7y11} that if \(\mathtt{x}\in \mathcal{I}\), then the sequence \(\lbrace\mathtt{x}_{k}\rbrace\) is not defined for all \(k\leq 0\). Indeed, For some \(k_{0}\leq 0\), we have \(\lbrace\mathtt{x}_{k_{0}}, \cdots, x_{-1}\rbrace\subset \mathcal{I}\) such that \(\mathtt{x}_{k_{0}-1}\) lies outside the \(\epsilon\)-balls around \(M^{s}_{1}\) or \(M^{s}_{2}\). This means that if \(\mathtt{x}\) belongs to \(\Lambda^{u}\), then it must belong to \(\mathcal{J}\). Therefore, \(\mathtt{x}\in \Lambda^{u}\) implies \(\mathtt{x}\in \Lambda_{\mathcal{J}}^{u}\). On the other hand, we know \(\Lambda_{\mathcal{J}}^{u}\subset\Lambda^{u}\). This proves the first part of the lemma. The proof of the other part is the same.
\end{proof}

By virtue of the preceding lemmas, we are now in a position to prove Theorem \ref{Thm6892jjdibbea}.

\begin{proof}[Proof of Theorem \ref{Thm6892jjdibbea}]
The local stable (resp. unstable) set of the homoclinic figure-eight \(\Gamma_{1}\cup\Gamma_{2}\), denoted by \(W^{s}_{\text{loc}}(\Gamma_{1}\cup\Gamma_{2})\) (resp. \(W^{u}_{\text{loc}}(\Gamma_{1}\cup\Gamma_{2})\)), is the union of \(\Gamma_{1}\cup\Gamma_{2}\) itself and the set of the points in a sufficiently small neighborhood \(\mathcal{V}\) of \(\Gamma_{1}\cup\Gamma_{2}\) whose forward (resp. backward) orbits lie in \(\mathcal{V}\) and their \(\omega\)-limit sets (resp. \(\alpha\)-limit sets) coincide with \(\Gamma_{1}\cup\Gamma_{2} \cup \lbrace O\rbrace\). By this definition, the intersection of \(W^{s}_{\text{loc}}(\Gamma_{1}\cup\Gamma_{2})\) and any of the cross-sections \(\Pi_{1}^{s}\) and \(\Pi_{2}^{s}\) must belong to \(\{M_{1}^{s}, M_{2}^{s}\}\cup\Lambda^{s}\). Similarly, the intersection of \(W^{u}_{\text{loc}}(\Gamma_{1}\cup\Gamma_{2})\) and the cross-sections \(\Pi_{1}^{s}\) and \(\Pi_{2}^{s}\) must belong to \(\{M_{1}^{s}, M_{2}^{s}\}\cup\Lambda^{u}\).

It follows from Lemma \ref{Lemouioqucrv84yrviq7y11} that in any cases except the case \(c_{1}d_{1} > 0\) and \(c_{2}d_{2} > 0\), the \(\omega\)-limit set of any orbit in \(\Lambda^{s}\) coincides with either \(\Gamma_{1} \cup \{O\}\) or \(\Gamma_{2} \cup \{O\}\). Therefore, in all of these cases, we have \(W^{s}_{\text{loc}}(\Gamma_{1}\cup\Gamma_{2}) = \Gamma_{1}\cup\Gamma_{2}\). 

Denote the flow of system (\ref{eq23000}) by \(\phi_{t}\). When \(c_{1}d_{1} > 0\) and \(c_{2}d_{2} > 0\), it follows from parts (\ref{Item22poiib97v86v37e11}) and (\ref{Item7783b8b7b7qqq477}) of Lemma \ref{Lemouioqucrv84yrviq7y11} that the set \(\Gamma_{1}\cup \Gamma_{2}\cup\phi_{t}\left(\Lambda^{s}\right)\) for \(t\geq 0\) is a 2-dimensional \(\mathcal{C}^{1}\) manifold, and the forward orbit of any point on this manifold converges to \(\Gamma_{1}\cup\{O\}\cup \Gamma_{2}\) as \(t\rightarrow \infty\). This means that this manifold is in fact the local stable set of the homoclinic figure-eight \(\Gamma_{1}\cup\Gamma_{2}\). The fact that this manifold is tangent to \(W^{s}_{\text{glo}}\left(O\right)\) at every point of \(\Gamma_{1}\cup\Gamma_{2}\) is an straightforward consequence of the discussion before Proposition \ref{Prop453ityuuuihbbdhjdndaa}.

The proof for the case of \(W^{u}_{\text{loc}}(\Gamma_{1}\cup\Gamma_{2})\) is the same. This ends the proof.
\end{proof}

\begin{mycor}\label{Corlub86276vcx6qrcyr6fc}
Let \(\phi_{t}\) be the flow of system (\ref{eq23000}). Then
\begin{equation*}
\begin{aligned}
W^{s}_{loc}\left(\Gamma_{1} \cup \Gamma_{2}\right) = \Gamma_{1} \cup \Gamma_{2} \cup \phi_{t}\left(\Lambda^{s}\right),\qquad\text{ for } t\geq 0,\\
W^{u}_{loc}\left(\Gamma_{1} \cup \Gamma_{2}\right) = \Gamma_{1} \cup \Gamma_{2} \cup \phi_{t}\left(\Lambda^{u}\right),\qquad\text{ for } t\leq 0.
\end{aligned}
\end{equation*}
\end{mycor}

Finally, we prove

\begin{proof}[Proof of Theorem \ref{Thmyipoinoib53}]
Denote the set \(W^{s}_{\text{loc}}\left(\Gamma_{1}\right) \cup W^{s}_{\text{loc}}\left(\Gamma_{2}\right) \cup W^{s}_{\text{loc}}\left(\Gamma_{1}\cup \Gamma_{2}\right)\) by \(\mathcal{W}^{s}\). By definition, the forward orbit of any point on \(W^{s}_{\mathcal{V}}\left(O\right)\cup\mathcal{W}^{s}\) lies in \(\mathcal{V}\). Consider a point in \(\mathcal{V}\setminus W_{\mathcal{V}}^{s}\left(O\right)\) whose forward orbit lies entirely in \(\mathcal{V}\). The forward orbit of this point must intersect \(\Pi_{1}^{s}\cup\Pi_{2}^{s}\) at \(\Lambda^{s}\). Therefore, it follows from the proof of Theorem \ref{thm89bqyvrtyvtv} (for the case \(\lambda_{2}<2\lambda_{1}\)) and Corollary \ref{Corlub86276vcx6qrcyr6fc} (for the case \(2\lambda_{1} < \lambda_{2}\)) that this point lies on \(\mathcal{W}^{s}\). This finishes the proof for the case of forward orbits.

The proof of the case of backward orbits is the same. This finishes the proof.
\end{proof}

\subsection{Dynamics near super-homoclinic orbits}\label{Superhomoclinicsubsection}

\begin{figure}[!tbh]
\centering
\includegraphics[scale=.24]{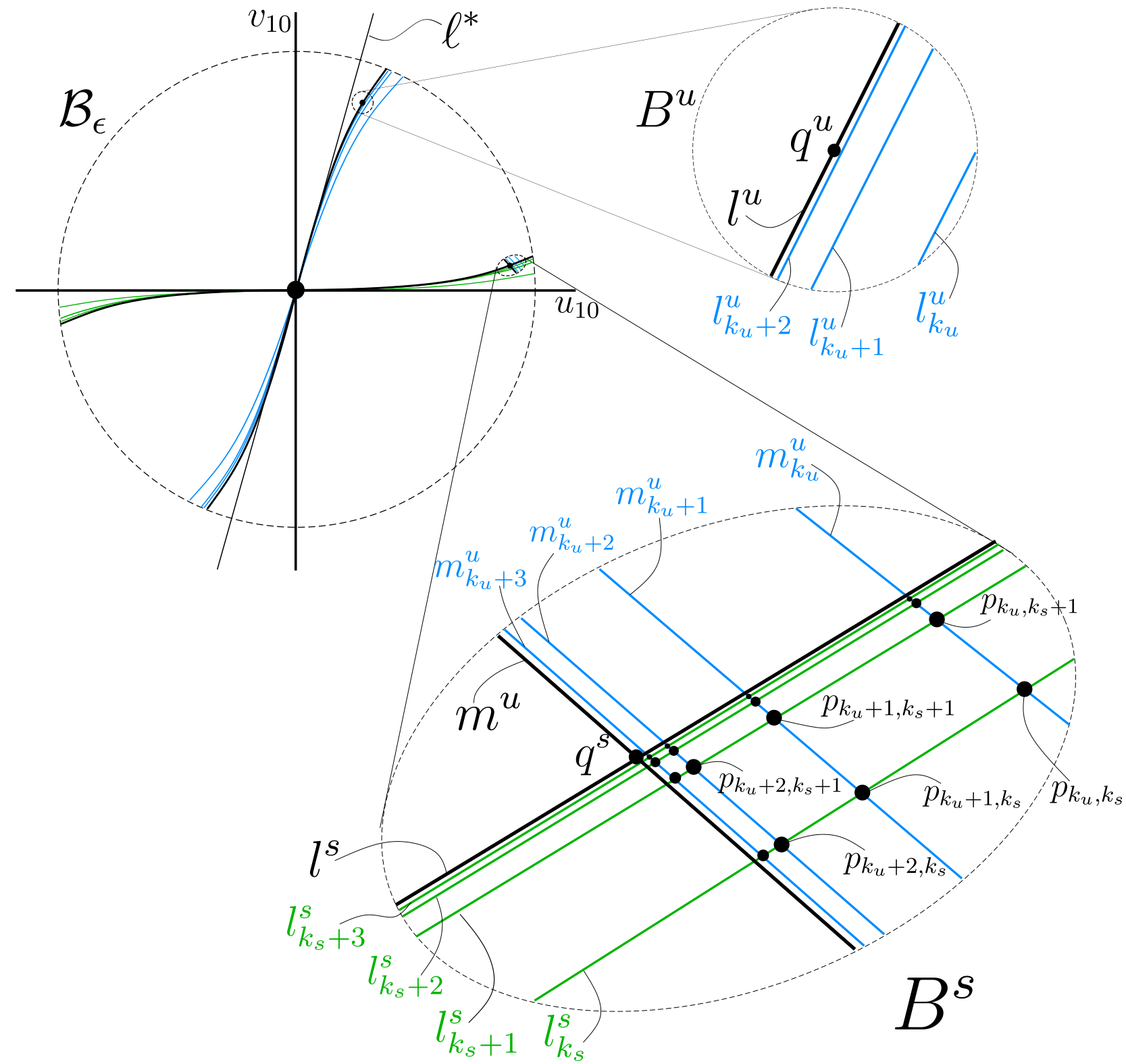}
\caption{\small The blue and green curves belong to the intersection of \(\Pi^{s}\) and the global unstable and stable invariant manifolds of the equilibrium \(O\), respectively. The blue curves accumulate to \(\mathcal{W}^{u}\) and the green curves accumulate to \(\mathcal{W}^{s}\), where \(\mathcal{W}^{u} = W^{u}_{\text{loc}}\left(\Gamma\right) \cap \Pi^{s}\) and \(\mathcal{W}^{s} = W^{s}_{\text{loc}}\left(\Gamma\right) \cap \Pi^{s}\). Let \(q^{u}\in \mathcal{W}^{u}\) and \(q^{s}\in\mathcal{W}^{s}\) be at the intersection of the super-homoclinic orbit and the cross-section \(\Pi^{s}\). The flow near the super-homoclinic orbit defines a map on a small neighborhood \(B^{u}\) of \(q^{u}\) onto a small neighborhood \(B^{s}\) of \(q^{s}\). This map maps the blue curves restricted to \(B^{u}\) to the blue curves in \(B^{s}\). The blue and green curves in \(B^{s}\) intersect transversely. Any point of these intersections belongs to both stable and unstable invariant manifolds of \(O\). Thus, the orbits passing through these points are homoclinic to \(O\).}
\label{Figure41414vi7v7tsrucy}
\end{figure}

In this section, we prove Theorem \ref{superhomoclinicthm}. The idea of the proof is to show that there exist sequences of curves \(\{l^{u}_{k}\}_{k=1}^{\infty} \subset W^{u}_{\text{glo}}\left(O\right) \cap \Pi^{s}\) and \(\{l^{s}_{k}\}_{k=1}^{\infty} \subset W^{s}_{\text{glo}}\left(O\right) \cap \Pi^{s}\) that accumulate to \(W^{u}_{\text{loc}}\left(\Gamma\right) \cap \Pi^{s}\) and \(W^{s}_{\text{loc}}\left(\Gamma\right) \cap \Pi^{s}\), respectively (see Figure \ref{Figure41414vi7v7tsrucy}). Then, the flow near the super-homoclinic orbit defines a map which maps the first sequence to a sequence of curves, denoted by \(\{m^{u}_{k}\}\) in Figure \ref{Figure41414vi7v7tsrucy}, such that each of the curves \(\{m^{u}_{k}\}\) intersects each of the curves \(\{l^{s}_{k}\}\) at a single point. Each of these intersection points correspond to a homoclinic orbit. The proof of Theorem \ref{superhomoclinicthmforfigure8} is exactly the same.

\begin{proof}[Proof of Theorem \ref{superhomoclinicthm}]
Let \(\mathcal{W}^{s} = W^{s}_{\text{loc}}\left(\Gamma\right) \cap \mathcal{D}_{1}\) and \(\mathcal{W}^{u} = W^{u}_{\text{loc}}\left(\Gamma\right) \cap \mathcal{D}_{2}\). We have shown in Section \ref{kjhliub812ftf878g3x7c} (after Proposition \ref{Cor12byeijheuije}) that \(T^{\text{glo}} \left(W^{u}_{\text{loc}}\left(O\right)\cap \Pi^{u}\right)\) intersects \(\Pi^{s}\) at a curve which is tangent to \(\ell^{*}\) at \(M^{s}\). For a sufficiently small \(\epsilon\), the restriction of this curve to \(\mathcal{B}_{\epsilon}\setminus\{M^{s}\}\) lies in \(\mathcal{D}_{2}\). Denote this restricted curve by \(L^{u}_{0}\), and let \(L^{u}_{k}\) (\(k\geq 1\)) be the restriction of \(T\left(L^{u}_{k-1}\right)\) to \(\mathcal{B}_{\epsilon}\setminus\{M^{s}\}\). By Remark \ref{rem800n09b98qba9gwakl}, the sequence \(\{L^{u}_{k}\}_{k=1}^{\infty}\) converges to \(\mathcal{W}^{u}\) uniformly.

Now, consider the restriction of \(W^{s}_{\text{loc}}\left(O\right) \cap \Pi^{s}\) to \(\mathcal{B}_{\epsilon}\setminus\{M^{s}\}\) and denote it by \(L^{s}_{0}\). We have \(L_{0}^{s}\subset K\), where \(K\) is as in (\ref{eq7b8qlbwltwyiu11gpoq}). Let \(L^{s}_{k}\) (\(k\geq 1\)) be the restriction of \(T^{-1}\left(L^{s}_{k-1}\right)\) to \(\mathcal{B}_{\epsilon}\setminus\{M^{s}\}\). By Remark \ref{rem78b78v6v2327872b78b}, the sequence \(\{L^{s}_{k}\}_{k=1}^{\infty}\) converges to \(\mathcal{W}^{s}\) uniformly.

The super-homoclinic orbit \(\mathcal{S}\) intersects \(\Pi^{s}\) at \(\mathcal{W}^{u}\) and \(\mathcal{W}^{s}\) infinitely many times. Denote  the furthest points of \(\mathcal{S} \cap \mathcal{W}^{u}\) and \( \mathcal{S} \cap \mathcal{W}^{s}\) from \(M^{s}\) by \(q^{u}\) and \(q^{s}\), respectively. Let \(B^{u}\) be a sufficiently small open ball in \(\mathcal{D}_{2}\) centered at \(q^{u}\). The orbits starting from \(B^{u}\) leave the small neighborhood \(\mathcal{U}\) of \(\Gamma\) and go along the super-homoclinic orbit \(\mathcal{S}\), and after a finite time, they come back and intersect \(\Pi^{s}\) at some points close to \(q^{s}\). These orbits induce a global map
\begin{equation*}
T_{\mathcal{S}}: B^{u} \subset \Pi^{s} \rightarrow B^{s}\subset \Pi^{s}
\end{equation*}
along the super-homoclinic orbit \(\mathcal{S}\), where \(B^{s} = T_{\mathcal{S}} \left(B^{u}\right)\) and \(T_{\mathcal{S}}\left(q^{u}\right) = q^{s}\). Since \(B^{u}\) is sufficiently small and the map \(T_{\mathcal{S}}\) is a diffeomorphism, the neighborhood \(B^{s}\) is small, connected and convex.

Define \(l^{u} = \mathcal{W}^{u} \cap B^{u}\) and \(l^{s} = \mathcal{W}^{s} \cap B^{s}\). Since the sequence \(\{L^{s}_{k}\}_{k=1}^{\infty}\) converges to \(\mathcal{W}^{s}\) uniformly, there exists a sufficiently large \(k_{s}\) such that for all \(k\geq k_{s}\), the curve \(L^{s}_{k}\) intersects \(B^{s}\). Let \(l^{s}_{k} = L^{s}_{k}\cap B^{s}\) for \(k\geq k_{s}\). This implies that \(l^{s}_{k}\xrightarrow{\text{unif}} l^{s}\). Similarly, for some sufficiently large \(k_{u}\), all the curves \(L^{u}_{k}\) for \(k\geq k_{u}\) intersect \(B^{u}\). Let \(l^{u}_{k} = L^{u}_{k}\cap B^{u}\) for \(k\geq k_{u}\). Therefore, \(l^{u}_{k}\xrightarrow{\text{unif}} l^{u}\).

The map \(T_{\mathcal{S}}\) maps \(B^{u}\) to \(B^{s}\). Thus, the curves \(l^{u}\) and \(l^{u}_{k}\) in \(B^{u}\) are mapped to some curves in \(B^{s}\) by \(T_{\mathcal{S}}\). Let \(m^{u} = T_{\mathcal{S}}\left(l^{u}\right)\) and \(m^{u}_{k} = T_{\mathcal{S}}\left(l^{u}_{k}\right)\) for \(k\geq k_{u}\). Since the super-homoclinic orbit \(\mathcal{S}\) is at the transverse intersection of the stable and unstable invariant manifolds of the homoclinic orbit \(\Gamma\), the curves \(m^{u}\) and \(l^{s}\) intersect each other transversely. On the other hand, the sequences of the curves \(m^{u}_{k}\) and \(l^{s}_{k}\) converge to \(m^{u}\) and \(l^s\), respectively. This implies that the curves \(m^{u}_{k}\) intersect the curves \(l^{s}_{k}\) transversely. Moreover, without loss of generality, we can assume that the integers \(k_{u}\) and \(k_{s}\) are large enough such that the curves \(m^{u}_{i}\) and \(l^{s}_{j}\) intersect each other at a unique point \(p_{i,j}\) for any \(i\geq k_{u}\) and \(j\geq k_{s}\). The orbits passing through the points \(p_{i,j}\) are the desired multi-pulse homoclinic orbits. This proves Theorem \ref{superhomoclinicthm}.
\end{proof}

\appendix

\addtocontents{toc}{\protect\setcounter{tocdepth}{1}}
% This command control what appears on the table of content.

\section{Proofs of Lemmas \ref{Nftheorem0}, \ref{Nftheorem} and \ref{Nftheorem2}}\label{Proofs of normal form Lemmas}

In this appendix, we prove normal form lemmas. We start with a brief discussion on some materials needed for the proofs of the lemmas, and then proceed to the proofs.

\subsection{Preliminaries}\label{Invariantmanifolds}

Consider the system
\begin{equation}\label{eqkuhufqpz00ruyf7y7348w7tv874}
\begin{aligned}
\dot{x} &= f\left(x,y\right),\\ 
\dot{y} &= g\left(x,y\right),
\end{aligned}
\end{equation}
where \(x\in \mathbb{R}^{m}\), \(y \in \mathbb{R}^{n}\) and \(f\left(0,0\right) = g\left(0,0\right) = 0\). Let \(\varphi: \mathbb{R}^{m} \rightarrow \mathbb{R}^{n}\) be a smooth mapping such that \(\varphi\left(0\right) = 0\) and \(\varphi^{\prime}\left(0\right) = 0\). Let the manifold \(\mathcal{M} = \lbrace \left(x, y\right): y = \varphi\left(x\right)\rbrace\) be invariant with respect to the flow of this system.
\begin{mydefn}
By straightening the invariant manifold \(\mathcal{M}\), we mean applying a change of coordinates of the form \(\left(\widetilde{x}, \widetilde{y}\right) = \left(x, y - \varphi\left(x\right)\right)\).
\end{mydefn}
Making this change of coordinates transforms the manifold \(\mathcal{M}\) to the linear subspace \(\lbrace \left(x, y\right): y = 0\rbrace\). Straightening an invariant manifold of the type \(\lbrace \left(x, y\right): x = \varphi\left(y\right)\rbrace\), where \(\varphi: \mathbb{R}^{n} \rightarrow \mathbb{R}^{m}\) is a smooth mapping such that \(\varphi\left(0\right) = 0\) and \(\varphi^{\prime}\left(0\right) = 0\), is defined analogously.

Again, consider system (\ref{eqkuhufqpz00ruyf7y7348w7tv874}) and let \(\varphi: \mathbb{R}^{m} \rightarrow \mathbb{R}^{n}\) and \(\psi: \mathbb{R}^{n} \rightarrow \mathbb{R}^{m}\) be some smooth maps such that \(\varphi\left(0\right) = 0\) and \(\psi\left(0\right) = 0\). According to \cite{BakraniPhDthesis}, we have
\begin{myprop}\label{invariancecondition}
The manifold \(\mathcal{M} = \lbrace \left(x, y\right): y = \varphi\left(x\right)\rbrace\) is invariant with respect to the flow of system (\ref{eqkuhufqpz00ruyf7y7348w7tv874}) if and only if
\begin{equation}\label{eq9uaaas8u843c75c83745b}
g\left(x,\varphi\left(x\right)\right) = \varphi^{\prime}\left(x\right)\cdot f\left(x,\varphi\left(x\right)\right).
\end{equation}
Analogously, the manifold \(\mathcal{N} = \lbrace \left(x,y\right): x=\psi\left(y\right)\rbrace\) is invariant with respect to the flow of system (\ref{eqkuhufqpz00ruyf7y7348w7tv874}) if and only if
\begin{equation}\label{eqb8b8baky10opqkniyfgubyfb}
f\left(\psi\left(y\right), y\right) = \psi^{\prime}\left(y\right)\cdot g\left(\psi\left(y\right), y\right).
\end{equation}
\end{myprop}

\begin{mydefn}\label{invarianceconditiondefinition}
We refer to (\ref{eq9uaaas8u843c75c83745b}) (resp. (\ref{eqb8b8baky10opqkniyfgubyfb})) as the condition of the invariance of the manifold \(\mathcal{M}\) (resp. \(\mathcal{N}\)) with respect to the flow of system (\ref{eqkuhufqpz00ruyf7y7348w7tv874}).
\end{mydefn}

\begin{myrem}\label{Cor99090q19n9be8u8y3f}
System (\ref{eqkuhufqpz00ruyf7y7348w7tv874}) has an equilibrium state at the origin. This equilibrium may possess strong stable \(W^{ss}\left(O\right)\), strong unstable \(W^{uu}\left(O\right)\), extended stable \(W^{sE}\left(O\right)\) and extended unstable invariant manifolds \(W^{uE}\left(O\right)\) (see \cite{Dimabook}). Let system (\ref{eqkuhufqpz00ruyf7y7348w7tv874}) be invariant with respect to some linear symmetry. Then, changes of coordinates that straighten the manifolds \(W^{ss}\left(O\right)\), \(W^{uu}\left(O\right)\), \(W^{sE}\left(O\right)\) and \(W^{uE}\left(O\right)\) commute with that symmetry.
\end{myrem}

\subsection{Proofs of Lemmas \ref{Nftheorem0}, \ref{Nftheorem} and \ref{Nftheorem2}}

\begin{proof}[Proof of Lemma \ref{Nftheorem0}]
To reduce system (\ref{eq300}) to the form (\ref{eq74wt366cuw6gt5uw6v01092}), we straighten the local stable and local unstable invariant manifolds of the equilibrium state \(O\), i.e. we apply a change of coordinates
\begin{equation}\label{eq0nu9wrub4c8wb8i8y275r6}
\begin{aligned}
\tilde{u}_{1} &=u_{1}-\varphi_{1 s}(v_{1}, v_{2}),\qquad
&&\tilde{u}_{2} =u_{2}-\varphi_{2 s}(v_{1}, v_{2}),\\
\tilde{v}_{1} &=v_{1}-\psi_{1 u}(u_{1}, u_{2}),\qquad
&&\tilde{v}_{2} =v_{2}-\psi_{2 u}(u_{1}, u_{2}),
\end{aligned}
\end{equation}
where \(\lbrace u_{1}=\varphi_{1 s}(v_{1}, v_{2}), u_{2}=\varphi_{2 s}(v_{1}, v_{2})\rbrace\) and \(\lbrace v_{1}=\psi_{1 u}(u_{1}, u_{2}), v_{2}=\psi_{2 u}(u_{1}, u_{2})\rbrace\) are the equations of the local stable and the local unstable invariant manifolds of \(O\), respectively. Thus, after applying (\ref{eq0nu9wrub4c8wb8i8y275r6}), the equations of the local stable and the local unstable manifolds of \(O\) become \(\lbrace v_{1} = v_{2} = 0\rbrace\) and \(\lbrace u_{1} = u_{2} = 0\rbrace\), respectively. This implies that system (\ref{eq300}) can be written in the form (\ref{eq74wt366cuw6gt5uw6v01092}) such that (\ref{equnsfuhnrx9uwn984urb348}) is satisfied. Notice that change of coordinates (\ref{eq0nu9wrub4c8wb8i8y275r6}) does not affect the quadratic part of (\ref{eq400}). Therefore, the updated first integral \(H\) keeps the form (\ref{eq400}).

Since \(H\) vanishes at every point of the local unstable invariant manifold \(\lbrace u_{1} = u_{2} = 0\rbrace\), it can be written as
\begin{equation}\label{eq0qi403n9ubv8b5v8b4}
H\left(u_{1}, u_{2}, v_{1}, v_{2}\right) = \lambda_{1}u_{1} \left[v_{1} + H_{1}\left(u_{1}, u_{2}, v_{1}, v_{2}\right)\right] - \lambda_{2}u_{2} \left[ v_{2} + H_{2}\left(u_{1}, u_{2}, v_{1}, v_{2}\right)\right],
\end{equation}
for some \(\mathcal{C}^{\infty}\)-smooth \(H_{1}, H_{2}: \mathbb{R}^{4} \rightarrow \mathbb{R}\) such that \(H_{1}\) and \(H_{2}\) and their first derivatives vanish at \(O\). On the other hand, \(H\) vanishes at every point of the local stable invariant manifold \(\lbrace v_{1} = v_{2} = 0\rbrace\). This implies
\begin{equation*}
\begin{aligned}
0 = H\left(u_{1}, u_{2}, 0, 0\right) = \lambda_{1}u_{1} H_{1}\left(u_{1}, u_{2},0, 0\right) - \lambda_{2}u_{2} H_{2}\left(u_{1}, u_{2}, 0, 0\right).
\end{aligned}
\end{equation*}
Therefore
\begin{equation*}
\begin{aligned}
H\left(u_{1}, u_{2}, v_{1}, v_{2}\right) =& \lambda_{1}u_{1} \left[v_{1} + H_{1}\left(u_{1}, u_{2}, v_{1}, v_{2}\right) - H_{1}\left(u_{1}, u_{2},0, 0\right)\right]\\
& - \lambda_{2}u_{2} \left[ v_{2} + H_{2}\left(u_{1}, u_{2}, v_{1}, v_{2}\right) - H_{2}\left(u_{1}, u_{2},0, 0\right)\right].
\end{aligned}
\end{equation*}
This suggests that, without loss of generality, we can assume that \(H_{1}\) and \(H_{2}\) vanish at \(\lbrace v_{1} = v_{2} = 0\rbrace\). Now, consider the change of coordinates
\begin{equation}\label{eq9u4nsr8uiubrw48yb4et8}
\begin{aligned}
\tilde{u}_{1} &= u_{1},\quad
&&\tilde{u}_{2} = u_{2},\\
\tilde{v}_{1} &= v_{1} + H_{1}\left(u_{1}, u_{2}, v_{1}, v_{2}\right),\quad
&&\tilde{v}_{2} = v_{2} + H_{2}\left(u_{1}, u_{2}, v_{1}, v_{2}\right).
\end{aligned}
\end{equation}
Since \(H_{1}\left(u_{1}, u_{2},0, 0\right) = H_{2}\left(u_{1}, u_{2},0, 0\right) = 0\), applying this change of coordinates on system (\ref{eq74wt366cuw6gt5uw6v01092}) keeps the local stable and local unstable invariant manifolds straightened and therefore keeps the form (\ref{eq74wt366cuw6gt5uw6v01092}) of the system such that (\ref{equnsfuhnrx9uwn984urb348}) still holds. However, this change of coordinates reduces the first integral \(H\) to the form (\ref{eq98n9uo8wbr8bau4b819ub92q1u8}).

It is a direct consequence of Remark \ref{Cor99090q19n9be8u8y3f} that change of coordinates (\ref{eq0nu9wrub4c8wb8i8y275r6}) preserves the symmetric structure of the system and the first integral. Concerning the change of coordinates (\ref{eq9u4nsr8uiubrw48yb4et8}), note that since \(H\) in (\ref{eq0qi403n9ubv8b5v8b4}) satisfies (\ref{eq426}), we have
\begin{equation*}
\begin{aligned}
H_{1}\left(-u_{1}, u_{2}, -v_{1}, v_{2}\right) &= - H_{1}\left(u_{1}, u_{2}, v_{1}, v_{2}\right),\\ H_{2}\left(-u_{1}, u_{2}, -v_{1}, v_{2}\right) &= H_{2}\left(u_{1}, u_{2}, v_{1}, v_{2}\right).
\end{aligned}
\end{equation*}
This implies that the change of coordinates (\ref{eq9u4nsr8uiubrw48yb4et8}) commutes with symmetry (\ref{eq425}), and therefore, preserves the invariance of the system with respect to symmetry (\ref{eq425}). This ends the proof of Lemma \ref{Nftheorem0}.
\end{proof}

Our proof of Lemma \ref{Nftheorem} is based on a theorem in \cite{Dimabook} (Theorem A.1). A special case of this theorem that we need for the proof of that lemma is stated below:

\begin{mylem}\label{Dimabooknormalform}(\cite{Dimabook}, Theorem A.1)
Consider system (\ref{eq74wt366cuw6gt5uw6v01092}) and assume \(\lambda_{1} < \lambda_{2}\). There exists a \({\mathcal{C}}^{\infty}\)-smooth change of coordinates which brings system (\ref{eq74wt366cuw6gt5uw6v01092}) to the form
\begin{equation}\label{eq0w9b38b47v8hbub8r57qw3t6n}
\begin{aligned}
\dot{u}_{1} &= -\lambda_{1}u_{1} + f_{11}(u_{1}, u_{2}, v_{1}, v_{2}) u_{1} + f_{12}(u_{1}, u_{2}, v_{1}, v_{2}) u_{2},\\
\dot{u}_{2} &= -\lambda_{2}u_{2} + f_{21}(u_{1}, u_{2}, v_{1}, v_{2}) u_{1} + f_{22}(u_{1}, u_{2}, v_{1}, v_{2}) u_{2},\\
\dot{v}_{1} &= +\lambda_{1}v_{1} + g_{11}(u_{1}, u_{2}, v_{1}, v_{2}) v_{1} + g_{12}(u_{1}, u_{2}, v_{1}, v_{2}) v_{2},\\
\dot{v}_{2} &= +\lambda_{2}v_{2} + g_{21}(u_{1}, u_{2}, v_{1}, v_{2}) v_{1} + g_{22}(u_{1}, u_{2}, v_{1}, v_{2}) v_{2},
\end{aligned}
\end{equation}
where the functions \(f_{ij}\), \(g_{ij}\) are \({\mathcal{C}}^{\infty}\)-smooth and
\begin{equation}\label{eq84un8byyyyybw7tttttttttyb7}
\begin{gathered}
f_{ij}\left(0,0,0,0\right) = 0,\quad f_{1i}\left(u_{1}, u_{2}, 0, 0\right) \equiv 0,\quad f_{j1}\left(0, 0, v_{1}, v_{2}\right) \equiv 0,\\
g_{ij}\left(0,0,0,0\right) = 0,\quad g_{1i}\left(0, 0, v_{1}, v_{2}\right) \equiv 0,\quad  g_{j1}\left(u_{1}, u_{2}, 0, 0\right) \equiv 0,\quad (i,j = 1,2).
\end{gathered}
\end{equation}
\end{mylem}

\begin{proof}
See \cite{Dimabook}.
\end{proof}

As a matter of comparison between this lemma and Lemma \ref{Nftheorem}, the functions \(f_{i1}\) and \(g_{i1}\) (\(i=1, 2\)) in (\ref{eq19000}) do not depend on  \(u_{2}\) and \(v_{2}\), respectively, and (\ref{eq20000}) includes all conditions (\ref{eq84un8byyyyybw7tttttttttyb7}) as well as two extra constraints
\begin{align}
f_{22}(0, v)\equiv 0,\label{eq28000}\\
g_{22}(u, 0)\equiv 0.\label{eq29000}
\end{align}

\begin{myrem}\label{rem8a7vcyut4l7acvtk3iyulvyt7w}
The desired change of coordinates in Lemma \ref{Dimabooknormalform}, denote it by \(\Phi\), is in fact a composition of several changes of coordinates, each describing some invariant manifolds. One can observe that each of these changes of coordinates commutes with symmetry (\ref{eq425}) (see \cite{BakraniPhDthesis}). Moreover, due to this symmetric property, each of these changes of coordinates can be written in the form
\begin{equation}\label{eq137010}
\begin{array}{ll}
\tilde{u}_{1} = u_{1}\left[1 + o\left(1\right)\right],\quad & \quad\tilde{u}_{2} = u_{2}\left[1 + o\left(1\right)\right] + u_{1}v_{1}O\left(1\right),\\
\tilde{v}_{1} = v_{1}\left[1 + o\left(1\right)\right],\quad & \quad\tilde{v}_{2} = v_{2}\left[1 + o\left(1\right)\right] + u_{1}v_{1}O\left(1\right),
\end{array}
\end{equation}
where \(o\left(1\right)\) and \(O\left(1\right)\) stand for \(\mathcal{C}^{\infty}\)-smooth functions of \(\left(u, v\right)\) which converge to zero  and are bounded above by a constant, respectively, as \(\left(u, v\right)\rightarrow O\). A straightforward calculation shows that making changes of coordinates of this form preserves the form (\ref{eq98n9uo8wbr8bau4b819ub92q1u8}) of the first integral \(H\). On the other hand, first integral (\ref{eq98n9uo8wbr8bau4b819ub92q1u8}) is already of the form (\ref{eq22000}). This implies that the change of coordinates \(\Phi\) transforms first integral (\ref{eq98n9uo8wbr8bau4b819ub92q1u8}) to the form (\ref{eq22000}).
\end{myrem}

\begin{proof}[Proof of Lemma \ref{Nftheorem}]
According to Lemmas \ref{Nftheorem0} and \ref{Dimabooknormalform}, there exists a change of coordinates which brings system (\ref{eq300}) to system (\ref{eq0w9b38b47v8hbub8r57qw3t6n}) where the functions \(f_{ij}\), \(g_{ij}\) are \({\mathcal{C}}^{\infty}\)-smooth and satisfy (\ref{eq84un8byyyyybw7tttttttttyb7}). We show that there exists a change of coordinates which brings system (\ref{eq0w9b38b47v8hbub8r57qw3t6n}) into the form (\ref{eq19000}), where \(f_{ij}\), \(g_{ij}\) are \({\mathcal{C}}^{\infty}\)-smooth and satisfy (\ref{eq20000}).

Consider system (\ref{eq0w9b38b47v8hbub8r57qw3t6n}) and, for \(i=1,2\), let
\begin{equation*}
\begin{gathered}
f_{i1}^{\text{new}}(u_{1}, v) = f_{i1}(u_{1}, 0, v_{1}, v_{2}),\qquad g_{i1}^{\text{new}}(u, v_{1}) = g_{i1}(u_{1}, u_{2}, v_{1}, 0),\\
f_{i2}^{\text{new}}(u_{1}, u_{2}, v) = \left[\frac{f_{i1}(u_{1}, u_{2}, v_{1}, v_{2}) - f_{i1}(u_{1}, 0, v_{1}, v_{2})}{u_{2}}\right] u_{1} + f_{i2}(u_{1}, u_{2}, v_{1}, v_{2}),\\
g_{i2}^{\text{new}}(u, v_{1}, v_{2}) = \left[\frac{g_{i1}(u_{1}, u_{2}, v_{1}, v_{2}) - g_{i1}(u_{1}, u_{2}, v_{1}, 0)}{v_{2}}\right] v_{1} + g_{i2}(u_{1}, u_{2}, v_{1}, v_{2}),
\end{gathered}
\end{equation*}
It is easily seen that \(\lbrace f_{ij}^{\text{new}}\rbrace\) and \(\lbrace g_{ij}^{\text{new}}\rbrace\) satisfy (\ref{eq84un8byyyyybw7tttttttttyb7}). Thus, by rewriting system (\ref{eq0w9b38b47v8hbub8r57qw3t6n}) with \(\lbrace f_{ij}^{\text{new}}\rbrace\) and \(\lbrace g_{ij}^{\text{new}}\rbrace\), this system takes the form (\ref{eq19000}) such that (\ref{eq84un8byyyyybw7tttttttttyb7}) holds.

Hereafter, we assume that (\ref{eq84un8byyyyybw7tttttttttyb7}) is satisfied for system (\ref{eq19000}). Write this system as
\begin{equation}\label{eq32000}
\begin{aligned}
\dot{u}_{1} &= -\lambda_{1}u_{1} + f_{11}(u_{1}, v) u_{1} + f_{12}(u_{1}, u_{2}, v) u_{2},\\
\dot{u}_{2} &= -\lambda_{2}u_{2} + f_{21}(u_{1}, v)u_{1} + J_{1}(u, v) u_{2} + \underline{J_{2}(v)u_{2}},\\
\dot{v}_{1} &= +\lambda_{1}v_{1} + g_{11}(u, v_{1}) v_{1} + g_{12}(u, v_{1}, v_{2}) v_{2},\\
\dot{v}_{2} &= +\lambda_{2}v_{2} + g_{21}(u, v_{1})v_{1} + J_{3}(u,v) v_{2} + \underline{J_{4}(u) v_{2}},
\end{aligned}
\end{equation}
where
\begin{equation*}
\begin{gathered}
J_{1}(u,v) = f_{22}(u, v)-f_{22}(0, v),\qquad J_{2}(v) = f_{22}(0, v),\\
J_{3}(u,v) = g_{22}(u, v)-g_{22}(u, 0),\qquad J_{4}(u) = g_{22}(u, 0).
\end{gathered}
\end{equation*}
In order to obtain conditions (\ref{eq28000}) and (\ref{eq29000}), we need to find a change of coordinates which eliminates the underlined terms in (\ref{eq32000}). We claim that this is possible by applying two consecutive \(\mathcal{C}^{\infty}\)-smooth changes of coordinates of the forms
\begin{equation}\label{eq30000}
\tilde{u}_{1} = u_{1},\qquad \tilde{u}_{2} = u_{2}+q_{1}(v_{1}, v_{2})u_{2},\qquad
\tilde{v}_{1} = v_{1},\qquad \tilde{v}_{2} = v_{2},
\end{equation}
and
\begin{equation}\label{eq31000}
\tilde{u}_{1} = u_{1},\qquad \tilde{u}_{2} = u_{2},\qquad
\tilde{v}_{1} = v_{1},\qquad \tilde{v}_{2} = v_{2}+q_{2}(u_{1}, u_{2})v_{2},
\end{equation}
where \(q_{1}\) and \(q_{2}\) are some functions such that \(q_{1}\left(0\right) = q_{2}\left(0\right) = 0\). We show that the underlined terms \(J_{2}(v) u_{2}\) and \(J_{4}(u) v_{2}\) can be eliminated by applying a change of coordinates of the forms (\ref{eq30000}) and (\ref{eq31000}), respectively.

Let \(z = \frac{u_{2}}{1+q_{1}(v)}\). Applying change of coordinates (\ref{eq30000}) brings system (\ref{eq32000}) to
\begin{equation}\label{eq34000}
\begin{aligned}
\dot{u}_{1} &= -\lambda_{1}u_{1} + \left[f_{11}(u_{1}, v)\right] u_{1} + \left[\frac{f_{12}\left(u_{1}, z, v\right)}{1+q_{1}(v)}\right] u_{2},\\
\dot{u}_{2} &= -\lambda_{2}u_{2} + \left[\left(1+q_{1}(v)\right) f_{21}\left(u_{1}, v\right)\right] u_{1} + Q_{1}(u,v)u_{2} + Q_{2}(v)u_{2}\\
\dot{v}_{1} &= \lambda_{1}v_{1} + \left[g_{11}\left(u_{1}, z, v_{1}\right)\right] v_{1} + \left[g_{12}\left(u_{1}, z, v_{1}, v_{2}\right)\right] v_{2},\\
\dot{v}_{2} &= \lambda_{2}v_{2} + \left[g_{21}\left(u_{1}, z, v_{1}\right)\right] v_{1} + \left[g_{22}\left(u_{1}, z, v\right)\right] v_{2},
\end{aligned}
\end{equation}
where
\begin{equation*}
\begin{aligned}
Q_{1}(u,v) =& J_{1}\left(u_{1}, z, v\right) + \frac{{q_{1}}_{v_{1}}\left(v\right)}{1+q_{1}(v)}\cdot\Big[g_{11}\left(u_{1}, z, v_{1}\right) v_{1} + g_{12}\left(u_{1}, z, v\right) v_{2}\Big]\\
+& \frac{{q_{1}}_{v_{2}}\left(v\right)}{1+q_{1}(v)}\cdot\Big[g_{21}\left(u_{1}, z, v_{1}\right)v_{1} + g_{22}\left(u_{1}, z, v\right)v_{2} - g_{21}\left(0, v_{1}\right)v_{1} -g_{22}\left(0, v\right)v_{2}\Big],\\
Q_{2}(v) =& J_{2}(v)+\frac{\lambda_{1}{q_{1}}_{v_{1}}(v) v_{1} + {q_{1}}_{v_{2}}(v)\big(\lambda_{2}v_{2} + g_{21}(0, v_{1})v_{1} + g_{22}(0,v) v_{2}\big)}{1+q_{1}(v)}.
\end{aligned}
\end{equation*}
It is easy to see that \(Q_{1}\) vanishes at \(u = 0\) and also the updated \(f_{ij}\) and \(g_{ij}\) in system (\ref{eq34000}) satisfy all the conditions (\ref{eq20000}) except for (\ref{eq28000}) and (\ref{eq29000}). In order to get (\ref{eq28000}), it is sufficient to find \(q_{1}(v)\) such that \(Q_{2}(v)\equiv 0\), i.e. \(q_{1}(v)\) satisfies the relation
\begin{equation}\label{eq37000}
-\left(1+q_{1}(v)\right)J_{2}(v)= {q_{1}}_{v_{1}}(v)\cdot\big[\lambda_{1} v_{1}\big] + {q_{1}}_{v_{2}}(v)\cdot\big[\lambda_{2}v_{2} + g_{21}(0, v_{1})v_{1} + g_{22}(0,v) v_{2}\big].
\end{equation}
Consider the \(\mathcal{C}^{\infty}\)-smooth system
\begin{equation}\label{eq38000}
\begin{aligned}
\dot{U} &= -\left(1+U\right)J_{2}(v),\\
\dot{v}_{1} &= \lambda_{1} v_{1},\\
\dot{v}_{2} &= \lambda_{2}v_{2} + g_{21}(0, v_{1})v_{1} + g_{22}(0,v) v_{2},
\end{aligned}
\end{equation}
where \(\left(U, v_{1}, v_{2}\right)\in \mathbb{R}^{3}\). The linear part of this system at the origin is
\begin{equation*}
\left(\begin{array}{ccc}
0 & \frac{\partial J_{2}}{\partial v_{1}}(0) & \frac{\partial J_{2}}{\partial v_{2}}(0)\\
0 & \lambda_{1} & 0 \\
0 & 0 & \lambda_{2}
\end{array}\right),
\end{equation*}
with the spectrum \(\lbrace 0, \lambda_{1}, \lambda_{2}\rbrace\). Therefore, this system has a \(\mathcal{C}^{\infty}\)-smooth 2-dimensional local unstable invariant manifold defined by the equation \(\lbrace U = q_{1}\left(v_{1}, v_{2}\right)\rbrace\) for some function \(q_{1}\). Moreover, this function satisfies (\ref{eq37000}) because this relation is nothing but the condition of the invariance of the local unstable invariant manifold with respect to the flow of system (\ref{eq38000}) (see Definition \ref{invarianceconditiondefinition}). Thereby, as we required, a \(\mathcal{C}^{\infty}\)-smooth function \(q_{1}(v_{1}, v_{2})\) that fulfills (\ref{eq37000}) exists.

We have now shown our first claim: the underlined term \(J_{2}(v) u_{2}\) in system (\ref{eq32000}) can be eliminated by performing a change of coordinates of the form (\ref{eq30000}). The proof for the second claim that the term \(J_{4}(u) v_{2}\) can be eliminated by applying a change of coordinates of the form (\ref{eq31000}) can be accomplished analogously (see \cite{BakraniPhDthesis}).

Note that, by Remark \ref{rem8a7vcyut4l7acvtk3iyulvyt7w}, since the desired change of coordinates in Lemma \ref{Dimabooknormalform} commutes with symmetry (\ref{eq425}), we have that system (\ref{eq32000}) is invariant with respect to this symmetry. Therefore, system (\ref{eq38000}) is invariant with respect to the symmetries \(v_{1} \leftrightarrow -v_{1}\). Thus, \(q_{1}\left(v_{1}, v_{2}\right) = q_{1}\left(-v_{1}, v_{2}\right)\). Analogously, \(q_{2}\left(u_{1}, u_{2}\right) = q_{2}\left(-u_{1}, u_{2}\right)\). This means that changes of coordinates (\ref{eq30000}) and (\ref{eq31000}) commute with symmetry (\ref{eq425}) too.

To finish the proof, note that changes of coordinates (\ref{eq30000}) and (\ref{eq31000}) are of the form (\ref{eq137010}). Therefore, by Remark \ref{rem8a7vcyut4l7acvtk3iyulvyt7w}, applying changes of coordinates (\ref{eq30000}) and (\ref{eq31000}) together with the change of coordinates used in Lemma \ref{Dimabooknormalform} transforms first integral (\ref{eq98n9uo8wbr8bau4b819ub92q1u8}) to the form (\ref{eq22000}). This ends the proof of Lemma \ref{Nftheorem}.
\end{proof}

\begin{proof}[Proof of Lemma \ref{Nftheorem2}]
By Lemma \ref{Nftheorem}, there exists a change of coordinates which commutes with symmetry (\ref{eq425}) and brings system (\ref{eq300}) and first integral (\ref{eq400}) to (\ref{eq19000}) and (\ref{eq22000}), respectively. System (\ref{eq19000}) possesses a \(\mathcal{C}^{q}\)-smooth three dimensional extended unstable invariant manifold \(W^{uE}\) defined by \(\lbrace \left(u,v\right): \, u_{2} = \phi^{uE}\left(u_{1}, v_{1}, v_{2}\right)\rbrace\), and a \(\mathcal{C}^{q}\)-smooth extended stable invariant manifold \(W^{sE}\) defined by \(\lbrace \left(u,v\right): \, v_{2} = \phi^{sE}\left(u_{1}, u_{2}, v_{1}\right)\rbrace\) (see \cite{Dimabook}). We claim that straightening \(W^{uE}\), i.e. applying the \(\mathcal{C}^{q}\)-smooth change of coordinates
\begin{equation}\label{eq43000}
\tilde{u}_{1} = u_{1},\qquad \tilde{u}_{2} = u_{2}-\phi^{uE}\left(u_{1}, v_{1}, v_{2}\right),\qquad
\tilde{v}_{1} = v_{1},\qquad \tilde{v}_{2} = v_{2},
\end{equation}
and straightening \(W^{sE}\), i.e. applying the \(\mathcal{C}^{q}\)-smooth change of coordinates
\begin{equation}\label{eq44000}
\tilde{u}_{1} = u_{1},\qquad \tilde{u}_{2} = u_{2},\qquad
\tilde{v}_{1} = v_{1},\qquad\tilde{v}_{2} = v_{2}-\phi^{sE}\left(u_{1}, u_{2}, v_{1}\right),
\end{equation}
reduce system (\ref{eq19000}) to system (\ref{eq23000}), where (\ref{eq20000}) is satisfied, and transforms first integral (\ref{eq22000}) to (\ref{eqhhuhut5989b8b8y4qlz37y7t}). On the other hand, by Remark \ref{Cor99090q19n9be8u8y3f}, straightening these manifolds keeps the invariance of system (\ref{eq19000}) and first integral (\ref{eq22000}) with respect to symmetry (\ref{eq425}). Thus, we are done once we prove this claim. To this end, we use the following lemma
\begin{mylem}\label{lem10000}
The following hold for the \(\mathcal{C}^{q}\)-smooth functions \(\phi^{uE}\) and \(\phi^{sE}\):\\
\begin{inparaenum}[(i)]
\item\label{Item87bibi7y7yi} \(\phi^{uE}\left(0, v_{1}, v_{2}\right) \equiv \phi^{uE}_{u_{1}}\left(0, v_{1}, v_{2}\right)\equiv 0\), \hspace{10mm}
\item\label{Item87bibi7y7yii} \(\phi^{sE}\left(u_{1}, u_{2}, 0\right) \equiv \phi^{sE}_{v_{1}}\left(u_{1}, u_{2}, 0\right)\equiv 0\).
\end{inparaenum}
\end{mylem}
The following are immediate consequences of this lemma:
\begin{mycor}\label{Cor1100}
We can write \(\phi^{uE}\) and \(\phi^{sE}\) as
\begin{equation*}
\begin{aligned}
\phi^{uE}\left(u_{1}, v_{1}, v_{2}\right) &= u_{1} p_{1}^{uE}(u_{1}, v_{1}, v_{2}) = u_{1}^{2} p_{2}^{uE}(u_{1}, v_{1}, v_{2}),\\
\phi^{sE}\left(u_{1}, u_{2}, v_{1}\right) &= v_{1} p_{1}^{sE}(u_{1}, u_{2}, v_{1}) = v_{1}^{2} p_{2}^{sE}(u_{1}, u_{2}, v_{1}),
\end{aligned}
\end{equation*}
where \(p_{1}^{sE}\) and \(p_{1}^{uE}\) are some \(\mathcal{C}^{q-1}\)-smooth functions and \(p_{2}^{sE}\) and \(p_{2}^{uE}\) are some \(\mathcal{C}^{q-2}\)-smooth functions such that \(p_{1}^{uE} = u_{1}p_{2}^{uE}\) and \(p_{1}^{sE} = v_{1} p_{2}^{sE}\).
\end{mycor}
\begin{mycor}\label{Cor001ibfuhg09ut2g}
We have\\
\begin{inparaenum}[(i)]
\item \(\phi^{uE}_{v_{1}}\left(0, v_{1}, v_{2}\right) \equiv \phi^{uE}_{v_{2}}\left(0, v_{1}, v_{2}\right)\equiv 0\), \hspace{10mm}
\item \(\phi^{sE}_{u_{1}}\left(u_{1}, u_{2}, 0\right) \equiv \phi^{sE}_{u_{2}}\left(u_{1}, u_{2}, 0\right)\equiv 0\).
\end{inparaenum}
\end{mycor}

We prove Lemma \ref{lem10000} later. Taking into account that as a result of applying change of coordinates (\ref{eq43000}), the equation of \(\dot{u}_{2}\) vanishes at \(\{u_{2} = 0\}\), one can easily see that change of coordinates (\ref{eq43000}) reduces system (\ref{eq19000}) to
\begin{equation}\label{eqwsadvfm98nub8ublkaigr1a10}
\begin{aligned}
\dot{u}_{1} &= -\lambda_{1}u_{1} + \tilde{f}_{11}\left(u_{1}, v\right) u_{1} + \tilde{f}_{12}\left(u_{1}, u_{2}, v\right) u_{2},\\
\dot{u}_{2} &= -\lambda_{2}u_{2} + \tilde{f}_{22}\left(u_{1}, u_{2}, v\right) u_{2},\\
\dot{v}_{1} &= +\lambda_{1}v_{1} + \tilde{g}_{11}\left(u, v_{1}\right) v_{1} + \tilde{g}_{12}\left(u, v_{1}, v_{2}\right) v_{2},\\
\dot{v}_{2} &= +\lambda_{2}v_{2} + \tilde{g}_{21}\left(u, v_{1}\right)v_{1} +  \tilde{g}_{22}\left(u, v_{1}, v_{2}\right) v_{2},
\end{aligned}
\end{equation}
where
\begin{equation*}
\begin{aligned}
\tilde{f}_{11}\left(u_{1}, v\right) &= f_{11}\left(u_{1}, v\right) + f_{12}\left(u_{1}, x, v\right)\, p_{1}^{uE}\left(u_{1}, v\right),\\
\tilde{f}_{12}\left(u, v\right) &= f_{12}\left(u_{1}, u_{2} + x, v\right) + P_{1}\left(u, v\right)x,\\
\tilde{f}_{22}\left(u, v\right) &= f_{22}\left(u_{1}, u_{2}+x, v\right) + P_{2}\left(u,v\right)x -\phi_{u_{1}}^{uE}\left(u_{1}, v\right)\left[f_{12}\left(u_{1}, u_{2} + x, v\right) + P_{3}\left(u, v\right)x\right]\\
-&\phi_{v_{1}}^{uE} \left(u_{1}, v\right)\left[P_{4}\left(u, v\right) v_{1} + P_{5}\left(u, v\right) v_{2}\right] - \phi_{v_{2}}^{uE}\left(u_{1}, v\right)\left[P_{6}\left(u, v\right) v_{1} + P_{7}\left(u, v\right) v_{2}\right],\\
\tilde{g}_{i1}\left(u, v_{1}\right) &= g_{i1}\left(u_{1}, u_{2} + x, v_{1}\right), \qquad\quad
\tilde{g}_{i2}\left(u, v\right) = g_{i2}\left(u_{1}, u_{2} + x, v\right), \qquad (i = 1,2).
\end{aligned}
\end{equation*}
Here, \(x:= \phi^{uE}\left(u_{1}, v\right)\), the functions \(f_{ij}\) and \(g_{ij}\) are as in (\ref{eq19000}), and \(P_{j}\left(u, v\right)\) (\(j=1,..,7\)) are some functions such that \(P_{1}\) vanishes at \(\{v = 0\}\) (see \cite{BakraniPhDthesis} for more details). Moreover, \(\tilde{f}_{ij}\) are \(\mathcal{C}^{q-1}\)-smooth and \(\tilde{g}_{ij}\) are \(\mathcal{C}^{q}\)-smooth. Using Lemma \ref{lem10000} and Corollaries \ref{Cor1100} and \ref{Cor001ibfuhg09ut2g} and taking into account that the expression \(u_{2} + \phi^{uE}\left(u_{1}, v\right)\) vanish at \(u=0\), and also the functions \(f_{ij}\) and \(g_{ij}\) satisfy (\ref{eq20000}), one can easily show that \(\tilde{f}_{ij}\) and \(\tilde{g}_{ij}\) satisfy (\ref{eq20000}) as well.

System (\ref{eqwsadvfm98nub8ublkaigr1a10}) is of the form (\ref{eq19000}) where   \(f_{21}\left(u_{1}, v\right)\equiv 0\). Similar to the case of straightening the extended unstable manifold, one can use Lemma \ref{lem10000} and Corollaries \ref{Cor1100} and \ref{Cor001ibfuhg09ut2g} and show that making change of coordinates (\ref{eq44000}) reduces system (\ref{eqwsadvfm98nub8ublkaigr1a10}) to system (\ref{eq23000}) where the corresponding \(f_{ij}\) and \(g_{ij}\) are \(\mathcal{C}^{q-1}\)-smooth and satisfy (\ref{eq20000}). This ends the proof of the first part of Lemma \ref{Nftheorem2}.

Denote the \(H_{1}\) and \(H_{2}\) in (\ref{eq22000}) by \(H_{1}^{\circ}\) and \(H_{2}^{\circ}\), respectively, and let \(x := \left(u_{1}, u_{2}, v_{1}\right)\), \(y := \left(u_{1}, v_{1}, v_{2} + \phi^{sE}\left(x\right)\right)\) and \(z:= \left(u_{1}, u_{2} + \phi^{uE}\left(y\right), v_{1}, v_{2} + \phi^{sE}\left(x\right)\right)\).
Applying changes of coordinates (\ref{eq43000}) and (\ref{eq44000}) brings (\ref{eq22000}) to
\begin{equation*}
\begin{aligned}
H =& \lambda_{1}u_{1}v_{1} \left[1 + H_{1}^{\circ}\left(z\right)\right] - \lambda_{2} \left(u_{2} + \phi^{uE}\left(y\right)\right) \left(v_{2} +\phi^{sE}\left(x\right)\right) \left[1 + H^{\circ}_{2}\left(z\right)\right],
\end{aligned}
\end{equation*}
which by Corollary \ref{Cor1100}, can be written in the form (\ref{eqhhuhut5989b8b8y4qlz37y7t}), for
\begin{equation*}
\begin{aligned}
& H_{1} = H^{\circ}_{1}\left(z\right) + \lambda_{2}\lambda_{1}^{-1} p_{1}^{uE}\left(y\right)p_{1}^{sE}\left(x\right) \left[1 + H_{2}^{\circ}\left(z\right)\right],\quad
&& H_{2} = H_{2}^{\circ},\\
& H_{3} = p_{2}^{sE}\left(x\right) \left[1 + H^{\circ}_{2}\left(z\right)\right],\quad
&& H_{4} = p_{2}^{uE}\left(y\right) \left[1 + H^{\circ}_{2}\left(z\right)\right].
\end{aligned}
\end{equation*}
This proves the second part of Lemma \ref{Nftheorem2}.

All that remains to finish the proof of Lemma \ref{Nftheorem2} is proving Lemma \ref{lem10000}. We only prove part (\ref{Item87bibi7y7yi}) of this lemma; the proof of part (\ref{Item87bibi7y7yii}) is the same.

The first identity 
\begin{equation}\label{eqj9srugsqyuplazmjc0s}
\phi^{uE}\left(0, v_{1}, v_{2}\right) \equiv 0,
\end{equation}
is an immediate consequence of the fact that the extended unstable invariant manifold \(W^{uE}\) contains the unstable invariant manifold \(\lbrace u_{1}=u_{2}=0\rbrace\) (see Section \cite{Dimabook}). Indeed, for any \(\left(v_{1},v_{2}\right)\), we have that \(\left(0, 0, v_{1}, v_{2}\right)\) belongs to \(\{\left(u,v\right): u_{2} = \phi^{uE}\left(u_{1}, v_{1}, v_{2}\right)\}\). This implies \(\phi^{uE}\left(0, v_{1}, v_{2}\right) = 0\), for any \(\left(v_{1},v_{2}\right)\), which yields (\ref{eqj9srugsqyuplazmjc0s}).

It is important to notice that relation (\ref{eqj9srugsqyuplazmjc0s}) is sufficient to obtain the statement of part I of Corollary \ref{Cor001ibfuhg09ut2g}. In other words, (\ref{eqj9srugsqyuplazmjc0s}) implies part I of Corollary \ref{Cor001ibfuhg09ut2g}.

To prove the identity
\begin{equation}\label{eq44500}
\phi^{uE}_{u_{1}}\left(0, v_{1}, v_{2}\right) \equiv 0,
\end{equation}
we consider the condition of the invariance of the manifold \(W^{uE}\) with respect to the flow of system (\ref{eq19000}) (see Definition \ref{invarianceconditiondefinition}), i.e.
\begin{equation*}
\begin{aligned}
-\lambda_{2}x + f_{21}&\left(u_{1}, v\right) u_{1} + f_{22}\left(u_{1}, x, v\right) x = \phi^{uE}_{u_{1}}\left(u_{1}, v\right)\big[-\lambda_{1}u_{1} + f_{11}\left(u_{1}, v\right) u_{1}\\
&+ f_{12}\left(u_{1}, x, v\right)x\big] +\phi^{uE}_{v_{1}}\left(u_{1}, v\right) \big[\lambda_{1}v_{1} + g_{11}\left(u_{1}, x, v_{1}\right) v_{1} + g_{12}\left(u_{1}, x, v\right) v_{2}\big]\\
& +\phi^{uE}_{v_{2}}\left(u_{1}, v\right) \big[\lambda_{2}v_{2} + g_{21}\left(u_{1}, x, v_{1}\right)v_{1} +  g_{22}\left(u_{1}, x, v\right) v_{2}\big],
\end{aligned}
\end{equation*}
where \(x:= \phi^{uE}\left(u_{1}, v\right)\). Both sides of this relation are \(\mathcal{C}^{q-1}\)-smooth (\(q\geq 2\) because \(2\lambda_{1} < \lambda_{2}\)) functions of \(u_{1}\), \(v_{1}\) and \(v_{2}\). Taking (\ref{eqj9srugsqyuplazmjc0s}) as well as conditions (\ref{eq20000}) and Corollary \ref{Cor001ibfuhg09ut2g} into account, we can differentiate this relation with respect to \(u_{1}\) at \(u_{1} = 0\) and obtain
\begin{equation}\label{eq45500}
\begin{aligned}
0 =& \left[\left(\lambda_{2}-\lambda_{1}\right)\phi^{uE}_{u_{1}}\left(0, v\right) + f_{12}(0,v)\left(\phi^{uE}_{u_{1}}\left(0, v\right)\right)^{2}\right] + \big[\lambda_{1}v_{1}\big]\phi^{uE}_{u_{1}v_{1}}\left(0, v\right)\\
&+ \big[\lambda_{2}v_{2} + g_{21}(0,v_{1})v_{1} + g_{22}(0,v)v_{2}\big]\phi^{uE}_{u_{1}v_{2}}\left(0, v\right).
\end{aligned}
\end{equation}
Define \(z = z(v) = \phi^{uE}_{u_{1}}(0,v)\). Then, (\ref{eq45500}) can be written as
\begin{equation}\label{eq52000}
\begin{aligned}
0 = \left[\left(\lambda_{2}-\lambda_{1}\right)z + f_{12}(0,v)z^{2}\right] + \big[\lambda_{1}v_{1}\big]\cdot\frac{\partial z(v)}{\partial v_{1}}
+ \big[\lambda_{2}v_{2} + g_{21}(0,v_{1})v_{1} + g_{22}(0,v)v_{2}\big]\cdot\frac{\partial z(v)}{\partial v_{2}},
\end{aligned}
\end{equation}
where \(z(0)=0\) (note that \(\phi^{uE}_{u_{1}}(0,0,0) = 0\)).

To get (\ref{eq44500}), we need to show \(z\left(v\right) \equiv 0\). First, note that \(z\left(v\right) \equiv 0\) satisfies (\ref{eq52000}). Thus, (\ref{eq44500}) holds if we show that \(z\equiv 0\) is the unique solution of (\ref{eq52000}). Note that, by Proposition \ref{invariancecondition}, \(z\left(v\right)\) satisfies (\ref{eq52000}) if and only if the 2-dimensional manifold
\begin{equation}\label{eq54000}
\{\left(v, z\right): z = z\left(v\right)\,\,\text{and}\,\, z\left(0\right) = 0 \}
\end{equation}
be invariant with respect to the flow of the \(\mathcal{C}^{q-1}\)-smooth system
\begin{equation}\label{eq53000}
\begin{aligned}
\dot{v}_{1} &= -\lambda_{1}v_{1},\\
\dot{v}_{2} &= -\lambda_{2}v_{2} - g_{21}(0,v_{1})v_{1} - g_{22}(0,v)v_{2},\\
\dot{z} &= \left(\lambda_{2}-\lambda_{1}\right)z + f_{12}(0,v)z^{2},
\end{aligned}
\end{equation}
which is defined on a small neighborhood of the origin in \(\mathbb{R}^{3}\). (Indeed, relation (\ref{eq52000}) is the condition of the invariance of (\ref{eq54000}) with respect to the flow of system (\ref{eq53000}).) Therefore, the uniqueness of the solution of (\ref{eq52000}) can be proved by showing that system (\ref{eq53000}) has a unique invariant manifold of the form (\ref{eq54000}). To do this, first, notice that this system possesses a unique two dimensional stable invariant manifold of form (\ref{eq54000}). Second, we observe that any orbit on manifold (\ref{eq54000}) converges to the origin of system (\ref{eq53000}): the first two equations in (\ref{eq53000}) are independent of \(z\) and have \(\left(v_{1}, v_{2}\right) = \left(0, 0\right)\) as an asymptotically stable equilibrium. Therefore, as \(t\rightarrow \infty\), an orbit \(\left(v(t), z(v(t))\right)\) of system (\ref{eq53000}) which belongs to invariant manifold (\ref{eq54000}) converges to \(\left(0, z\left(0\right)\right)\). Since \(z\left(0\right)=0\), this means that any invariant manifold of the form (\ref{eq53000}) must be a subset of the stable manifold of system (\ref{eq53000}). However, since both manifolds are 2-dimensional, they must be the same. Therefore, system (\ref{eq53000}) has a unique invariant manifold of the form (\ref{eq54000}) which is in fact its stable invariant manifold. This ends the proof of Lemma \ref{lem10000} and hence the proof of Lemma \ref{Nftheorem2}.
\end{proof}

\section{Proofs of Lemmas \ref{flowlemmaresonant}, \ref{flowlemma} and \ref{flowlemma0}}\label{ProofsOfBoundaryValueProblems}

We only prove Lemma \ref{flowlemmaresonant}. The proofs of Lemmas \ref{flowlemma} and \ref{flowlemma0} are similar. We refer the reader to \cite{BakraniPhDthesis} for the proofs of these lemmas. We start with a discussion on the method of the proof and then proceed to the proof of Lemma \ref{flowlemmaresonant}.

\subsection{The method of the proof}

Here, we present the main procedure which is used in the proofs of Lemmas \ref{flowlemmaresonant}, \ref{flowlemma}, \ref{flowlemma0} and also Lemma \ref{Derivativeflowlemma} (see Appendix \ref{Appendix989897g86g7656u}). Consider system (\ref{eqbwiebro4ibyw43yyyiw0290}) and denote its unique solution that satisfies boundary condition (\ref{eq800020}) by  
\(\left(u^{*}, v^{*}\right)\), where \(u^{*} = (u_{1}^{*},u_{2}^{*})\) and \(v^{*} = (v_{1}^{*},v_{2}^{*})\). We may also write this as
\begin{equation}\label{eq800050}
\left(u^{*}, v^{*}\right) = \left(u^{*}\left(t\right), v^{*}\left(t\right)\right) = \big(u^{*}\left(t, \tau, u_{10}, u_{20}, v_{1\tau}, v_{2\tau}\right), v^{*}\left(t, \tau, u_{10}, u_{20}, v_{1\tau}, v_{2\tau}\right)\big),
\end{equation}
to emphasise that in addition to time variable \(t\), this solution explicitly depends on \(\tau\) and the boundary conditions \(u_{10}\), \(u_{20}\), \(v_{1\tau}\) and \(v_{2\tau}\) as well. It is easy to see that \(\left(u^{*}(t), v^{*}(t)\right)\) is a solution of this system with boundary conditions (\ref{eq800020}) if and only if
\begin{equation}\label{eq800025}
\begin{aligned}
u_{1}^{*}(t) =& e^{-\lambda_{1}t} u_{10} + \int_{0}^{t} e^{\lambda_{1}(s-t)} F_{1}\left(u^{*}\left(s\right), v^{*}\left(s\right)\right) ds,\\
u_{2}^{*}(t) =& e^{-\lambda_{2}t} u_{20} + \int_{0}^{t} e^{\lambda_{2}(s-t)} F_{2}\left(u^{*}\left(s\right), v^{*}\left(s\right)\right) ds,\\
v_{1}^{*}(t) =& e^{-\lambda_{1}\left(\tau-t\right)} v_{1\tau} - \int_{t}^{\tau} e^{-\lambda_{1}(s-t)} G_{1}\left(u^{*}\left(s\right), v^{*}\left(s\right)\right) ds,\\
v_{2}^{*}(t) =& e^{-\lambda_{2}\left(\tau-t\right)} v_{2\tau} - \int_{t}^{\tau} e^{-\lambda_{2}(s-t)} G_{2}\left(u^{*}\left(s\right), v^{*}\left(s\right)\right) ds.
\end{aligned}
\end{equation}
For a given \(\tau\), denote by \(\mathcal{I}\) the set of all vector valued functions \(\left(u_{1}\left(t\right), u_{2}\left(t\right), v_{1}\left(t\right), v_{2}\left(t\right)\right)\) defined for \(t\in [0, \tau]\) on some small neighborhood of the origin in \(\mathbb{R}^{4}\). Then, the right-hand side of (\ref{eq800025}) defines an integral operator on \(\mathcal{I}\), denote it by \(\mathfrak{T}\), as follows:
\begin{equation*}
\mathfrak{T}: \left(u_{1}\left(t\right), u_{2}\left(t\right), v_{1}\left(t\right), v_{2}\left(t\right)\right) \mapsto \left(\overline{u}_{1}\left(t\right), \overline{u}_{2}\left(t\right), \overline{v}_{1}\left(t\right), \overline{v}_{2}\left(t\right)\right),
\end{equation*}
where
\begin{equation*}
\begin{aligned}
\overline{u}_{1}\left(t\right) &= e^{-\lambda_{1}t} u_{10} + \int_{0}^{t} e^{\lambda_{1}(s-t)} F_{1}\left(u\left(s\right), v\left(s\right)\right) ds,\\
\overline{u}_{2}\left(t\right) &= e^{-\lambda_{2}t} u_{20} + \int_{0}^{t} e^{\lambda_{2}(s-t)} F_{2}\left(u\left(s\right), v\left(s\right)\right) ds,\\
\overline{v}_{1}\left(t\right) &= e^{-\lambda_{1}\left(\tau -t\right)} v_{1\tau} - \int_{t}^{\tau} e^{-\lambda_{1}(s-t)} G_{1}\left(u\left(s\right), v\left(s\right)\right) ds,\\
\overline{v}_{2}\left(t\right) &= e^{-\lambda_{2}\left(\tau -t\right)} v_{2\tau} - \int_{t}^{\tau} e^{-\lambda_{2}(s-t)} G_{2}\left(u\left(s\right), v\left(s\right)\right) ds.
\end{aligned}
\end{equation*}

The solution \(\left(u^{*}(t), v^{*}(t)\right)\) is in fact the fixed point of this integral operator. According to \cite{Dimabook} (Theorems 2.9 and 5.11), this integral operator is a contraction and its fixed point is the limit of the sequence of successive approximations
\begin{equation*}
\Big\{ \left(u^{\left(n\right)}(t), v^{\left(n\right)}(t)\right) = \left(u_{1}^{\left(n\right)}(t), u_{2}^{\left(n\right)}(t), v_{1}^{\left(n\right)}(t), v_{2}^{\left(n\right)}(t)\right)\Big\}_{n=0}^{n=\infty},
\end{equation*}
where \(\left(u^{\left(0\right)}, v^{\left(0\right)}\right) \equiv \left(0,0\right)\) and
\begin{equation*}
\left(u^{\left(n+1\right)}(t), v^{\left(n+1\right)}(t)\right) = \mathfrak{T} \left(u^{\left(n\right)}\left(t\right), v^{\left(n\right)}\left(t\right)\right),\quad \forall n\geq 0.
\end{equation*}

Let \(\mathcal{A}\) be a closed subset of \(\mathcal{I}\) such that \(\left(u\left(t\right), v\left(t\right)\right) \equiv \left(0, 0\right) \in \mathcal{A}\) and \(\mathfrak{T}\left(\mathcal{A}\right) \subset \mathcal{A}\). Since \(\left(u^{\left(0\right)}, v^{\left(0\right)}\right) \equiv \left(0,0\right)\in \mathcal{A}\), the invariance of \(\mathcal{A}\) implies that \(\left(u^{\left(n\right)}(t), v^{\left(n\right)}(t)\right)\) belongs to \(\mathcal{A}\) for all \(n> 0\), and so does the solution  \(\left(u^{*}(t), v^{*}(t)\right)\).
\begin{myrem}\label{Rem9n9iqnubxfybxbxkiqub}
Assume that there exists a {\it 'certain estimate'} which for any arbitrary \(\left(u\left(t\right), v\left(t\right)\right)\in \mathcal{A}\), its image \(\mathfrak{T}\left(u\left(t\right), v\left(t\right)\right)\) satisfies. Therefore, since \(\mathfrak{T}\left(u^{*}, v^{*}\right) = \left(u^{*}, v^{*}\right)\), the solution \(\left(u^{*}, v^{*}\right)\) itself satisfies that certain estimate as well.
\end{myrem}

Our approach for proving Lemmas \ref{flowlemmaresonant}, \ref{flowlemma} and \ref{flowlemma0} (and Lemma \ref{Derivativeflowlemma}) is based on this remark. We construct the integral operator, introduce the invariant set \(\mathcal{A}\) and find an estimate for the image of the elements of this set under \(\mathfrak{T}\). Then, this estimate holds for the solution \(\left(u^{*}, v^{*}\right)\) too.

\subsection{Proof of Lemma \ref{flowlemmaresonant}}

Throughout, we use the following notation: \(x=\left(u, v\right)\), \(x\left(t\right) = \left(u\left(t\right), v\left(t\right)\right)\), and \(x\left(s\right) = \left(u\left(s\right), v\left(s\right)\right)\).
Recast system (\ref{eq74wt366cuw6gt5uw6v01092}) into the form (\ref{eqbwiebro4ibyw43yyyiw0290}), where
\begin{equation}\label{eqgyyty736q487eg6347ctsg}
\begin{aligned}
F_{i}\left(x\right) =& \mathtt{f}_{i1}\left(x\right) u_{1}^{2} + \mathtt{f}_{i2}\left(x\right) u_{1} u_{2} + \mathtt{f}_{i3}\left(x\right) u_{1} v_{1} + \mathtt{f}_{i4}\left(x\right) u_{1} v_{2} + \mathtt{f}_{i5}\left(x\right) u_{2}^{2}\\ &+ \mathtt{f}_{i6}\left(x\right) u_{2} v_{1} + \mathtt{f}_{i7}\left(x\right) u_{2} v_{2},\\
G_{i}\left(x\right) =& \mathtt{g}_{i1}\left(x\right) v_{1}^{2} + \mathtt{g}_{i2}\left(x\right) v_{1} v_{2} + \mathtt{g}_{i3}\left(x\right) v_{1} u_{1} + \mathtt{g}_{i4}\left(x\right) v_{1} u_{2} + \mathtt{g}_{i5}\left(x\right) v_{2}^{2}\\ &+ \mathtt{g}_{i6}\left(x\right) v_{2} u_{1} + \mathtt{g}_{i7}\left(x\right) v_{2} u_{2},
\end{aligned}
\end{equation}
for \(i=1,2\), and some continuous functions \(\mathtt{f}_{ij}\) and \(\mathtt{g}_{ij}\). Let \(\Omega\) be a small compact neighborhood of \(O\) and define
\begin{equation}\label{eqiw76iebyac7t4jv6}
M^{*} := \sup_{\left(u, v\right)\in \Omega} \big\{ \lvert \mathtt{f}_{ij}\left(u, v\right)\rvert, \lvert \mathtt{g}_{ij}\left(u, v\right)\rvert\big\}.
\end{equation}
Let \(\delta > 0\) be small and consider the set
\begin{equation}\label{eqmjhscfnhf6ba7i6w38nrgfjsyg}
\mathcal{A} = \big\{x\left(t\right): \quad \lvert u_{1}(t)\rvert, \lvert u_{2}(t)\rvert \leq 2e^{-\lambda t}\delta,\quad \lvert v_{1}(t)\rvert, \lvert v_{2}(t)\rvert \leq 2e^{-\lambda\left(\tau-t\right)}\delta\big\},
\end{equation}
where \(x\left(t\right)\) is any continuous function defined on \(\Omega\) for \(t \in \left[0, \tau\right]\).

We first show that \(\mathcal{A}\) is invariant with respect to the integral operator \(\mathfrak{T}\), i.e. \(\mathfrak{T}\left(\mathcal{A}\right) \subseteq \mathcal{A}\). By (\ref{eqgyyty736q487eg6347ctsg}), (\ref{eqiw76iebyac7t4jv6}) and (\ref{eqmjhscfnhf6ba7i6w38nrgfjsyg}), for any \(\left(u_{1}(t), u_{2}(t), v_{1}(t), v_{2}(t)\right)\) in \(\mathcal{A}\), we have
\begin{equation*}
\begin{aligned}
\max\{\left\vert F_{1}\left(x\left(t\right)\right)\right\vert, \left\vert F_{2}\left(x\left(t\right)\right)\right\vert\} &\leq M^{*}\left(12 e^{-2\lambda t}\delta^{2} + 16 e^{-\lambda \tau} \delta^{2}\right),\\
\max \{\left\vert G_{1}\left(x\left(t\right)\right)\right\vert, \left\vert G_{2}\left(x\left(t\right)\right)\right\vert\} &\leq M^{*}\left(12 e^{-2\lambda \left(\tau - t\right)}\delta^{2} + 16 e^{-\lambda \tau} \delta^{2}\right).
\end{aligned}
\end{equation*}
Let \(M = 32 M^{*}\lambda^{-1}\). For \(i= 1, 2\), we have
\begin{equation*}
\begin{gathered}
\left\vert \overline{u}_{i}\left(t\right) - e^{-\lambda t} u_{i0}\right\vert \leq \int_{0}^{t} e^{\lambda(s-t)} \left\vert F_{i}\left(x\left(s\right)\right)\right\vert ds
\leq 16 M^{*} \delta^{2}\int_{0}^{t} e^{\lambda(s-t)} \left(e^{-2\lambda s} + e^{-\lambda \tau}\right) ds
\leq M e^{-\lambda t}\delta^{2},\\
\begin{aligned}
\lvert \overline{v}_{i}\left(t\right) - e^{-\lambda\left(\tau-t\right)} v_{i\tau}\rvert
\leq & \int_{t}^{\tau} e^{\lambda(t-s)} \lvert G_{i}\left(x\left(s\right)\right)\rvert ds
\leq 16 M^{*}\delta^{2}\int_{t}^{\tau} e^{\lambda(t-s)}\left(e^{-2\lambda \left(\tau -s\right)} + e^{-\lambda \tau}\right)\\
\leq & M e^{-\lambda\left(\tau-t\right)}\delta^{2}.
\end{aligned}
\end{gathered}
\end{equation*}

Choose \(\delta\) sufficiently small such that \(M\delta < 1\). Taking into account that \(\lvert u_{10}\rvert\), \(\lvert u_{20}\rvert\), \(\lvert v_{1\tau}\rvert\), \(\lvert v_{2\tau}\rvert\) are all bounded by \(\delta\), we have \(\max\{\left\vert\overline{u}_{1}\left(t\right)\right\vert, \left\vert\overline{u}_{2}\left(t\right)\right\vert\} \leq 2e^{-\lambda t}\delta\) and \(\max\{\left\vert\overline{v}_{1}\left(t\right)\right\vert, \left\vert\overline{v}_{2}\left(t\right)\right\vert\} \leq 2e^{-\lambda \left(\tau - t\right)}\delta\). Thus, \(\left(\overline{u}_{1}\left(t\right), \overline{u}_{2}\left(t\right), \overline{v}_{1}\left(t\right), \overline{v}_{2}\left(t\right)\right) \in \mathcal{A}\), as desired.

Meanwhile, we have shown that the image of any element of \(\mathcal{A}\) under \(\mathfrak{T}\) can be written in the form (\ref{eq7bwwkauser6518092hz9ub}) such that the corresponding \(\xi_{1}\), \(\xi_{2}\), \(\zeta_{1}\) and \(\zeta_{2}\) satisfy the estimates given in the statement of the lemma. However, since \(\left(u^{\left(0\right)}, v^{\left(0\right)}\right) = \left(0,0\right)\in \mathcal{A}\), it follows from Remark \ref{Rem9n9iqnubxfybxbxkiqub} that the same holds for the solution \(\left(u\left(t\right), v\left(t\right)\right)\) that satisfies boundary condition (\ref{eq800020}). This ends the proof of Lemma \ref{flowlemmaresonant}.

\section{Proofs of Lemmas \ref{State920jj3nu3uhwu} and \ref{Stat4u74i4jm1ml2mo3p}}\label{Appendix989897g86g7656u}

The main part of this appendix is the proof of the following lemma:
\begin{mylem}\label{Derivativeflowlemma}
Let (\ref{eq68960}) be the local map of system (\ref{eq23000}) and suppose \(\left(u_{10}, v_{10}\right)\in \mathcal{D}_{2}\). Write \(x :=\left(u_{10}, v_{10}\right)\). We have
\begin{equation}\label{eq69000}
\begin{aligned}
\frac{\partial \eta_{1}}{\partial u_{1}}\left(x\right) &= \left(1+\gamma\right) e^{-\lambda_{1}\tau}\left[1 + O\left(\delta\right)\right],\quad
&&\frac{\partial \eta_{1}}{\partial v_{1}}\left(x\right) = \gamma \frac{u_{10}}{v_{10}}\cdot e^{-\lambda_{1}\tau}\left[1 + O\left(\delta\right)\right],\\
\frac{\partial \eta_{2}}{\partial u_{1}}\left(x\right) &= -\gamma \frac{v_{10}}{u_{10}}\cdot e^{\lambda_{1}\tau}\left[1 + O\left(\delta\right)\right],\quad
&&\frac{\partial \eta_{2}}{\partial v_{1}}\left(x\right) = \left(1-\gamma\right) e^{\lambda_{1}\tau}\left[1 + O\left(\delta\right)\right].
\end{aligned}
\end{equation}
\end{mylem}

It is straightforward to derive Lemmas \ref{State920jj3nu3uhwu} and \ref{Stat4u74i4jm1ml2mo3p} from Lemma \ref{Derivativeflowlemma} (see \cite{BakraniPhDthesis}). Indeed, by this lemma, an estimate for \(\frac{\partial \left(\overline{u}_{10}, \overline{v}_{10}\right)}{\partial \left(u_{10}, v_{10}\right)}\) can be obtained. One can use this estimate and the implicit relations between \(\left(u_{10}, v_{10}\right)\), \(\left(w, z\right)\), \(\left(\overline{u}_{10}, \overline{v}_{10}\right)\) and \(\left(\overline{w}, \overline{z}\right)\) to prove Lemmas \ref{State920jj3nu3uhwu} and \ref{Stat4u74i4jm1ml2mo3p}. We refer the reader to \cite{BakraniPhDthesis} for further details.

\begin{proof}[Proof of Lemma \ref{Derivativeflowlemma}]
Let (\ref{eq90ioi48iskju389874jksh}) be the solution of system (\ref{eq23000}) that satisfies boundary conditions (\ref{eq800020}), where \(u_{20} = v_{2\tau} = \delta\). When the point \(\left(u_{10}, \delta, v_{10}, v_{20}\right)\) on \(\Pi^{s}\) reaches the cross-section \(\Pi^{u}\) at \(\left(u_{1\tau}, u_{2\tau}, v_{1\tau}, \delta\right)\), the corresponding flight time \(\tau\) is uniquely determined by \(u_{10}\) and \(v_{10}\), i.e. \(\tau = \tau\left(u_{10}, v_{10}\right)\), for some function \(\tau\). Thus, by (\ref{eqin98nwkndoib112}), we have
\begin{gather}
v_{10} = v_{1}^{*}\left(0, \tau\left(u_{10}, v_{10}\right), u_{10}, \delta, \eta_{2}\left(u_{10},v_{10}\right), \delta\right),\label{eq900040}\\
v_{20} = v_{2}^{*}\left(0, \tau\left(u_{10}, v_{10}\right), u_{10}, \delta, \eta_{2}\left(u_{10},v_{10}\right), \delta\right)\label{eq900045},\\
\eta_{1}\left(u_{10},v_{10}\right) = u_{1}^{*}\left(\tau\left(u_{10},v_{10}\right), \tau\left(u_{10},v_{10}\right), u_{10}, \delta, \eta_{2}\left(u_{10},v_{10}\right), \delta\right).\label{eq900080}
\end{gather}
Recall that, by Lemma \ref{Cor518b8b89or48v10010}, \(v_{20}\) is a function of \(u_{10}\) and \(v_{10}\) which we denote it by \(\kappa\left(u_{10}, v_{10}\right)\).
Both sides of (\ref{eq900040}), (\ref{eq900045}) and (\ref{eq900080}) are functions of \(u_{10}\) and \(v_{10}\). Differentiating these three relations with respect to \(u_{10}\) and \(v_{10}\) gives the following identities
\begin{gather}
0 = \frac{\partial v^{*}_{1}}{\partial\tau}\bigg\vert_{t=0}\cdot\frac{\partial \tau}{\partial u_{10}} + \frac{\partial v^{*}_{1}}{\partial u_{10}}\bigg\vert_{t=0} + \frac{\partial v^{*}_{1}}{\partial v_{1\tau}}\bigg\vert_{t=0}\cdot\frac{\partial \eta_{2}}{\partial u_{10}},\label{eq900050}\\
1 = \frac{\partial v^{*}_{1}}{\partial\tau}\bigg\vert_{t=0}\cdot\frac{\partial \tau}{\partial v_{10}} + \frac{\partial v^{*}_{1}}{\partial v_{1\tau}}\bigg\vert_{t=0}\cdot\frac{\partial \eta_{2}}{\partial v_{10}},\label{eq900060}\\
\frac{\partial \kappa}{\partial u_{10}} = \frac{\partial v^{*}_{2}}{\partial\tau}\bigg\vert_{t=0}\cdot\frac{\partial \tau}{\partial u_{10}} + \frac{\partial v^{*}_{2}}{\partial u_{10}}\bigg\vert_{t=0} + \frac{\partial v^{*}_{2}}{\partial v_{1\tau}}\bigg\vert_{t=0}\cdot\frac{\partial \eta_{2}}{\partial u_{10}},\label{eq7u48ui9i49i90oq0q}\\
\frac{\partial \kappa}{\partial v_{10}} = \frac{\partial v^{*}_{2}}{\partial\tau}\bigg\vert_{t=0}\cdot\frac{\partial \tau}{\partial v_{10}} + \frac{\partial v^{*}_{2}}{\partial v_{1\tau}}\bigg\vert_{t=0}\cdot\frac{\partial \eta_{2}}{\partial v_{10}},\label{equijlms093kdmnjks}\\
\frac{\partial \eta_{1}}{\partial u_{10}} = \frac{\partial u_{1}^{*}}{\partial t}\bigg\vert_{t=\tau}\cdot\frac{\partial \tau}{\partial u_{10}} + \frac{\partial u_{1}^{*}}{\partial \tau}\bigg\vert_{t=\tau}\cdot\frac{\partial \tau}{\partial u_{10}} + \frac{\partial u_{1}^{*}}{\partial u_{10}}\bigg\vert_{t=\tau} + \frac{\partial u_{1}^{*}}{\partial v_{1\tau}}\bigg\vert_{t=\tau}\cdot\frac{\partial \eta_{2}}{\partial u_{10}},\label{eq900090}\\
\frac{\partial \eta_{1}}{\partial v_{10}} = \frac{\partial u_{1}^{*}}{\partial t}\bigg\vert_{t=\tau}\cdot\frac{\partial \tau}{\partial v_{10}} + \frac{\partial u_{1}^{*}}{\partial \tau}\bigg\vert_{t=\tau}\cdot\frac{\partial \tau}{\partial v_{10}} + \frac{\partial u_{1}^{*}}{\partial v_{1\tau}}\bigg\vert_{t=\tau}\cdot\frac{\partial \eta_{2}}{\partial v_{10}}.\label{eq900100}
\end{gather}

To obtain the estimates in (\ref{eq69000}): we first estimate the following expressions\\
\begin{inparaenum}[(i)]
\item\label{Item87fi6f5d65d6bwti} \(\frac{\partial u_{1}^{*}}{\partial t}\Big\vert_{t=\tau}\), \hspace{14mm}
\item\label{Item87fi6f5d65d6bwtii} \(\frac{\partial \kappa}{\partial u_{10}},\, \frac{\partial \kappa}{\partial v_{10}}\), \hspace{14mm}
\item\label{Item87fi6f5d65d6bwtiii} \(\frac{\partial u^{*}_{1}}{\partial u_{10}}\Big\vert_{t=\tau},\,
\frac{\partial v^{*}_{1}}{\partial u_{10}}\Big\vert_{t=0},\,
\frac{\partial v^{*}_{2}}{\partial u_{10}}\Big\vert_{t=0}\),\vspace{1mm}\\
\item\label{Item87fi6f5d65d6bwtiv} \(\frac{\partial u^{*}_{1}}{\partial v_{1\tau}}\Big\vert_{t=\tau},\, \frac{\partial v^{*}_{1}}{\partial v_{1\tau}}\Big\vert_{t=0},\,
\frac{\partial v^{*}_{2}}{\partial v_{1\tau}}\Big\vert_{t=0}\),\hspace{20mm}
\item\label{Item87fi6f5d65d6bwtv} \(\frac{\partial u_{1}^{*}}{\partial \tau}\Big\vert_{t=\tau},\,
\frac{\partial v_{1}^{*}}{\partial \tau}\Big\vert_{t=0},\,
\frac{\partial v_{2}^{*}}{\partial \tau}\Big\vert_{t=0}\).\vspace{1mm}\\
\end{inparaenum}

{\it (\ref{Item87fi6f5d65d6bwti}) Estimate for \(\frac{\partial u_{1}^{*}}{\partial t}\Big\vert_{t=\tau}\)}: By (\ref{eq20000}), (\ref{eq24500}) and the first equation of (\ref{eq23000}), we have
\begin{equation*}
\frac{\partial u_{1}^{*}}{\partial t}\Big\vert_{t=\tau} = -\lambda_{1} u_{1\tau} + O\left(u_{1\tau}^{2}\right) + O\left(u_{1\tau}u_{2\tau}\right) + O\left(v_{1\tau}u_{2\tau}\right),
\end{equation*}
and by virtue of (\ref{eqi903lokdjkflk}), (\ref{eq09o02oknpqlamncccfg}) and (\ref{equ9sjhe9opqk39kxnj}), for \(\left(u_{10}, v_{10}\right)\in \mathcal{D}_{1}\cup \mathcal{D}_{2}\), we have
\begin{equation}\label{eq920llsjcijndyudnnd}
\frac{\partial u_{1}^{*}}{\partial t}\Big\vert_{t=\tau} = -\lambda_{1} e^{-\lambda_{1}\tau}u_{10}\left[1 + O\left(\delta\right)\right].
\end{equation}

{\it (\ref{Item87fi6f5d65d6bwtii}) Estimates for \(\frac{\partial \kappa}{\partial u_{10}}\) and \(\frac{\partial \kappa}{\partial v_{10}}\)}: Following Lemma \ref{Cor518b8b89or48v10010}, \(\kappa\) is a \(\mathcal{C}^q\)-smooth (\(q\geq 2\)) function of \((u_{10}, v_{10})\) which is defined on an open neighborhood of \(M^{s}\in \Pi^{s}\). Since its restriction to \(\mathcal{D}_{1} \cup \mathcal{D}_{2}\) is of the form (\ref{eq09uis4kkn0cusbgdj}), we have
\begin{equation*}
\kappa\left(0, 0\right) = \frac{\partial \kappa}{\partial u_{10}}\left(0, 0\right) = \frac{\partial \kappa}{\partial v_{10}}\left(0, 0\right) = \frac{\partial^{2} \kappa}{\partial u_{10}^{2}}\left(0, 0\right) = \frac{\partial^{2} \kappa}{\partial v_{10}^{2}}\left(0, 0\right) = 0, \quad \frac{\partial^{2} \kappa}{\partial u_{10} v_{10}}\left(0, 0\right) = \frac{\gamma}{\delta}.
\end{equation*}
Since \(v_{10} = O\left(u_{10}\right)\) and \(u_{10} = O\left(v_{10}\right)\), by Taylor theorem, for \((u_{10}, v_{10})\in \mathcal{D}_{2}\), we derive
\begin{equation*}
\frac{\partial \kappa}{\partial u_{10}}\left(u_{10}, v_{10}\right) = \frac{\gamma}{\delta}v_{10}\left[1 + o\left(1\right)\right],\qquad
\frac{\partial \kappa}{\partial v_{10}}\left(u_{10}, v_{10}\right) = \frac{\gamma}{\delta}u_{10}\left[1 + o\left(1\right)\right].
\end{equation*}

In order to get estimates for \(\frac{\partial u^{*}_{1}}{\partial \theta}\), \(\frac{\partial v^{*}_{1}}{\partial \theta}\) and \(\frac{\partial v^{*}_{2}}{\partial \theta}\), where \(\theta = u_{10}\), \(v_{1\tau}\) and \(\tau\), we solve some boundary value problems. Let (\ref{eq90ioi48iskju389874jksh}) be the solution of system (\ref{eq23000}) which satisfies boundary conditions (\ref{eq800020}). By writing system (\ref{eq23000}) in the form (\ref{eqbwiebro4ibyw43yyyiw0290}), i.e.
\begin{equation}\label{eq890lnn3jknnshjnbs}
\begin{gathered}
F_{1}\left(u,v\right) = f_{11}\left(u_{1}, v\right) u_{1} + f_{12}\left(u_{1}, u_{2}, v\right) u_{2},\qquad
F_{2}\left(u,v\right) = f_{22}\left(u_{1}, u_{2}, v\right) u_{2},\\
G_{1}\left(u,v\right) = g_{11}\left(u, v_{1}\right) v_{1} + g_{12}\left(u, v_{1}, v_{2}\right) v_{2},\qquad
G_{2}\left(u,v\right) = g_{22}\left(u, v_{1}, v_{2}\right) v_{2},
\end{gathered}
\end{equation}
where \(f_{ij}\) and \(g_{ij}\) satisfy (\ref{eq20000}) and (\ref{eq24500}), we have
\begin{equation}\label{eq900540}
\begin{aligned}
\dot{u}_{k}^{*} &= -\lambda_{k}u_{k}^{*} + F_{k}(u_{1}^{*}, u_{2}^{*}, v_{1}^{*}, v_{2}^{*}),\\
\dot{v}_{k}^{*} &= +\lambda_{k}v_{k}^{*} + G_{k}(u_{1}^{*}, u_{2}^{*}, v_{1}^{*}, v_{2}^{*}),\qquad (k=1,2).
\end{aligned}
\end{equation}
Differentiating (\ref{eq900540}) with respect to \(\theta\), where \(\theta = u_{10}\), \(v_{1\tau}\) and \(\tau\), gives
\begin{equation}\label{eq900560}
\left(\begin{array}{c}
\dot{U}\\
\dot{V}
\end{array}\right) = \textbf{diagonal}\left(-\lambda_{1}, -\lambda_{2}, \lambda_{1}, \lambda_{2}\right)\cdot\left(\begin{array}{c}
U\\
V
\end{array}\right) + \mathbf{M}(t)\cdot\left(\begin{array}{c}
U\\
V
\end{array}\right),
\end{equation}
where \(U = \left(U_{1}, U_{2}\right)\), \(V = \left(V_{1}, V_{2}\right)\), \(\mathbf{M}(t) = \frac{\partial \left(F_{1}, F_{2}, G_{1}, G_{2}\right)}{\partial \left(u_{1}, u_{2}, v_{1}, v_{2}\right)}\Big\vert_{\left(u^{*}, v^{*}\right)}\) and, for \(\left(i=1, 2\right)\),
\begin{equation}\label{eq900570}
U_{i}(t) = \frac{\partial u^{*}_{i}(t, \tau, u_{10}, u_{20}, v_{1\tau}, v_{2\tau})}{\partial \theta},\qquad
V_{i}(t) = \frac{\partial v^{*}_{i}(t, \tau, u_{10}, u_{20}, v_{1\tau}, v_{2\tau})}{\partial \theta}.
\end{equation}
The solution \(\left(U\left(t\right), V\left(t\right)\right)\) of system (\ref{eq900560}) that satisfies the boundary conditions
\begin{equation}\label{equisnmlkdnnfk0910}
U_{1}(0) = U_{10},\quad U_{2}(0) = U_{20},\quad V_{1}(\tau) = V_{1\tau},\quad V_{2}(\tau) = V_{2\tau}
\end{equation}
is in fact the fixed point of the integral operator
\begin{equation}\label{eqijomk7890ojijgpolkaanbzm891}
\mathfrak{T}: \Big( U_{1}\left(t\right), U_{2}\left(t\right), V_{1}\left(t\right), V_{2}\left(t\right)\Big) \mapsto \left(\overline{U}_{1}\left(t\right), \overline{U}_{2}\left(t\right), \overline{V}_{1}\left(t\right), \overline{V}_{2}\left(t\right)\right),
\end{equation}
such that
\begin{equation*}
\overline{U}_{i}\left(t\right) = e^{-\lambda_{i}t} U_{i0} + \int_{0}^{t} e^{\lambda_{i}(s-t)} P_{i}\left(s\right) ds,\qquad
\overline{V}_{i}\left(t\right) = e^{-\lambda_{i}(\tau -t)} V_{i\tau} + \int_{t}^{\tau} e^{\lambda_{i}(t-s)} Q_{i}\left(s\right) ds,
\end{equation*}
where
\begin{equation*}
\begin{aligned}
P_{i}\left(t\right) =& {F_{i}}_{u_{1}}\left(x^{*}\left(t\right)\right)\cdot U_{1}\left(t\right) + {F_{i}}_{u_{2}}\left(x^{*}\left(t\right)\right)\cdot U_{2}\left(t\right) + {F_{i}}_{v_{1}}\left(x^{*}\left(t\right)\right)\cdot V_{1}\left(t\right) + {F_{i}}_{v_{2}}\left(x^{*}\left(t\right)\right)\cdot V_{2}\left(t\right),\\
Q_{i}\left(t\right) =& {G_{i}}_{u_{1}}\left(x^{*}\left(t\right)\right)\cdot U_{1}\left(t\right) + {G_{i}}_{u_{2}}\left(x^{*}\left(t\right)\right)\cdot U_{2}\left(t\right) + {G_{i}}_{v_{1}}\left(x^{*}\left(t\right)\right)\cdot V_{1}\left(t\right) + {G_{i}}_{v_{2}}\left(x^{*}\left(t\right)\right)\cdot V_{2}\left(t\right),
\end{aligned}
\end{equation*}
for \(i=1, 2\) and \(x^{*}\left(t\right) = \left(u^{*}\left(t\right), v^{*}\left(t\right)\right)\) (see \cite{Dimabook}). Moreover, this integral operator is a contraction and its fixed point \(\left(U\left(t\right), V\left(t\right)\right)\) is the limit of the sequence
\begin{equation*}
\Big\{\left(U^{\left(n\right)}(t), V^{\left(n\right)}(t)\right) = \left(U_{1}^{\left(n\right)}(t), U_{2}^{\left(n\right)}(t), V_{1}^{\left(n\right)}(t), V_{2}^{\left(n\right)}(t)\right)\Big\}_{n=0}^{n=\infty},
\end{equation*}
where \(\left(U^{\left(0\right)}, V^{\left(0\right)}\right) = \left(0,0\right)\) and \(\left(U^{\left(n+1\right)}(t), V^{\left(n+1\right)}(t)\right) = \mathfrak{T} \left(U^{\left(n\right)}\left(t\right), V^{\left(n\right)}\left(t\right)\right)\), \(\forall n\geq 0\).

{\it (\ref{Item87fi6f5d65d6bwtiii}) Estimates for \(\frac{\partial u^{*}_{1}}{\partial u_{10}}\Big\vert_{t=\tau}, \frac{\partial v^{*}_{1}}{\partial u_{10}}\Big\vert_{t=0}\) and \(\frac{\partial v^{*}_{2}}{\partial u_{10}}\Big\vert_{t=0}\)}: Let \(\left(U_{1}, U_{2}, V_{1}, V_{2}\right)\) be the solution of system (\ref{eq900560}), i.e. the fixed point of (\ref{eqijomk7890ojijgpolkaanbzm891}), where
\begin{equation*}
U_{i}(t) = \frac{\partial u^{*}_{i}(t, \tau, u_{10}, u_{20}, v_{1\tau}, v_{2\tau})}{\partial u_{10}},\qquad
V_{i}(t) = \frac{\partial v^{*}_{i}(t, \tau, u_{10}, u_{20}, v_{1\tau}, v_{2\tau})}{\partial u_{10}},\quad \left(i=1, 2\right).
\end{equation*}
Taking into account that (\ref{eq800025}) holds for the solution \(\left(u^{*}, v^{*}\right)\) of system (\ref{eq23000}), we have
\begin{equation}\label{eqiomsnmkc89enhheu}
U_{1}(0) = U_{10} = 1,\quad U_{2}(0) = U_{20} = 0,\quad V_{1}(\tau) = V_{1\tau} = 0,\quad V_{2}(\tau) = V_{2\tau} = 0.
\end{equation}
We claim that the solution \(\left(U, V\right)\) that satisfies (\ref{eqiomsnmkc89enhheu}) is of the form
\begin{equation}\label{eq900840}
\begin{array}{ll}
U_{1}(t) = e^{-\lambda_{1}t} \left[1 + O\left(\delta \right)\right],\quad &
U_{2}(t) = e^{-\lambda_{2}t} O\left(\delta \right),\\
V_{1}(t) = e^{-\lambda_{1}\left(\tau-t\right)} O\left(\delta \right),\quad &
V_{2}(t) = e^{-\lambda_{2}\left(\tau-t\right)} O\left(\delta\right).
\end{array}
\end{equation}
To prove the claim, let us first show that the set
\begin{equation*}
\begin{aligned}
\mathcal{A} = \Big\{\big(U_{1}\left(t\right), U_{2}\left(t\right), V_{1}\left(t\right), V_{2}\left(t\right)\big):&\quad \lvert U_{1}(t)\rvert \leq 2 e^{-\lambda_{1}t},\quad \lvert U_{2}(t)\rvert \leq e^{-\lambda_{2}t}, \\
&\quad \lvert V_{1}(t)\rvert \leq e^{-\lambda_{1}\left(\tau-t\right)},\quad \lvert V_{2}(t)\rvert \leq e^{-\lambda_{2}\left(\tau-t\right)}\Big\},
\end{aligned}
\end{equation*}
where \(\left(U_{1}(t), U_{2}(t), V_{1}(t), V_{2}(t)\right)\) is any continuous function defined on \(t \in \left[0, \tau\right]\), is invariant with respect to integral operator (\ref{eqijomk7890ojijgpolkaanbzm891}). Note that since \(f_{ij}\) and \(g_{ij}\) in (\ref{eq23000}) are \(\mathcal{C}^{q-1}\) (\(q\geq 2\)) and fulfill (\ref{eq20000}) and (\ref{eq24500}), the first derivatives of \(F_{i}\) and \(G_{i}\) can be written as
\begin{equation*}
\begin{gathered}
{F_{1}}_{u_{1}}\left(x\right) = \mathtt{f}^{1}_{11}\left(x\right) u_{1} + \mathtt{f}^{1}_{12}\left(x\right) u_{2},\qquad
{F_{1}}_{u_{2}}\left(x\right) = \mathtt{f}^{2}_{11}\left(x\right) v_{1} + \mathtt{f}^{2}_{12}\left(x\right) v_{2},\\
{F_{1}}_{v_{1}}\left(x\right) = \mathtt{f}^{3}_{11}\left(x\right) u_{1} + \mathtt{f}^{3}_{12}\left(x\right) u_{2}, \qquad 
{F_{1}}_{v_{2}}\left(x\right) = \mathtt{f}^{4}_{11}\left(x\right) u_{1} + \mathtt{f}^{4}_{12}\left(x\right) u_{2},\\
{F_{2}}_{u_{1}}\left(x\right) = \mathtt{f}^{1}_{21}\left(x\right) u_{2},\qquad {F_{2}}_{u_{2}}\left(x\right) = \mathtt{f}^{2}_{21}\left(x\right) u_{1} + \mathtt{f}^{2}_{22}\left(x\right) u_{2},\\
{F_{2}}_{v_{1}}\left(x\right) = \mathtt{f}^{3}_{21}\left(x\right) u_{2}\qquad
{F_{2}}_{v_{2}}\left(x\right) = \mathtt{f}^{4}_{21}\left(x\right) u_{2},\qquad
{G_{1}}_{u_{1}}\left(x\right) = \mathtt{g}^{1}_{11}\left(x\right) v_{1} + \mathtt{g}^{1}_{12}\left(x\right) v_{2},\\
{G_{1}}_{u_{2}}\left(x\right) = \mathtt{g}^{2}_{11}\left(x\right) v_{1} + \mathtt{g}^{2}_{12}\left(x\right) v_{2},\qquad
{G_{1}}_{v_{1}}\left(x\right) = \mathtt{g}^{3}_{11}\left(x\right) v_{1} + \mathtt{g}^{3}_{12}\left(x\right) v_{2},\\
{G_{1}}_{v_{2}}\left(x\right) = \mathtt{g}^{4}_{11}\left(x\right) u_{1} + \mathtt{g}^{4}_{12}\left(x\right) u_{2},\qquad
{G_{2}}_{u_{1}}\left(x\right) = \mathtt{g}^{1}_{21}\left(x\right) v_{2},\\
{G_{2}}_{u_{2}}\left(x\right) = \mathtt{g}^{2}_{21}\left(x\right) v_{2},\qquad {G_{2}}_{v_{1}}\left(x\right) = \mathtt{g}^{3}_{21}\left(x\right) v_{2},\qquad
{G_{2}}_{v_{2}}\left(x\right) = \mathtt{g}^{4}_{21}\left(x\right) v_{1} + \mathtt{g}^{4}_{22}\left(x\right) v_{2},
\end{gathered}
\end{equation*}
where \(x=\left(u, v\right)\), and \(\mathtt{f}^{k}_{ij}\) and \(\mathtt{g}^{k}_{ij}\) are some continuous functions. Consider the constant \(M\) given by Lemma \ref{flowlemma0}. Recall that \(\delta\) is sufficiently small such that \(M\delta < 1\). Let \(M^{\dagger} = \max\lbrace 3, 3M\rbrace\). Then, for the solution \(\left(u^{*}\left(t\right), v^{*}\left(t\right)\right)\) of system (\ref{eq23000}), we have
\begin{equation*}
\lvert u^{*}_{i}\left(t\right)\rvert \leq M^{\dagger} e^{-\lambda_{i}t} \delta,\qquad 
\lvert v^{*}_{i}\left(t\right)\rvert \leq M^{\dagger} e^{-\lambda_{i}\left(\tau-t\right)}\delta,\qquad (i=1,2).
\end{equation*}
Let \(\Omega\) be a small compact neighborhood of the equilibrium \(O\) of system (\ref{eq23000}). Define \(M^{*} := \sup_{\left(u, v\right)\in \Omega} \big\{ \lvert \mathtt{f}^{k}_{ij}\left(u, v\right)\rvert, \lvert \mathtt{g}^{k}_{ij}\left(u, v\right)\rvert\big\}\), and \(M^{\ddagger} := M^{*}M^{\dagger}\). We have
\begin{equation*}
\begin{gathered}
\big\vert {F_{1}}_{u_{1}}\left(u^{*},v^{*}\right)\big\vert, \big\vert {F_{1}}_{v_{1}}\left(u^{*},v^{*}\right)\big\vert, \big\vert {F_{1}}_{v_{2}}\left(u^{*},v^{*}\right)\big\vert, \big\vert {F_{2}}_{u_{2}}\left(u^{*},v^{*}\right)\big\vert, \big\vert {G_{1}}_{v_{2}}\left(u^{*},v^{*}\right)\big\vert \leq M^{\ddagger} e^{-\lambda_{1}t}\delta\\
\big\vert {F_{2}}_{u_{1}}\left(u^{*},v^{*}\right)\big\vert, \big\vert {F_{2}}_{v_{1}}\left(u^{*},v^{*}\right)\big\vert, \big\vert {F_{2}}_{v_{2}}\left(u^{*},v^{*}\right)\big\vert \leq M^{\ddagger} e^{-\lambda_{2}t}\delta\\
\big\vert {F_{1}}_{u_{2}}\left(u^{*},v^{*}\right)\big\vert, \big\vert {G_{1}}_{u_{1}}\left(u^{*},v^{*}\right)\big\vert, \big\vert {G_{1}}_{u_{2}}\left(u^{*},v^{*}\right)\big\vert, \big\vert {G_{1}}_{v_{1}}\left(u^{*},v^{*}\right)\big\vert, \big\vert {G_{2}}_{v_{2}}\left(u^{*},v^{*}\right)\big\vert \leq M^{\ddagger} e^{-\lambda_{1}\left(\tau - t\right)}\delta\\
\big\vert {G_{2}}_{u_{1}}\left(u^{*},v^{*}\right)\big\vert, \big\vert {G_{2}}_{u_{2}}\left(u^{*},v^{*}\right)\big\vert, \big\vert {G_{2}}_{v_{1}}\left(u^{*},v^{*}\right)\big\vert \leq M^{\ddagger} e^{-\lambda_{2}\left(\tau - t\right)}\delta.
\end{gathered}
\end{equation*}
This implies
\begin{equation*}
\begin{array}{cc}
\lvert P_{1}\left(t\right)\rvert \leq M^{\ddagger} \delta\left[3 e^{-2\lambda_{1}t} + 2 e^{-\lambda_{1}\tau}\right],
&\lvert P_{2}\left(t\right)\rvert \leq M^{\ddagger}\delta \left[3 e^{-\left(\lambda_{1} + \lambda_{2}\right)t} + 2 e^{-\lambda_{2}t - \lambda_{1}\left(\tau - t\right)}\right],\\
\lvert Q_{1}\left(t\right)\rvert \leq M^{\ddagger} \delta\left[3 e^{-\lambda_{1}\tau} + 2 e^{-2\lambda_{1}\left(\tau - t\right)}\right],
& \lvert Q_{2}\left(t\right)\rvert \leq M^{\ddagger}\delta \left[3e^{-\lambda_{1}t -\lambda_{2}\left(\tau - t\right)} + 2e^{-\left(\lambda_{1} +\lambda_{2}\right)\left(\tau - t\right)}\right].
\end{array}
\end{equation*}
Define \(M = 6 M^{\ddagger} {\lambda_{1}}^{-1}\). Using the above relations, we have
\begin{equation*}
\begin{gathered}
\Big\vert\overline{U}_{1}(t) - e^{-\lambda_{1}t}\Big\vert \leq \int_{0}^{t} \Big\vert e^{\lambda_{1}(s-t)} P_{1}\left(s\right)\Big\vert ds
\leq 3 M^{\ddagger}\delta \int_{0}^{t} e^{\lambda_{1}(s-t)} \left[e^{-2\lambda_{1}s} + e^{-\lambda_{1}\tau}\right] ds
\leq M e^{-\lambda_{1}t}\delta,\\
%%%%%%%%%%%%%%%%%%%%%%%%%%%
\Big\vert\overline{U}_{2}(t)\Big\vert \leq \int_{0}^{t} \Big\vert e^{\lambda_{2}(s-t)} P_{2}\left(s\right)\Big\vert ds
\leq 3 M^{\ddagger}\delta \int_{0}^{t} e^{\lambda_{2}(s-t)} \left[e^{-\left(\lambda_{1} + \lambda_{2}\right)s} + e^{-\lambda_{2}s - \lambda_{1}\left(\tau - s\right)}\right] ds
\leq M e^{-\lambda_{2}t}\delta,
\end{gathered}
\end{equation*}
%%%%%%%%%%%%%%%%%%%%%%%%%%%
\begin{equation*}
\begin{gathered}
\Big\vert\overline{V}_{1}(t) \Big\vert \leq \int_{t}^{\tau}\Big\vert e^{\lambda_{1}(t-s)} Q_{1}\left(s\right)\Big\vert ds
\leq 3 M^{\ddagger}\delta \int_{t}^{\tau} e^{\lambda_{1}(t-s)} \left[e^{-\lambda_{1}\tau} + e^{-2\lambda_{1}\left(\tau - s\right)}\right] ds \leq M e^{-\lambda_{1}\left(\tau - t\right)}\delta,\\
%%%%%%%%%%%%%%%%%%%%%%%%%%%
\begin{aligned}
\Big\vert\overline{V}_{2}(t)\Big\vert \leq & \int_{t}^{\tau}\Big\vert e^{\lambda_{2}(t-s)} Q_{2}\left(s\right)\Big\vert ds
\leq 3 M^{\ddagger}\delta \int_{t}^{\tau} e^{\lambda_{2}(t-s)} \left[e^{-\lambda_{1}s -\lambda_{2}\left(\tau - s\right)} + e^{-\left(\lambda_{1} +\lambda_{2}\right)\left(\tau - s\right)}\right] ds\\
\leq & M e^{-\lambda_{2}\left(\tau - t\right)}\delta.
\end{aligned}
\end{gathered}
\end{equation*}
Choosing \(\delta\) sufficiently small such that \(M\delta < 1\), the above relation immediately implies \(\left(\overline{U}_{1}\left(t\right), \overline{U}_{2}\left(t\right), \overline{V}_{1}\left(t\right), \overline{V}_{2}\left(t\right)\right)\in \mathcal{A}\), as desired.

Meanwhile, we have shown that the image of any element of \(\mathcal{A}\) under \(\mathfrak{T}\) is of the form (\ref{eq900840}). However, since \(\left(U^{\left(0\right)}, V^{\left(0\right)}\right) \equiv \left(0,0\right)\in \mathcal{A}\), it follows from Remark \ref{Rem9n9iqnubxfybxbxkiqub} that the same holds for the solution \(\left(U\left(t\right), V\left(t\right)\right)\) that satisfies boundary condition (\ref{equisnmlkdnnfk0910}). This gives (\ref{eq900840}) and therefore,
\begin{equation*}
\begin{aligned}
\frac{\partial u^{*}_{1}}{\partial u_{10}}\bigg\vert_{t=\tau} = e^{-\lambda_{1}\tau} \left[1 + O\left(\delta \right)\right],\quad
\frac{\partial v^{*}_{1}}{\partial u_{10}}\bigg\vert_{t=0} = e^{-\lambda_{1}\tau} O\left(\delta \right),\quad
\frac{\partial v^{*}_{2}}{\partial u_{10}}\bigg\vert_{t=0} = e^{-\lambda_{2}\tau} O\left(\delta \right).
\end{aligned}
\end{equation*}

{\it (\ref{Item87fi6f5d65d6bwtiv}) Estimates for \(\frac{\partial u^{*}_{1}}{\partial v_{1\tau}}\Big\vert_{t=\tau}, \frac{\partial v^{*}_{1}}{\partial v_{1\tau}}\Big\vert_{t=0}\) and \(\frac{\partial v^{*}_{2}}{\partial v_{1\tau}}\Big\vert_{t=0}\)}: Let \(\left(U_{1}, U_{2}, V_{1}, V_{2}\right)\) be the solution of system (\ref{eq900560}), where \(U_{i}\) and \(V_{i}\) are as in (\ref{eq900570}) for \(\theta = v_{1\tau}\). With the same method that we derived (\ref{eq900840}), one can prove that when \(\left(u_{10}, v_{10}\right)\in \mathcal{D}_{2}\), the solution \(\left(U, V\right)\) is of the form
\begin{equation*}
\begin{array}{ll}
U_{1}(t) = e^{-\lambda_{1}\left(\tau+t\right)} O\left(\delta \right),\quad &
U_{2}(t) = e^{-\lambda_{2}t} O\left(\delta \right),\\
V_{1}(t) = e^{-\lambda_{1}\left(\tau-t\right)} \left[1 + O\left(\delta \right)\right],\quad &
V_{2}(t) = e^{-\lambda_{2}\left(\tau-t\right)} O\left(\delta\right),
\end{array}
\end{equation*}
(see \cite{BakraniPhDthesis}). Therefore, when \(\left(u_{10}, v_{10}\right)\in \mathcal{D}_{2}\), we have
\begin{equation*}
\begin{aligned}
\frac{\partial u^{*}_{1}}{\partial v_{1\tau}}\Big\vert_{t=\tau} = e^{-2\lambda_{1}\tau} O\left(\delta\right),\quad
\frac{\partial v^{*}_{1}}{\partial v_{1\tau}}\Big\vert_{t=0} = e^{-\lambda_{1}\tau} \left[1 + O\left(\delta \right)\right],\quad
\frac{\partial v^{*}_{2}}{\partial v_{1\tau}}\Big\vert_{t=0} = e^{-\lambda_{2}\tau} O\left(\delta\right).
\end{aligned}
\end{equation*}

{\it (\ref{Item87fi6f5d65d6bwtv}) Estimates for \(\frac{\partial u_{1}^{*}}{\partial \tau}\Big\vert_{t=\tau}\), \(\frac{\partial v_{1}^{*}}{\partial \tau}\Big\vert_{t=0}\) and \(\frac{\partial v_{2}^{*}}{\partial \tau}\Big\vert_{t=0}\)}: Let \(\left(U_{1}, U_{2}, V_{1}, V_{2}\right)\) be the solution of system (\ref{eq900560}), where \(U_{i}\) and \(V_{i}\) are as in (\ref{eq900570}) for \(\theta = \tau\). With the same method that we derived (\ref{eq900840}), one can prove that when \(\left(u_{10}, v_{10}\right)\in \mathcal{D}_{2}\), the solution \(\left(U, V\right)\) is of the form
\begin{equation*}
\begin{array}{ll}
U_{1}(t) = e^{-\lambda_{1}t} O\left(\delta \lvert u_{10}\rvert\right),\quad &
U_{2}(t) = e^{-\lambda_{2}t} O\left(\delta^{2} \right),\\
V_{1}(t) = e^{-\lambda_{1}\left(\tau-t\right)} v_{1\tau}\left[-\lambda_{1} + O\left(\delta \right)\right],\quad &
V_{2}(t) = e^{-\lambda_{2}\left(\tau-t\right)}\delta \left[-\lambda_{2} + O\left(\delta\right)\right],
\end{array}
\end{equation*}
(see \cite{BakraniPhDthesis}). Therefore, when \(\left(u_{10}, v_{10}\right)\in \mathcal{D}_{2}\), we have
\begin{equation*}
\begin{gathered}
\frac{\partial u_{1}^{*}}{\partial \tau}\Big\vert_{t=\tau} = e^{-\lambda_{1}\tau} O\left(\delta \lvert u_{10}\rvert\right),\quad
\frac{\partial v_{1}^{*}}{\partial \tau}\Big\vert_{t=0} = -\lambda_{1} e^{-\lambda_{1}\tau} v_{1\tau}\left[1 + O\left(\delta \right)\right],\\
\frac{\partial v_{2}^{*}}{\partial \tau}\Big\vert_{t=0} = -\lambda_{2} e^{-\lambda_{2}\tau} \delta\left[1 + O\left(\delta \right)\right].
\end{gathered}
\end{equation*}

So far, we have obtained all the estimates that we required. Substituting these estimates into (\ref{eq900060}) and (\ref{equijlms093kdmnjks}) gives
\begin{gather}
1 = -\lambda_{1}e^{-\lambda_{1}\tau} v_{1\tau}\left[1 + O\left(\delta \right)\right]\cdot\frac{\partial \tau}{\partial v_{10}} + e^{-\lambda_{1}\tau} \left[1 + O\left(\delta \right)\right]\cdot\frac{\partial \eta_{2}}{\partial v_{10}}\label{eq67uiwo01893uu3nj}\\
\frac{\gamma u_{10}}{\delta}\left[1 + O\left(\delta \right)\right] = -\lambda_{2} e^{-\lambda_{2}\tau} \delta\left[1 + O\left(\delta \right)\right]\cdot\frac{\partial \tau}{\partial v_{10}} + e^{-\lambda_{2}\tau} O\left(\delta\right)\cdot\frac{\partial \eta_{2}}{\partial v_{10}}.\label{eqyusoll45ljsojhs}
\end{gather}
Relation (\ref{eqyusoll45ljsojhs}) implies
\begin{equation}\label{eq55yu72900ii2po}
\frac{\partial \tau}{\partial v_{10}} = \frac{-1}{\lambda_{2} \delta}\left(e^{\lambda_{2}\tau} \frac{\gamma u_{10}}{\delta} + O\left(\delta\right)\cdot\frac{\partial \eta_{2}}{\partial v_{10}}\right) \left[1 + O\left(\delta \right)\right].
\end{equation}
By substituting this into (\ref{eq67uiwo01893uu3nj}), we have
\begin{equation*}
\frac{\partial \eta_{2}}{\partial v_{10}} = \left(e^{\lambda_{1}\tau} - \frac{\gamma^{2}}{\delta^{2}}\cdot e^{\lambda_{2}\tau} u_{10}v_{1\tau}\right)\left[1 + O\left(\delta \right)\right] \xlongequal{\text{Proposition } \ref{Cor67ujokjsuijwww}} \frac{\partial \eta_{2}}{\partial v_{10}} = \left(1-\gamma\right) e^{\lambda_{1}\tau}\Big[1 + O\left(\delta\right)\Big],
\end{equation*}
as desired in (\ref{eq69000}). By Proposition \ref{Cor67ujokjsuijwww}, substituting this estimate into (\ref{eq55yu72900ii2po}) yields
\begin{equation*}
\frac{\partial \tau}{\partial v_{10}} = -\frac{1}{\lambda_{2}}\cdot\frac{1}{v_{10}} \left[1 + O(v_{1\tau})\right] =  -\frac{1}{\lambda_{2}}\cdot\frac{1}{v_{10}} \left[1 + O(\delta)\right].
\end{equation*}

Similarly, we can estimate derivatives of \(\tau\) and \(\eta_{2}\) with respect to \(u_{10}\). By substituting the obtained estimates into (\ref{eq7u48ui9i49i90oq0q}), we have
\begin{equation}\label{eq67hjuoe912p4}
\frac{\partial \tau}{\partial u_{10}} = \frac{-1}{\lambda_{2} \delta} \left(\frac{\gamma}{\delta} e^{\lambda_{2}\tau}v_{10} + O\left(\delta \right) + O\left(\delta\right)\cdot\frac{\partial \eta_{2}}{\partial u_{10}}\right)\left[1 + O\left(\delta \right)\right].
\end{equation}
Substituting these estimates into (\ref{eq900050}) and simplifying the result using Proposition \ref{Cor67ujokjsuijwww} give the desired estimate for \(\frac{\partial \eta_{2}}{\partial u_{10}}\). In addition, we have
\begin{equation*}
\frac{\partial \tau}{\partial u_{10}} = -\frac{1}{\lambda_{2}}\cdot\frac{1}{u_{10}} \left[1 + O(\delta)\right].
\end{equation*}

It is easily seen that substituting the estimates that we have derived so far into (\ref{eq900090}) and (\ref{eq900100}) gives the desired estimates for \(\frac{\partial \eta_{1}}{\partial u_{10}}\) and \(\frac{\partial \eta_{1}}{\partial v_{10}}\). This ends the proof.
\end{proof}

\begin{myrem}\label{rem90876c5xxz5ktc}
In the case of homoclinic figure-eight, the estimates given by Lemma \ref{Derivativeflowlemma} also hold for the local maps \(T_{1}^{\text{loc}}\) (on \(\mathcal{D}^{1}_{2}\)), \(T_{12}^{\text{loc}}\) (on \(\mathbb{D}^{1}_{2}\)), \(T_{21}^{\text{loc}}\) (on \(\mathbb{D}^{2}_{2}\)) and \(T_{2}^{\text{loc}}\) (on \(\mathcal{D}^{2}_{2}\)). For instance, applying Lemma \ref{Derivativeflowlemma} on the local map \(T^{\text{loc}}\) on \(\mathcal{D}_{2}\) of the system which is derived from system (\ref{eq23000}) by applying the linear change of coordinates \(\left(\tilde{u}_{1}, \tilde{u}_{2}, \tilde{v}_{1}, \tilde{v}_{2}\right) = \left(u_{1}, u_{2}, -v_{1}, -v_{2}\right)\) gives the estimates in Lemma \ref{Derivativeflowlemma} for \(T_{12}^{\text{loc}}\) on \(\mathbb{D}^{1}_{2}\).
\end{myrem}

\section{Invariant manifolds theory for cross-maps}\label{Invariantmanifoldscrossmaps}

In this appendix, we briefly discuss the method of cross-maps which is used in this paper to prove the existence of the invariant manifolds of the Poincar\'e map along the homoclinic orbits. We start with the formal definition of cross-maps:
\begin{mydefn}\label{Defno87iki6jaytxeyt}
Let \(\left(X^{*}, \left\Vert \cdot\right\Vert_{X^{*}}\right)\) and \(\left(Y^{*}, \left\Vert \cdot\right\Vert_{Y^{*}}\right)\) be two Banach spaces, and \(\mathcal{U}\) be a subset of \(X^{*}\times Y^{*}\). Let
\begin{equation*}
\begin{aligned}
T: \mathcal{U} &\rightarrow T\left(\mathcal{U}\right)\\
\left(x, y\right) &\mapsto \left(\overline{x}, \overline{y}\right)
\end{aligned},
\end{equation*}
be a map. We say \(T\) can be written in cross-form if and only if
\begin{equation}\label{eq4000}
\begin{aligned}
\overline{x} &= F\left(x, \overline{y}\right),\\
y &= G\left(x,\overline{y}\right),
\end{aligned}
\end{equation}
holds for some functions \(F\) and \(G\). The map defined by (\ref{eq4000}) (which maps \(\left(x,\overline{y}\right)\) to \(\left(\overline{x}, y\right)\)), is called the cross-map of \(T\) and denoted by \(T^{\times}\).
\end{mydefn}

In general, the composition of two maps which each can be written in cross-form cannot necessarily be written in cross-form. Here we provide a specific setting in which the property of 'being written in cross-form' can transfer to the composition map: let \(\left(\mathcal{B}, \left\Vert \cdot\right\Vert\right)\) be a Banach space, and \(X_{1}\), \(X_{2}\), \(Y_{1}\) and \(Y_{2}\) be convex subsets of \(\mathcal{B}\). Consider the maps \(T_{1}: X_{1}\times Y_{1} \rightarrow X_{2}\times Y_{2}\) and \(T_{2}: X_{2}\times Y_{2} \rightarrow X_{1}\times Y_{1}\) and suppose that both of them can be written in cross-form in the following way:
\begin{equation*}
\left(\overline{x}, \overline{y}\right) = T_{1}\left(x, y\right) \quad \text{ if and only if } \quad  \overline{x} = p_{1}\left(x, \overline{y}\right) \,\text{ and }\, y = q_{1}\left(x, \overline{y}\right),
\end{equation*}
and
\begin{equation*}
\left(\hat{x}, \hat{y}\right) = T_{2}\left(\overline{x}, \overline{y}\right) \quad \text{ if and only if } \quad  \hat{x} = p_{2}\left(\overline{x}, \hat{y}\right) \,\text{ and }\, \overline{y} = q_{2}\left(\overline{x}, \hat{y}\right),
\end{equation*}
where \(p_{1}: X_{1} \times Y_{2} \rightarrow X_{2}\), \(q_{1}: X_{1} \times Y_{2} \rightarrow Y_{1}\), \(p_{2}: X_{2} \times Y_{1} \rightarrow X_{1}\) and \(q_{2}: X_{2} \times Y_{1} \rightarrow Y_{2}\) are some smooth functions. Let
\begin{equation*}
\begin{aligned}
&\max \left\{\left\|\frac{\partial p_{1}}{\partial x}\right\|,\left\|\frac{\partial p_{1}}{\partial \overline{y}}\right\|,\left\|\frac{\partial q_{1}}{\partial x}\right\|,\left\|\frac{\partial q_{1}}{\partial \overline{y}}\right\|\right\} \leq K_{1},\\
&\max \left\{\left\|\frac{\partial p_{2}}{\partial \overline{x}}\right\|,\left\|\frac{\partial p_{2}}{\partial \hat{y}}\right\|,\left\|\frac{\partial q_{2}}{\partial \overline{x}}\right\|,\left\|\frac{\partial q_{2}}{\partial \hat{y}}\right\|\right\} \leq K_{2}
\end{aligned}
\end{equation*}
for some constants \(K_{1}\) and \(K_{2}\).
\begin{mylem}\label{Lem22009o0ub8yvy7v21w}
(\cite{Turaev2014}, Lemma 4) Define \(T:= T_{2}\circ T_{1}: X_{1}\times Y_{1}\rightarrow X_{1}\times Y_{1}\). If \(K_{1}K_{2} < 1\), then
\begin{enumerate}[(i)]
\item the map \(T\) can be written in cross-form, i.e. there exist functions \(p\) and \(q\) such that 
\begin{equation*}
\left(\hat{x}, \hat{y}\right) = T\left(x, y\right) \quad \text{ if and only if } \quad  \hat{x} = p\left(x, \hat{y}\right) \,\text{ and }\, y = q\left(x, \hat{y}\right).
\end{equation*}
Moreover, the functions \(p\) and \(q\) are smooth and defined everywhere on \(X_{1}\times Y_{1}\).
\item Equip \(X_{1}\times Y_{1}\) with the norm \(\left\| \left(x, y\right)\right\|_{*} = \max \{\sqrt{K_{1}} \Vert x\Vert, \sqrt{K_{2}} \Vert y\Vert\}\). We have
\end{enumerate}
\begin{equation*}
\left\|\frac{\partial\left(p, q\right)}{\partial(x, \hat{y})}\right\|_{*} \leq \frac{\sqrt{K_{1} K_{2}}}{1-\sqrt{K_{1} K_{2}}}.
\end{equation*}
\end{mylem}

The next theorem provides a setting in which if a map \(T\) possesses a cross-map \(T^{\times}\) which satisfies certain properties, then it has an invariant manifold that contains \(\omega\)-limit points of every forward orbit of the domain. This theorem becomes powerful when one is looking for the invariant manifolds of a non-smooth map whose cross-map is smooth. This result was first obtained by Afraimovich and Shilnikov \cite{AnnulusPrinciple} for maps defined on an annulus. The following formulation of this result which holds for arbitrary Banach spaces is stated in \cite{Dimabook}.
\begin{mythm}\label{thm3000}
(\cite{Dimabook}, Theorem 4.3) With the setting in Definition \ref{Defno87iki6jaytxeyt}, let \(X\) and \(Y\) be two convex closed subsets of \(X^{*}\) and \(Y^{*}\), respectively, such that \(\mathcal{R} = X\times Y\subset \mathcal{U}\), \(T^{\times}\) is defined on \(\mathcal{R}\) and \(T^{\times}\left(\mathcal{R}\right) \subset \mathcal{R}\). Let \(F\) and \(G\) in (\ref{eq4000}) be \(\mathcal{C}^{1}\)-smooth and satisfy
\begin{equation*}
\sqrt{\sup_{(x, \overline{y})\in X\times Y} \Bigg\lbrace \bigg\Vert \frac{\partial F}{\partial x}\bigg\Vert \cdot \bigg\Vert \frac{\partial G}{\partial \overline{y}}\bigg\Vert\Bigg\rbrace} + \sqrt{\bigg\Vert \frac{\partial F}{\partial \overline{y}}\bigg\Vert_{\circ}\cdot\bigg\Vert \frac{\partial G}{\partial x}\bigg\Vert_{\circ}} < 1
\end{equation*}
and
\begin{equation*}
\bigg\Vert \frac{\partial F}{\partial x}\bigg\Vert_{\circ} + \sqrt{\bigg\Vert \frac{\partial F}{\partial \overline{y}}\bigg\Vert_{\circ}\cdot\bigg\Vert \frac{\partial G}{\partial x}\bigg\Vert_{\circ}} < 1,
\end{equation*}
where \(\left\Vert \varphi\left(x, \overline{y}\right)\right\Vert_{\circ} = \sup_{(x, \overline{y})\in X\times Y} \left\Vert \varphi\left(x, \overline{y}\right)\right\Vert\) for any vector-valued or matrix-valued function \(\varphi\). Then
\begin{enumerate}[(i)]
\item the map \(T\) has a \(\mathcal{C}^{1}\)-smooth invariant manifold \(M^{*} = \lbrace \left(x,y\right)\in \mathcal{R}: x = h^{*}\left(y\right)\rbrace\), where \(h^{*}: Y\rightarrow X\) is a Lipschitz function with the Lipschitz constant
\begin{equation*}
\mathcal{L}=\sqrt{\left\Vert\frac{\partial F}{\partial \overline{y}}\right\Vert_{\circ}\left(\left\Vert\frac{\partial G}{\partial x}\right\Vert_{\circ}\right)^{-1}}.
\end{equation*}
\item for any \(\mathtt{x} = \left(x, y\right)\in \mathcal{R}\) and any arbitrary \(\epsilon > 0\), there exists an integer \(\mathcal{N}_{\epsilon}^{\mathtt{x}} \in \mathbb{N}\) such that for any \(n > \mathcal{N}_{\epsilon}^{\mathtt{x}}\) if \(\lbrace T^{i}\left(\mathtt{x}\right)\rbrace_{i=0}^{i=n} \subset \mathcal{R}\), then \(\mathrm{dist}\left(T^{n}\left(\mathtt{x}\right), M^{*}\right) < \epsilon\). In particular, \(M^{*}\) contains the \(\omega\)-limit set of any point of \(\mathcal{R}\) whose forward orbit lies entirely in \(\mathcal{R}\).
\item\label{Item53iniubyv7t65czcr} if \(\mathcal{R}\) is bounded, then the integer \(\mathcal{N}^{\epsilon}_{\mathtt{x}}\) given above can be chosen independent of \(\mathtt{x}\), i.e. for any arbitrary \(\epsilon > 0\), there exists an integer \(\mathcal{N}_{\epsilon} \in \mathbb{N}\) such that for any \(n > \mathcal{N}_{\epsilon}\) and any \(\mathtt{x}\in \mathcal{R}\) if \(\lbrace T^{i}\left(\mathtt{x}\right)\rbrace_{i=0}^{i=n} \subset \mathcal{R}\), then \(\mathrm{dist}\left(T^{n}\left(\mathtt{x}\right), M^{*}\right) < \epsilon\).
\item let \(M\) be a \(\mathcal{L}\)-surface (i.e. \(M\) is the graph of some \(\mathcal{L}\)-Lipschitz function \(h: Y\rightarrow X\)). Then \(T\left(M\right)\vert_{X\times Y}\) is a \(\mathcal{L}\)-surface as well. Moreover, the sequence \(\lbrace T^{n}\left(M\right)\vert_{X\times Y}\rbrace\) converges to \(M^{*}\).
\end{enumerate}
\end{mythm}
\begin{proof}
See \cite{Dimabook}, Theorem 4.3 as well as Theorem 4.2 and its proof.
\end{proof}

\begin{myprop}\label{Prop7jnhiojhui22bh}
With the setting of Theorem \ref{thm3000}, if \(\mathcal{R}\) is bounded, \(T^{-1}\) exists and the backward orbit of a point \(\mathtt{x}\in\mathcal{R}\) lies entirely in \(\mathcal{R}\) then \(\mathtt{x}\in M^{*}\).
\end{myprop}
\begin{proof}
The proof is by contradiction. Assume \(\mathtt{x}\notin M^{*}\). This implies \(\text{dist}\left(\mathtt{x}, M^{*}\right) > 0\). Choose an \(0 < \epsilon < \text{dist}\left(\mathtt{x}, M^{*}\right)\) and consider \(\mathcal{N}_{\epsilon}\) given by Theorem \ref{thm3000}. We have
\(\text{dist}\left(M^{*}, T^{\mathcal{N}_{\epsilon}}\left(T^{-\mathcal{N}_{\epsilon}}\left(\mathtt{x}\right)\right)\right) < \epsilon\) and thereby
\(\text{dist}\left(M^{*}, \mathtt{x}\right) < \epsilon\), which is a contradiction.
\end{proof}

\bibliography{Myreferences}{}
\bibliographystyle{alpha}

\end{document}